\documentclass{amsart} 
\usepackage{amsfonts}
\usepackage{amsopn}
\usepackage{epsfig}
\usepackage{amsmath}
\usepackage{epic,eepic,latexsym}
\usepackage{color}

\newtheorem{theorem}{Theorem}[section]
\newtheorem{proposition}[theorem]{Proposition}
\newtheorem{lemma}[theorem]{Lemma}

\newtheorem{conjecture}[theorem]{Conjecture}

\theoremstyle{definition}

\newtheorem{question}[theorem]{Question}

\newtheorem{example}[theorem]{Example}

\newcommand{\PP}{\mathbb{P}}
\newcommand{\proj}{\mathbb{P}}
\newcommand{\FF}{\mathbb{F}}

\newcommand{\ZZ}{\mathbb{Z}}
\newcommand{\CC}{\mathbb{C}}

\newcommand{\GG}{\mathbb{G}}

\newcommand{\Z}{\mathbb{Z}}

\newcommand{\Flag}{Fl(k_1, \dots, k_r ; n)}
\newcommand{\Fl}{Fl(k_1, k_2 ; n)}

\newcommand{\Kgnb}[1]{\overline{\mathcal{M}}_{#1}}

\newcommand{\remind}[1]{{}}
\renewcommand{\remind}[1]{{\bf [#1]}} 
\newcommand{\lremind}[1]{{}}
\newcommand{\kbb}{L}
\newcommand{\al}{\alpha}

\newcommand{\be}{\beta}
\newcommand{\ga}{\gamma}

\newcommand{\si}{\sigma}

\newcommand{\codim}{\operatorname{codim}}
\newcommand{\oh}{\mathcal{O}}
\newcommand{\vv}{\vec{v}}
\newcommand{\C}{\mathbb{C}}
\newcommand{\oO}{\overline{\Omega}}
\newcommand{\dis}{\operatorname{dis}}

\title[Geometric Schubert calculus]{Geometric positivity in the cohomology of
  homogeneous spaces and generalized Schubert calculus}
\author[Coskun]{Izzet Coskun}\address{Department of Mathematics,
M.I.T., Cambridge, MA 02139}
\email{coskun@math.mit.edu} 
\author[Vakil]{Ravi Vakil} \address{Department of Mathematics,
Stanford University, Stanford, CA 94305-2125}
\email{vakil@math.stanford.edu}
\date{Tuesday, October 17, 2006.}
 \thanks{During the
preparation of this article the second author was partially supported
by NSF Grant DMS-0238532.}

\begin{document}

\begin{abstract}
  We describe recent work on positive descriptions of the structure
  constants of the cohomology of homogeneous spaces such as the
  Grassmannian, by degenerations and related methods. We give various
  extensions of these rules,  some new and conjectural, to $K$-theory,
  equivariant cohomology, equivariant $K$-theory, and quantum
  cohomology.
\end{abstract}

\maketitle
\tableofcontents

\section{Introduction} \label{sec-intro}

\lremind{sec-intro}In this article we describe recent work on positive
algorithms for computing the structure constants of the cohomology of
homogeneous varieties. We also discuss extensions of these rules to
$K$-theory, equivariant cohomology, equivariant $K$-theory and quantum
cohomology. We have two aims. First, we would like to provide the
experts in the field with a compendium of recent results and
references regarding positivity in Schubert calculus. Second, we would
like to present many examples so that the casual user can
perform basic calculations that occur in concrete problems. 

Homogeneous varieties are ubiquitous in mathematics, playing an
important role in representation theory, algebraic and differential
geometry, combinatorics and the theory of symmetric functions.  The
structure constants (``Littlewood-Richardson coefficients'') of the
cohomology rings of homogeneous varieties exhibit a rich and
surprising structure.  For fundamental geometric reasons, the
Littlewood-Richardson coefficients and their generalizations tend to
be positive (interpreted appropriately).  In cohomology, this is a
consequence of Kleiman's Transversality Theorem~\ref{thm-kleiman}
(sometimes called the Kleiman-Bertini Theorem).  In $K$-theory, the
most notable general positivity result is a theorem of Brion
\cite{brion:positivity}, and in equivariant cohomology, the key result
is due to Graham \cite{graham:equivariant},
confirming a conjecture of D.~Peterson.
Such positivity suggests
that these coefficients have a combinatorial interpretation (a
``Littlewood-Richardson rule''), and that such an interpretation
should be geometric in nature.  In recent years positive algorithms
for computing these constants have unraveled some of this beautiful
structure. We now survey the techniques used in obtaining positive
geometric algorithms for determining these structure constants
starting with the case of the ordinary Grassmannians.  

In Part~\ref{partone} we consider the case of type $A$, following a
particular series of degenerations (which first arose in
\cite{vakil:checkers}, and are described in \S
\ref{degenerationorder}) that seems to be particularly fruitful in
resolving a product of Schubert classes into a combination of other
Schubert classes.  In most cases, we may at least conjecturally
interpret these degenerations in terms of generalizations of the
elegant puzzles of A. Knutson and T. Tao.  Many of the new conjectural
statements in this Part are joint work of Knutson and the second
author.


In Part \ref{parttwo}, instead of degenerating according to a fixed
order, we adapt our degeneration order to the Schubert problem at
hand. This flexibility allows us to simplify the geometry. In addition
to a new rule for Grassmannians (\S \ref{secflag}), we obtain
Littlewood-Richardson rules for two-step flag varieties (\S \ref{pmt}) and
the quantum cohomology of Grassmannians (\S \ref{qcg}). Furthermore, the same
degeneration technique can be applied to Fano varieties of quadric
hypersurfaces, thus yielding a method to calculate certain
intersections in Type $B$ and $D$ Grassmannians (\S \ref{sec-orth}).

We will present lots of explicit examples, which are the heart of the
article, as they illustrate the general techniques.  In some sense,
every example is a generalization of a single classical example, \S
\ref{fourlines}.

We emphasize that this is an active, burgeoning area of research, and
we are presenting only a sample of recent work.  We wish to at least
advertise several results closely related to those discussed here.
(i) L. Mihalcea has an explicit statement of a geometric
Littlewood-Richardson rule for the Lagrangian Grassmannian
(type $C$), verified
computationally in a convincing number of cases.  He has proved part
of the statement and is working on the rest.  
(There is already a non-geometric Littlewood-Richardson rule in this
case, due to Stembridge.)
(ii) D. Davis is
pursuing ongoing work on an analogous geometric rule in type $B$.
(iii) W. Graham and S. Kumar have explicitly computed the structure
co-efficients of equivariant $K$-theory of projective space \cite{grahamkumar}.
The description is surprisingly rich and nontrivial given the
``simplicity'' of the space.  Earlier approaches to the equivariant
$K$-theory of $G/B$ are due to Griffeth and Ram \cite{gr} and Willems
\cite{w}.  We haven't yet attempted to relate these results to the
conjecture of \S \ref{KT}.  (iv) Thomas and Yong have recently given
a root-theoretically uniform generalization of the Littlewood-Richardson rule
for intersection numbers of Schubert varieties in minuscule and cominuscule
flag varieties \cite{thomasyong}.  It would be very interesting to 
understand their work geometrically.

There is an even larger amount of work on positive rules without (yet)
direct geometric interpretation.  A discussion of these ideas would
triple the length of the article, so we content ourselves with listing a
sampling of the authors who have contributed to the area: N.~Bergeron,
Fomin, Gelfand, Knutson, Lenart, Postnikov, Robinson, Sottile, Yong,
\dots .

\subsection{(Type $A$) Grassmannians and positivity} Let $G(k,n)$
denote the Grassmannian that parametrizes $k$-dimensional subspaces
of an $n$-dimensional vector space $W$. It is sometimes more
convenient to interpret $G(k,n)$ as the parameter space of
$(k-1)$-dimensional projective linear spaces in $\PP^{n-1}$. When we
wish to emphasize this point of view, we will denote $G(k,n)$ by
$\GG(k-1,n-1)$.

Fix a complete flag $F_{\bullet}$
$$0 =F_0 \subset F_1 \subset \cdots \subset F_n = W.$$ Let $\lambda$
be a partition with $k$ parts satisfying $n-k \geq \lambda_1 \geq
\cdots \geq \lambda_k \geq 0.$ The {\em Schubert variety}
$\Omega_{\lambda}(F_{\bullet})$ of type $\lambda$ associated to the
flag $F_{\bullet}$ is defined by
$$\Omega_{\lambda}(F_{\bullet}) := \{ \ [V] \in G(k,n) \ : \ \dim(
V \cap F_{n-k+i-\lambda_i}) \geq i \ \}.$$ 
Those partitions corresponding to  (non-empty) Schubert
varieties in $G(k,n)$ are readily seen to be those  contained in a
$k \times (n-k)$ rectangle. 
The parts $\lambda_i$
that are zero are often omitted from the notation.

 For example, in
$G(2,4)= \GG(1,3)$ a complete flag is given by specifying a point $p$,
contained in a line $l$, contained in a plane $P$, contained in
$\PP^3$.  The Schubert variety $\Omega_{1}$ parametrizes the lines in
$\PP^3$ that intersect the line $l$. The Schubert variety
$\Omega_{2,1}$ parametrizes the 
lines in $\PP^3$ that are contained in the plane $P$ and contain the
point $p$.

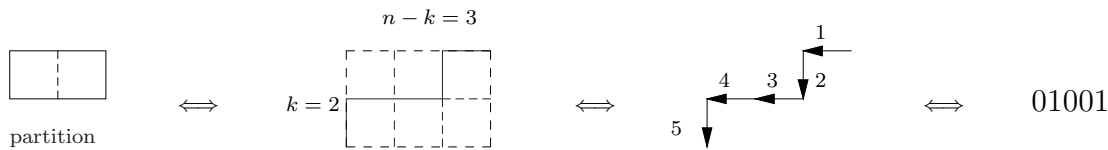
\begin{figure}[ht]
\begin{center} 
\setlength{\unitlength}{0.00083333in}
\begingroup\makeatletter\ifx\SetFigFont\undefined%
\gdef\SetFigFont#1#2#3#4#5{%
  \reset@font\fontsize{#1}{#2pt}%
  \fontfamily{#3}\fontseries{#4}\fontshape{#5}%
  \selectfont}%
\fi\endgroup%
{\renewcommand{\dashlinestretch}{30}
\begin{picture}(7033,885)(0,-10)
\put(12,36){\makebox(0,0)[lb]{{\SetFigFont{8}{9.6}{\rmdefault}{\mddefault}{\updefault}partition}}}
\dashline{60.000}(2712,36)(2712,336)(3012,336)
\dashline{60.000}(2412,636)(2412,36)
\path(2112,36)(2112,336)(2712,336)
	(2712,636)(3012,636)
\blacken\path(4482.000,366.000)(4362.000,336.000)(4482.000,306.000)(4482.000,366.000)
\path(4362,336)(4662,336)
\blacken\path(4782.000,366.000)(4662.000,336.000)(4782.000,306.000)(4782.000,366.000)
\path(4662,336)(4962,336)
\blacken\path(4932.000,456.000)(4962.000,336.000)(4992.000,456.000)(4932.000,456.000)
\path(4962,336)(4962,636)
\blacken\path(5082.000,666.000)(4962.000,636.000)(5082.000,606.000)(5082.000,666.000)
\path(4962,636)(5262,636)
\blacken\path(4332.000,156.000)(4362.000,36.000)(4392.000,156.000)(4332.000,156.000)
\path(4362,36)(4362,336)
\path(612,636)(12,636)(12,336)
	(612,336)(612,636)
\dashline{60.000}(312,336)(312,636)
\put(1737,261){\makebox(0,0)[lb]{{\SetFigFont{8}{9.6}{\rmdefault}{\mddefault}{\updefault}$k=2$}}}
\put(2337,786){\makebox(0,0)[lb]{{\SetFigFont{8}{9.6}{\rmdefault}{\mddefault}{\updefault}$n-k=3$}}}
\put(3537,261){\makebox(0,0)[lb]{{\SetFigFont{8}{9.6}{\rmdefault}{\mddefault}{\updefault}$\Longleftrightarrow$}}}
\put(6387,261){\makebox(0,0)[lb]{{\SetFigFont{12}{14.4}{\rmdefault}{\mddefault}{\updefault}$01001$}}}
\put(5712,261){\makebox(0,0)[lb]{{\SetFigFont{8}{9.6}{\rmdefault}{\mddefault}{\updefault}$\Longleftrightarrow$}}}
\put(4137,111){\makebox(0,0)[lb]{{\SetFigFont{8}{9.6}{\rmdefault}{\mddefault}{\updefault}$5$}}}
\put(4437,411){\makebox(0,0)[lb]{{\SetFigFont{8}{9.6}{\rmdefault}{\mddefault}{\updefault}$4$}}}
\put(4737,411){\makebox(0,0)[lb]{{\SetFigFont{8}{9.6}{\rmdefault}{\mddefault}{\updefault}$3$}}}
\put(5037,411){\makebox(0,0)[lb]{{\SetFigFont{8}{9.6}{\rmdefault}{\mddefault}{\updefault}$2$}}}
\put(5037,711){\makebox(0,0)[lb]{{\SetFigFont{8}{9.6}{\rmdefault}{\mddefault}{\updefault}$1$}}}
\put(1062,261){\makebox(0,0)[lb]{{\SetFigFont{8}{9.6}{\rmdefault}{\mddefault}{\updefault}$\Longleftrightarrow$}}}
\dashline{60.000}(2112,636)(3012,636)(3012,36)
	(2112,36)(2112,636)
\end{picture}
}

\end{center}
\caption{The bijection between partitions contained in a $k \times (n-k)$
rectangle and strings of $k$ $1$'s and $n-k$ $0$'s.
\lremind{seattlerec}
} \label{seattlerec}
\end{figure}

An alternate indexing set for Schubert varieties of $G(k,n)$ consists
of the size $k$ subsets of $\{ 1, \dots, n \}$, presented as a string
of $n$ digits, of which $k$ are $1$ and $n-k$ are $0$.  This is known
as {\em string notation}.  The bijection is as follows.  The $1$'s are
placed at those positions $j \in \{ 1, \dots, n \}$ where $\dim V \cap
F_j > \dim V \cap F_{j-1}$.  A more visual description is given in
Figure~\ref{seattlerec}: to obtain the string of $0$'s and $1$'s,
consider the path from the northeast corner to the southwest corner of
the rectangle, along the border of the partition.  This path consists
of $n$ segments.  If the $j$th step is west (resp.\ south), the $j$th
element of the string is $0$ (resp.\ $1$).

We will use string notation interchangeably with partition notation
(as the index set of Schubert varieties of the Grassmannian), and we
apologize for any confusion this might cause.  Partitions (of length
greater than  $1$) will have commas, and strings will have no commas.


The homology class of a Schubert variety is independent of the
defining flag, and depends only on the partition. We will denote the
Poincar\'e dual of the class of $\Omega_{\lambda}$ by
$\sigma_{\lambda}$.  The Poincar\'e duals of the classes of Schubert
varieties give an additive basis for the cohomology of $G(k,n)$.
Therefore, given two Schubert cycles $\sigma_{\lambda}$,
$\sigma_{\mu}$, their product is a $\Z$-linear combination of Schubert
cycles $\sum_{\nu} c_{\lambda, \mu}^{\nu} \sigma_{\nu}$. The structure
constants $c_{\lambda, \mu}^{\nu}$ of the cohomology ring with respect
to the Schubert basis are called {\em Littlewood-Richardson
  coefficients}. 


The basic problem we would like to address is finding positive
geometric algorithms for computing the Littlewood-Richardson
coefficients. A positive combinatorial rule giving  these
coefficients is called a {\em Littlewood-Richardson rule}. 


 The underlying strategy for the geometric rules is as follows. One
 begins with the intersection of two Schubert varieties defined with
 respect to two general flags. Each rule has a recipe for making the
 flags more special via codimension one degenerations. As one makes
 the flags more special, the intersection of the two Schubert
 varieties defined with respect to the flags become more special,
 sometimes breaking into a number of pieces, each of which is analyzed
 separately in the same way.  The pieces at the end are Schubert
 varieties.  The combinatorial objects such as puzzles or checkers
 encode the varieties that arise as a result of the degenerations. The
 following fundamental example illustrates the strategy.


 \begin{example} \label{fourlines}\lremind{fourlines}How many lines in
   $\PP^3$ intersect four general given lines $l_1, \dots, l_4$?  If
   the lines are in general position, the answer to this question is
   not immediately clear. However, if the lines are in a special
   position, then the answer might be easier to see. Suppose two of
   the lines $l_1$ and $l_2$ intersect at a point $q$. There are two
   ways that a line $l$ can intersect both $l_1$ and $l_2$. If $l$
   passes through the point $q$, then $l$ intersects both. If not, $l$
   must lie in the plane $P$ spanned by $l_1$ and $l_2$. In the first
   case the only line that can also intersect $l_3$ and $l_4$, is the
   intersection of the two planes spanned by $q$ and $l_3$ and $q$ and
   $l_4$. In the second case the only line that intersects $l_3$ and
   $l_4$ is the line that joins $P \cap l_3$ and $P \cap l_4$. We see
   that when the lines are in a slightly special position, then the
   answer to the question is $2$ (see Figure \ref{lines}). 


\begin{figure}[htbp]
\begin{center}
\begin{picture}(0,0)%
\includegraphics{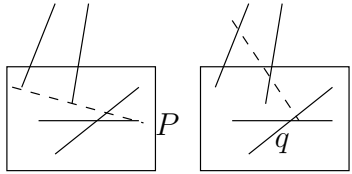}%
\end{picture}%
\setlength{\unitlength}{3947sp}%
\begingroup\makeatletter\ifx\SetFigFont\undefined%
\gdef\SetFigFont#1#2#3#4#5{%
  \reset@font\fontsize{#1}{#2pt}%
  \fontfamily{#3}\fontseries{#4}\fontshape{#5}%
  \selectfont}%
\fi\endgroup%
\begin{picture}(2173,1074)(781,-1198)
\put(2476,-1036){\makebox(0,0)[lb]{\smash{{\SetFigFont{12}{14.4}{\rmdefault}{\mddefault}{\updefault}{\color[rgb]{0,0,0}$q$}%
}}}}
\put(1726,-961){\makebox(0,0)[lb]{\smash{{\SetFigFont{12}{14.4}{\familydefault}{\mddefault}{\updefault}{\color[rgb]{0,0,0}$P$}%
}}}}
\end{picture}%

\end{center}
\caption{There are two lines in $\PP^3$ that intersect four general
lines.\lremind{lines}}
\label{lines}
\end{figure}

An additional argument is required to deduce that if the lines $l_1,
\dots, l_4$ are in general position, the answer is still $2$. Suppose
we begin with the lines in general position and rotate $l_1$ until it
intersects $l_2$. By the principle of conservation of number, the
number of lines intersecting $l_1, \dots, l_4$ counted with
multiplicity remains constant. One may check that the multiplicities
are one by a tangent space calculation that shows that the Schubert
varieties of lines intersecting $l_1$ and $l_2$ intersect generically
transversally even when $l_1$ and $l_2$ intersect at a point.

\end{example}

\subsection{Advantages of geometric rules}
There are many presentations of the cohomology ring of the
Grassmannian. For instance, Pieri's rule gives an easy method to
multiply special Schubert cycles, those Schubert cycles where
$\lambda_i = 0$ except for $i=1$. Let $\sigma_{\lambda}$ be a special
Schubert cycle. Suppose $\sigma_{\mu}$ is any Schubert cycle with
parts $\mu_1, \dots, \mu_k$. Then
$$
\sigma_{\lambda} \cdot \sigma_{\mu} = \sum_{\substack{ \mu_i \leq
\nu_i \leq \mu_{i-1} \\ \sum \nu_i = \lambda + \sum
\mu_i}}\sigma_{\nu}.
$$
Giambelli's determinantal
formula expresses any Schubert cycle as a determinant of a matrix
consisting of only special Schubert cycles: 
$$\sigma_{\lambda_1, \dots, \lambda_k} =
  \begin{vmatrix} \sigma_{\lambda_1} & \sigma_{\lambda_1 + 1} &
  \sigma_{\lambda_1 + 2} & \dots & \sigma_{\lambda_1 + k -1} \\
  \sigma_{\lambda_2 - 1} & \sigma_{\lambda_2} & \sigma_{\lambda_2 + 1}
  & \dots & \sigma_{\lambda_2 + k -2} \\ \vdots & \vdots & \vdots & \ddots & \vdots  \\
  \sigma_{\lambda_k -k+ 1} & \sigma_{\lambda_k - k + 2} &
  \sigma_{\lambda_k - k + 3} & \dots &
  \sigma_{\lambda_k}\end{vmatrix}.$$ These two formulae taken together
  give a presentation for the cohomology ring of $G(k,n)$. 
 However, a geometric Littlewood-Richardson rule has some
advantages over these
  presentations.


  Positive combinatorial rules give more efficient algorithms for
  computing structure constants than determinantal formulae. More
  importantly, some of the structure of the Littlewood-Richardson
  coefficients is very hard to see from determinantal formulae. For
  instance, Horn's conjecture was resolved through the combined work
  of Klyachko, Knutson, and Tao, and the key final step was Knutson
  and Tao's proof of the Saturation conjecture \cite{knutson:puzzle}
  using honeycombs, which are intimately related to puzzles. An
  important consequence of their work is the following property of
  Littlewood-Richardson coefficients, which was originally known as
  Fulton's conjecture:
 
\begin{theorem}(Knutson-Tao-Woodward, \cite{ktw})
If $c_{\lambda,
  \mu}^{\nu} = 1$, then $c_{N
  \lambda, N\mu}^{N \nu} = 1$.
\end{theorem}
\noindent This kind of structure is very hard to prove using
determinantal formulae, but is often immediate from positive
Littlewood-Richardson rules. 


Another advantage of geometric Littlewood-Richardson rules is that
they apply over arbitrary fields.  Over algebraically closed fields of
characteristic zero the connection between cohomology and enumerative
geometry is provided by Kleiman's Transversality Theorem
\cite{kleiman:transverse}.

\begin{theorem}(Kleiman) \label{thm-kleiman}
\lremind{thm-kleiman}Let $G$ be an integral algebraic group scheme, $X$ an integral
algebraic scheme with a transitive $G$ action. Let $f: Y \rightarrow
X$ and $g: Z \rightarrow X$ be two maps of integral algebraic
schemes. For each rational element of $g \in G$, denote by $gY$ the
$X$-scheme given by $y \mapsto g f(y)$. Then there exists a dense open
subset $U$ of $G$ such that for every rational element $g \in U$, the
fiber product $(gY) \times_X Z$ is either empty or equidimensional of
the expected dimension $$\dim Y + \dim Z - \dim X.$$
Furthermore, if $Y$ and $Z$ are regular, then for a dense open set
this fibered product is regular.
\end{theorem} 

\noindent In particular, applying Theorem \ref{thm-kleiman} in the
case when $X$ is a homogeneous variety and $f$ and $g$ are the
inclusions of Schubert subvarieties, we conclude that for general
translates the intersections will be generically transverse of the
expected dimension provided the intersection is
non-empty. Unfortunately, Kleiman's theorem does not hold when the
characteristic of the ground field is not zero, or when the
ground field is not algebraically closed. 

The geometric Littlewood-Richardson rules apply even when the ground
field has positive characteristic. The following is a
characteristic-free version of the Kleiman-Bertini Theorem proved in
\cite{vakil:induction}.

\begin{theorem}(Generic smoothness)
Suppose $Q \subset G(k,n)$ is a subvariety such that $(Q \times Fl(n))
\cap \Omega_{\lambda}(F_{\bullet}) \rightarrow Fl(n)$ is generically
smooth for all $\lambda$. Then $$(Q \cap Fl(n)^m) \cap
\pi_1^*\Omega_{\lambda_1}(F_{\bullet}^1) \cap \cdots \cap
\pi_m^*\Omega_{\lambda_m} (F_{\bullet}^m) \rightarrow Fl(n)^m$$ is
also generically smooth.
\end{theorem}

Let $F_{\bullet}^1, \dots, F_{\bullet}^m$ be $m$ general flags in
$\CC^n$. Suppose $\Omega_{\lambda_1}(F_{\bullet}^1), \dots,
\Omega_{\lambda_m}(F_{\bullet}^m)$ are $m$ Schubert varieties in
$G(k,n)$ whose dimension of intersection is zero. The corresponding
Schubert problem asks for the cardinality of the intersection of these
varieties. More generally, a Schubert problem asks for the cardinality
of the intersection of Schubert varieties defined with respect to
general flags in case the dimension of intersection is zero. Over the
complex numbers this cardinality can be determined by computing the
degree of the intersection in the cohomology ring by Kleiman's
theorem. However, a priori, over other fields it is not clear that the
cardinality will equal the degree.  We say that a Schubert problem is
enumerative over a field $\kbb$ if there exists flags defined over
$\kbb$ such that the cardinality of the $\kbb$-points in the
intersection of the Schubert varieties defined with respect to these
flags is equal to the degree of intersection of the Schubert
varieties. 

One application of geometric Littlewood-Richardson rules has been to
resolve the question of whether Schubert problems are enumerative over
the real numbers, finite fields and algebraically closed fields of any
characteristic. The theorems may be found in \cite{vakil:induction}.

\begin{theorem}
\begin{enumerate}
\item All Schubert problems for the Grassmannians are enumerative over
  the real numbers (in fact for any field satisfying the implicit
  function theorem).

\item All Schubert problems are enumerative over algebraically closed
  fields (of arbitrary characteristic).

\end{enumerate}
\end{theorem}

Despite the progress many natural questions still remain.

\begin{question}
Are Schubert problems enumerative over the rational numbers?
\end{question}

There are many other applications of geometric Littlewood-Richardson
rules.  For instance, they can be used to compute monodromy groups of
Schubert problems. One may construct interesting examples of Schubert
problems whose monodromy group is not the full symmetric group. For
precise details about this and other applications see
\cite{vakil:induction}. 

One final appeal of the geometric rules over other combinatorial
positive rules is that the techniques extend to other homogeneous
varieties. We will now detail the extension of these ideas to two-step
flag varieties and to orthogonal and symplectic Grassmannians.

\subsection{Flag varieties}

Let $0 < k_1 < \cdots < k_r < n$ be an increasing sequence of $r$
positive integers.  Let $\Flag$ denote the $r$-step flag variety of
$r$-tuples of linear subspaces $(V_1, \cdots, V_r)$ of an
$n$-dimensional vector space $W$, where $V_i$ are $k_i$-dimensional
linear spaces and $V_i \subset V_{i+1}$ for all $1 \leq i \leq r-1$.
When we would like to consider the flag variety as a parameter space
for nested sequences of linear subspaces of projective space, we will
use the notation $\FF l(k_1 - 1, \dots, k_r -1; n-1)$.

The cohomology groups of flag varieties are also generated by Schubert
classes. Usually in the literature the Schubert cycles are
parametrized by certain permutations. More precisely, Schubert
varieties are parametrized by permutations $\omega$ of length $n$ for
which $\omega(i) < \omega(i+1)$ whenever $i \notin \{k_1, \dots, k_r
\}$. More explicitly, the Schubert variety $X_{\omega} (F_{\bullet})$
is defined by
\begin{align*}
  X_{\omega} (F_{\bullet}) : = \{ \ ( V_1, \dots, V_r) \in \Flag \ |
   \dim (V_i \cap F_j) \geq \# \{ \alpha \leq i : \omega(\alpha) > n-j
   \} \ \forall \ i,j \ \}.
\end{align*}
The Poincar\'e duals of the classes of all the Schubert varieties form
an additive basis for the cohomology of the flag variety.

For our purposes two other notations for Schubert varieties of
$r$-step flag varieties are useful.

First, in analogy with the partition notation for the Grassmannians we
will use the notation $\sigma_{\lambda_1, \cdots, \lambda_{k_r}}^{
  \delta_1, \cdots, \delta_{k_r}}$. The bottom row denotes the usual
partition corresponding to the $k_r$-plane $V_r$ in $W$ treated as a
Schubert cycle in $G(k_r, n)$. The numbers $\delta_i$ are integers
between $1$ and $r$. For a Schubert cycle in $\Flag$, $k_1$ of the
upper indices will be $1$ and $k_i - k_{i-1}$ of them will be $i$. The
flag $F_{\bullet}$ induces a complete flag $G_{\bullet}$ on the
largest vector space $V_r$.  For each $j$, there exists a smallest $i$
such that $$\dim (V_i \cap G_j) = \dim (V_i \cap G_{j-1}) + 1.
$$
For a
Zariski-open subset of the Schubert variety this index will be
constant. In that case we write $i$ on top of $\lambda_j$. In the case
of Grassmannians this notation reduces to the ordinary notation with a
sequence of $1$'s on the top row. For the complete flag variety the top
row becomes the permutation defining the Schubert cycle. 

Second, there is a string notation for partial flag varieties.
Schubert varieties of $\Flag$ are
indexed by $n$-tuples of numbers $0$, \dots, $r-1$, where $k_{i+1} -
k_i$ of the digits are $i$ (where $k_0 := 0$).  
Taking
the convention $k_{r+1} = n$, 
for each $j$, there exists a smallest $i$
such that $$\dim (V_i \cap F_j) = \dim (V_i \cap F_{j-1}) + 1.
$$
Then in position $j$, we place $r+1-i$.  The reader may verify that
this generalizes the string notation for the Grassmannian.

We will use both notations interchangeably. As with the Grassmannian
case, strings will have no commas.

\begin{example} Fix a flag $F_1 \subset \cdots \subset F_6$ in $W^6$.
  The Schubert cycle $\sigma_{2,1,0}^{2,1,2}$ in $Fl(1,3;6)$ denotes
  the pairs of subspaces $V_1 \subset V_2$ where $V_1$ has dimension
  one and $V_2$ has dimension 3. $V_2$ is required to meet $F_2$ in
  dimension one, $F_4$ in dimension $2$ and be contained in $F_6$.
  $V_1$ lies in the intersection of $V_2$ with $F_4$.  In string
  notation, this is written $\si_{010201}$.
\end{example}

Using the ideas of geometric rules for Grassmannians it is possible to
give an explicit rule for two-step flag varieties (see
\cite{coskun1:LR}, discussed in \S \ref{secflag}). Previously there
were no proven rules for partial flag varieties. Knutson conjectured a
rule for two-step flag varieties in terms of puzzles
(Conjecture~\ref{allentwostep}).  A. Buch extended the conjecture to
three-step flag varieties, described in Section \ref{s:anders}.  We
give a conjectural geometric rule 
(due in part to Knutson) based on Knutson's two-step
conjecture in Section \ref{s:twostep}.

\subsection{Other groups}\label{othergps}

So far we have discussed the ordinary Grassmannians and flag
varieties. These varieties are associated to the Lie group $SL(n)$
(Type $A$). One can consider homogeneous varieties associated to other
classical groups such as $SO(n)$, $Sp(2n)$ and the exceptional groups.
In this paper we will have nothing to say about the exceptional Lie
groups. However, we would like to explain the geometric point of view
for the other infinite families of Lie groups.  

Let $Q$ be a non-degenerate, symmetric (resp., alternating) bilinear
form on an $n$-dimensional vector space $W$. The isotropic
Grassmannians $OG(k,n)$ (resp., $SG(k,n)$) parametrize
$k$-dimensional subspaces of $W$ which are isotropic with respect to
$Q$. When $n$ is even, $k=n/2$ and $Q$ is symmetric, the isotropic
subspaces form two isomorphic irreducible families. In that case, the
orthogonal Grassmannian is one of the irreducible components.  The
varieties $OG(k,n)$ are quotients of the Lie group $SO(n)$ by maximal
parabolic subgroups. They are classified into Type $B$ and Type $D$
depending on the parity of $n$. The varieties $SG(k,n)$ are quotients
of $Sp(n)$ by maximal parabolic subgroups.  Since there cannot be a
non-degenerate alternating form on an odd-dimensional vector space,
here $n$ has to be even. This is the Type $C$
case.

The cohomology of isotropic Grassmannians is generated by Schubert
varieties. The Littlewood-Richardson coefficients of Type $B$ and Type
$C$ flag varieties are equal up to an explicit power of $2$.  The
relation, which we will summarize in Lemma \ref{sot}, is discussed in
\cite{sottile:pieri} p. 17 or \cite{billey:polynomials}. Let $u$, $v$
and $w$ be three permutations. For a permutation $u$ denote $s(u)$ the
number of sign changes of the permutation.

\begin{lemma}\label{sot}
The Littlewood-Richardson coefficients of Type $B$ and Type $C$ full-flag
varieties satisfy 
$$2^{s(u)+s(v)} b_{u,v}^w = 2^{s(w)} c_{u,v}^w,$$ where $b_{u,v}^w$
(resp., $ c_{u,v}^w$) denotes the structure coefficients of Type $B$
(resp., Type $C$) full-flag variety.  \end{lemma} In view of Lemma
\ref{sot}, we can restrict our attention to orthogonal Grassmannians.
There are minor differences in the description of the Schubert
varieties depending on the type of the isotropic Grassmannian. For
ease of exposition, we will describe the geometric viewpoint for the
case $OG(k,2m+1)$.  

Let $s \leq m$ be a non-negative integer. Let $\lambda$ denote a
strictly decreasing partition $$ m \geq \lambda_1 > \lambda_2 > \cdots
> \lambda_s > 0.$$ Given $\lambda$, there is an associated partition
$\tilde{\lambda}$
$$ m-1 \geq \tilde{\lambda}_{s+1} > \tilde{\lambda}_{s+2} > \cdots >
\tilde{\lambda}_m \geq 0 $$ defined by the requirement that there does
not exist any parts $\lambda_i$ for which $\tilde{\lambda}_j +
\lambda_i = m$. In other words, the associated partition is obtained
by removing the integers $m-\lambda_1, \dots, m-\lambda_s$ from the
sequence $m-1, m-2, \dots, 0$. For example, if $m=6$, then the
partition associated to $(6,4)$ is $(5,4,3,1)$. The Schubert varieties
in $OG(k,2m+1)$ are parametrized by pairs $(\lambda, \mu)$, where
$\lambda$ is a strictly decreasing partition of length $s$ and $\mu$
$$\mu_{s+1} > \mu_{s+2} > \cdots > \mu_k \geq 0$$ is a subpartition of
$\tilde{\lambda}$ of length $k-s$. Given a pair $(\lambda, \mu)$ the
{\em discrepancy} $\dis(\lambda,\mu)$ of the pair is defined as follows:
Since $\mu$ is a subpartition of $\tilde{\lambda}$ we can assume that
it parts occur as $\tilde{\lambda}_{i_{s+1}}, \cdots,
\tilde{\lambda}_{i_k}$. The {discrepancy} is given by
$$\dis(\lambda,\mu) = (m-k)s + \sum_{j=s+1}^k (m -k + j - i_j).$$ 


Fix an isotropic flag $F_{\bullet}$
$$0 \subset F_1 \subset F_2 \subset \cdots \subset F_m \subset
F_{m-1}^{\perp} \subset \cdots \subset F_1^{\perp} \subset W.$$ Here
$F_i^{\perp}$ denotes the orthogonal complement of $F_i$ with respect
to the bilinear form. The Schubert variety
$\Omega_{\lambda,\mu}(F_{\bullet})$ is defined as the closure of the
locus $$ \{ [V] \in OG(k,2m+1)\ | \ \dim( V \cap F_{m+1 - \lambda_i})
= i, \ \dim( V \cap F_{\mu_j}^{\perp}) = j \ \}. $$ The codimension of
a Schubert variety is given by $\sum_{i=1}^s \lambda_i +
\dis(\lambda,\mu).$ We will denote the Poincar\'e dual of the
cohomology class of $\Omega_{\lambda, \mu}$ by $\sigma_{\lambda,\mu}$.
Observe that for maximal isotropic Grassmannians $OG(m,2m+1)$, the
partition $\mu$ is uniquely determined by the partition $\lambda$.
Consequently, in the literature the sequence $\mu$ is often omitted
form the notation.  

Geometrically, the orthogonal Grassmannian $OG(k,2m+1)$ may be
interpreted as the Fano variety of $(k-1)$-dimensional projective
linear spaces on a smooth quadric hypersurface in $\PP^{2m}$. The
non-degenerate form $Q$ defines the smooth quadric hypersurface
$\mathrm{Q}$. A linear space $F_i$ is isotropic with respect to $Q$ if
and only if its projectivization is contained in $\mathrm{Q}$. The
projectivization of the orthogonal complement $F_i^{\perp}$
corresponds to the linear space of codimension $i$ everywhere tangent
to $\mathrm{Q}$ along the projectivization of $F_i$. This geometric
reinterpretation allows us to apply modifications of the Mondrian
tableaux rule to perform calculations in the cohomology of
$OG(k,2m+1)$ (see \S \ref{sec-orth}).  

In order to adapt the previous discussion to Grassmannians $OG(k,2m)$,
we have to account for the existence of two isomorphic irreducible
families of $m$-dimensional isotropic linear spaces. The
$m$-dimensional linear spaces belong to the same irreducible family if
and only if their dimension of intersection has the same parity as $m$. 
Let $\lambda$ be a strictly decreasing partition
$$m-1 \geq \lambda_1 > \cdots > \lambda_s \geq 0,$$ where $s$ has
the same parity as $m$. Define the associated
partition as those integers $$m-1 \geq \tilde{\lambda}_{s+1} > \cdots
> \lambda_m \geq 0,$$ whose sums with any $\lambda_i$ is not equal to
$m-1$.  The Schubert varieties in $OG(k,2m)$ are parametrized by
pairs of strictly decreasing partitions $(\lambda, \mu)$, where $\mu$
is a subpartition of $\tilde{\lambda}$ of length $k-s$. With these
modifications in the numerics, the discussion of the Schubert
varieties of $OG(k,2m+1)$ carries over to the case of $OG(k,2m)$. We
will leave the details to the reader.

\noindent {\bf Acknowledgments.}  We thank Sara Billey, Anders Buch,
Diane Davis, William Graham, Joe Harris, Allen Knutson, Shrawan Kumar,
Robert MacPherson, Leonardo Mihalcea, Harry Tamvakis, and the
organizers of the 2005 Seattle conference in algebraic geometry.  In
particular, we are grateful to Buch and Knutson for allowing us to
describe their conjectures, and Graham and Kumar for allowing us to
present their as-yet-unpublished work \cite{grahamkumar}.  We thank
Buch, Knutson, Mihalcea, and Tamvakis for detailed helpful comments on the
manuscript.

\part{TYPE $A$ RULES, USING A SPECIFIC DEGENERATION ORDER}

\label{partone}
In the next section, we will describe the geometric
Littlewood-Richardson rule \cite{vakil:checkers}, but in the language
of puzzles of Knutson and Tao.  In later sections in
Part~\ref{partone}, we will use the same degeneration order to extend
these ideas conjecturally to much more general situations ($K$-theory,
equivariant cohomology, equivariant $K$-theory, $2$-step flag
varieties, and partial flag varieties in general).  The key
construction will be varieties associated to partially completed puzzles
(\eqref{puzzlevariety}, and a variation \eqref{eqPV} in the
equivariant case).  The puzzles are very friendly to use, and the
reader is very strongly encouraged to work through the examples.

\section{The Grassmannian}

\label{hhh}\lremind{hhh}We now describe the geometric
Littlewood-Richardson rule for the Grassmannian $G(k,n)$, following
\cite{vakil:checkers}.  We begin by describing the sequence of
degenerations.  We will be considering the intersection of Schubert
varieties with respect to two transverse flags $F_{\bullet}$ and
$M_{\bullet}$ in $Fl(n)$ (where $F$ and $M$ stand for $F$ixed and
$M$oving flags, respectively).  We describe a series of ``codimension
$1$'' degenerations of $M_{\bullet}$, beginning with $M_{\bullet}$
transverse to $F_{\bullet}$, and ending with $M_{\bullet} =
F_{\bullet}$.  This degeneration order is very special; it appears
that the good behavior described below happens in general only for
this order (and the ``dual'' order).

\subsection{A key example}
\label{akeyexample}
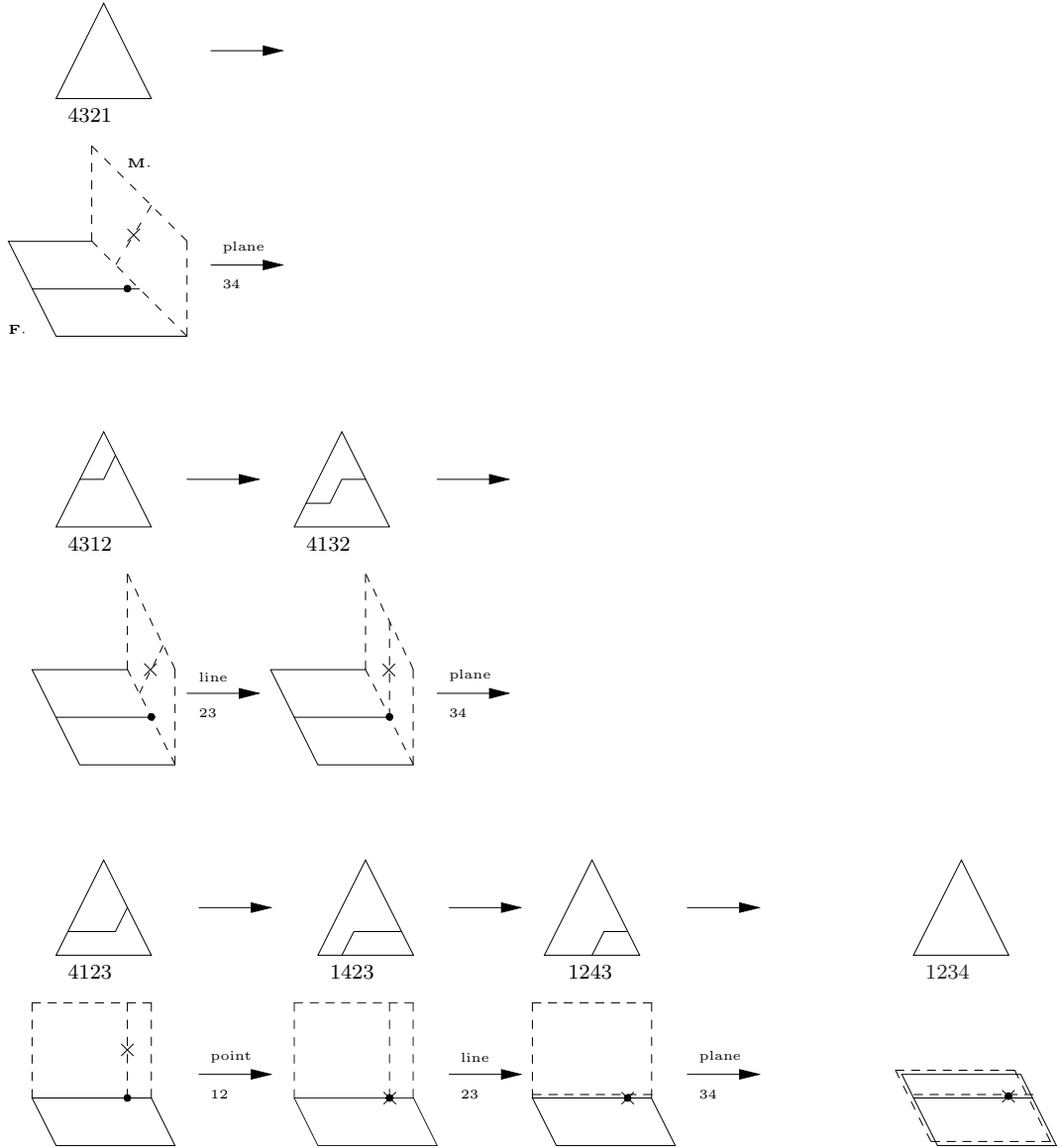
\begin{figure}[htbp]
\begin{center}
\setlength{\unitlength}{0.00083333in}
\begingroup\makeatletter\ifx\SetFigFont\undefined%
\gdef\SetFigFont#1#2#3#4#5{%
  \reset@font\fontsize{#1}{#2pt}%
  \fontfamily{#3}\fontseries{#4}\fontshape{#5}%
  \selectfont}%
\fi\endgroup%
{\renewcommand{\dashlinestretch}{30}
\begin{picture}(6624,7239)(0,-10)
\put(5787,1062){\makebox(0,0)[lb]{{\SetFigFont{8}{9.6}{\rmdefault}{\mddefault}{\updefault}$1234$}}}
\path(1812,312)(1962,12)(2712,12)
	(2562,312)(1812,312)
\dashline{60.000}(1812,312)(1812,912)(2562,912)(2562,312)
\dashline{60.000}(2412,912)(2412,312)
\put(3912,312){\blacken\ellipse{40}{40}}
\put(3912,312){\ellipse{40}{40}}
\path(3312,312)(3462,12)(4212,12)
	(4062,312)(3312,312)
\dashline{60.000}(3312,312)(3312,912)(4062,912)(4062,312)
\dashline{60.000}(3315,334)(4065,334)
\path(764,5788)(839,5713)
\path(839,5788)(764,5713)
\path(869,3050)(944,2975)
\path(944,3050)(869,2975)
\path(2370,3047)(2445,2972)
\path(2445,3047)(2370,2972)
\path(723,653)(798,578)
\path(798,653)(723,578)
\path(2376,349)(2451,274)
\path(2451,349)(2376,274)
\path(3874,349)(3949,274)
\path(3949,349)(3874,274)
\path(6272,361)(6347,286)
\path(6347,361)(6272,286)
\put(6310,324){\blacken\ellipse{40}{40}}
\put(6310,324){\ellipse{40}{40}}
\path(5637,462)(5862,12)(6612,12)
	(6387,462)(5637,462)
\dashline{60.000}(5715,334)(6465,334)
\path(5712,312)(6462,312)
\dashline{60.000}(5598,487)(5823,37)(6573,37)
	(6348,487)(5598,487)
\put(762,5412){\blacken\ellipse{40}{40}}
\put(762,5412){\ellipse{40}{40}}
\put(912,2712){\blacken\ellipse{40}{40}}
\put(912,2712){\ellipse{40}{40}}
\put(2412,2712){\blacken\ellipse{40}{40}}
\put(2412,2712){\ellipse{40}{40}}
\put(762,312){\blacken\ellipse{40}{40}}
\put(762,312){\ellipse{40}{40}}
\path(4287,1512)(4737,1512)
\blacken\path(4617.000,1482.000)(4737.000,1512.000)(4617.000,1542.000)(4617.000,1482.000)
\path(4287,462)(4737,462)
\blacken\path(4617.000,432.000)(4737.000,462.000)(4617.000,492.000)(4617.000,432.000)
\path(2787,462)(3237,462)
\blacken\path(3117.000,432.000)(3237.000,462.000)(3117.000,492.000)(3117.000,432.000)
\path(2787,1512)(3237,1512)
\blacken\path(3117.000,1482.000)(3237.000,1512.000)(3117.000,1542.000)(3117.000,1482.000)
\path(1212,1512)(1662,1512)
\blacken\path(1542.000,1482.000)(1662.000,1512.000)(1542.000,1542.000)(1542.000,1482.000)
\path(1212,462)(1662,462)
\blacken\path(1542.000,432.000)(1662.000,462.000)(1542.000,492.000)(1542.000,432.000)
\path(1287,6912)(1737,6912)
\blacken\path(1617.000,6882.000)(1737.000,6912.000)(1617.000,6942.000)(1617.000,6882.000)
\path(1287,5562)(1737,5562)
\blacken\path(1617.000,5532.000)(1737.000,5562.000)(1617.000,5592.000)(1617.000,5532.000)
\path(1137,4212)(1587,4212)
\blacken\path(1467.000,4182.000)(1587.000,4212.000)(1467.000,4242.000)(1467.000,4182.000)
\path(1137,2862)(1587,2862)
\blacken\path(1467.000,2832.000)(1587.000,2862.000)(1467.000,2892.000)(1467.000,2832.000)
\path(2712,4212)(3162,4212)
\blacken\path(3042.000,4182.000)(3162.000,4212.000)(3042.000,4242.000)(3042.000,4182.000)
\path(2712,2862)(3162,2862)
\blacken\path(3042.000,2832.000)(3162.000,2862.000)(3042.000,2892.000)(3042.000,2832.000)
\path(537,5712)(12,5712)(312,5112)(1137,5112)
\dashline{60.000}(912,5937)(687,5562)
\dashline{60.000}(537,5712)(1137,5112)(1137,5712)
	(537,6312)(537,5712)
\path(162,5412)(837,5412)
\path(762,3012)(162,3012)(462,2412)(1062,2412)
\dashline{60.000}(762,3012)(1062,2412)(1062,3012)
	(762,3612)(762,3012)
\path(312,2712)(912,2712)
\dashline{60.000}(987,3162)(837,2862)
\path(2262,3012)(1662,3012)(1962,2412)(2562,2412)
\dashline{60.000}(2262,3012)(2562,2412)(2562,3012)
	(2262,3612)(2262,3012)
\path(1812,2712)(2412,2712)
\dashline{60.000}(2412,3312)(2412,2712)
\path(162,312)(312,12)(1062,12)
	(912,312)(162,312)
\dashline{60.000}(162,312)(162,912)(912,912)(912,312)
\dashline{60.000}(762,912)(762,312)
\path(912,6612)(612,7212)(312,6612)(912,6612)
\path(912,3912)(612,4512)(312,3912)(912,3912)
\path(2412,3912)(2112,4512)(1812,3912)(2412,3912)
\path(912,1212)(612,1812)(312,1212)(912,1212)
\path(2562,1212)(2262,1812)(1962,1212)(2562,1212)
\path(3987,1212)(3687,1812)(3387,1212)(3987,1212)
\path(6312,1212)(6012,1812)(5712,1212)(6312,1212)
\path(687,4362)(612,4212)(462,4212)
\path(2262,4212)(2112,4212)(2037,4062)(1887,4062)
\path(762,1512)(687,1362)(387,1362)
\path(2487,1362)(2187,1362)(2112,1212)
\path(3912,1362)(3762,1362)(3687,1212)
\put(1287,537){\makebox(0,0)[lb]{{\SetFigFont{5}{6.0}{\rmdefault}{\mddefault}{\updefault}point}}}
\put(2862,537){\makebox(0,0)[lb]{{\SetFigFont{5}{6.0}{\rmdefault}{\mddefault}{\updefault}line}}}
\put(4362,537){\makebox(0,0)[lb]{{\SetFigFont{5}{6.0}{\rmdefault}{\mddefault}{\updefault}plane}}}
\put(387,1062){\makebox(0,0)[lb]{{\SetFigFont{8}{9.6}{\rmdefault}{\mddefault}{\updefault}$4123$}}}
\put(2037,1062){\makebox(0,0)[lb]{{\SetFigFont{8}{9.6}{\rmdefault}{\mddefault}{\updefault}$1423$}}}
\put(3537,1062){\makebox(0,0)[lb]{{\SetFigFont{8}{9.6}{\rmdefault}{\mddefault}{\updefault}$1243$}}}
\put(1287,312){\makebox(0,0)[lb]{{\SetFigFont{5}{6.0}{\rmdefault}{\mddefault}{\updefault}$12$}}}
\put(2862,312){\makebox(0,0)[lb]{{\SetFigFont{5}{6.0}{\rmdefault}{\mddefault}{\updefault}$23$}}}
\put(4362,312){\makebox(0,0)[lb]{{\SetFigFont{5}{6.0}{\rmdefault}{\mddefault}{\updefault}$34$}}}
\put(762,6162){\makebox(0,0)[lb]{{\SetFigFont{5}{6.0}{\rmdefault}{\mddefault}{\updefault}${\mathbf{M}}_{\cdot}$}}}
\put(387,6462){\makebox(0,0)[lb]{{\SetFigFont{8}{9.6}{\rmdefault}{\mddefault}{\updefault}$4321$}}}
\put(12,5112){\makebox(0,0)[lb]{{\SetFigFont{5}{6.0}{\rmdefault}{\mddefault}{\updefault}${\mathbf{F}}_{\cdot}$}}}
\put(1362,5637){\makebox(0,0)[lb]{{\SetFigFont{5}{6.0}{\rmdefault}{\mddefault}{\updefault}plane}}}
\put(1362,5412){\makebox(0,0)[lb]{{\SetFigFont{5}{6.0}{\rmdefault}{\mddefault}{\updefault}$34$}}}
\put(1212,2937){\makebox(0,0)[lb]{{\SetFigFont{5}{6.0}{\rmdefault}{\mddefault}{\updefault}line}}}
\put(2787,2937){\makebox(0,0)[lb]{{\SetFigFont{5}{6.0}{\rmdefault}{\mddefault}{\updefault}plane}}}
\put(387,3762){\makebox(0,0)[lb]{{\SetFigFont{8}{9.6}{\rmdefault}{\mddefault}{\updefault}$4312$}}}
\put(1887,3762){\makebox(0,0)[lb]{{\SetFigFont{8}{9.6}{\rmdefault}{\mddefault}{\updefault}$4132$}}}
\put(1212,2712){\makebox(0,0)[lb]{{\SetFigFont{5}{6.0}{\rmdefault}{\mddefault}{\updefault}$23$}}}
\put(2787,2712){\makebox(0,0)[lb]{{\SetFigFont{5}{6.0}{\rmdefault}{\mddefault}{\updefault}$34$}}}
\put(2412,312){\blacken\ellipse{40}{40}}
\put(2412,312){\ellipse{40}{40}}
\end{picture}
}

\end{center}
\caption{The degeneration order for $n=4$ (shown projectively,
in $\PP^3$).  The puzzles will be explained shortly.\lremind{degenorder}
}
\label{degenorder}
\end{figure}

\lremind{akeyexample}Before describing the order of degenerations in general, we begin with an example,
shown in Figure~\ref{degenorder}, for $n=4$, but shown projectively
(in $\PP^3$) for convenience of visualization. 
 In the first
degeneration, the moving plane moves (and all other
parts of the moving flag are stationary) until ``something unusual''
happens, which is when it contains the fixed point.  Then the moving
line (in the moving plane, containing the moving point) moves until it too
contains the fixed point.  Then the moving plane moves again (around
the moving line, which remains stationary) until it contains the fixed
line.  Then the moving point moves down to agree with the fixed point,
then the moving line pivots to agree with the fixed line, and finally
the moving plane closes like a book to agree with the fixed plane.

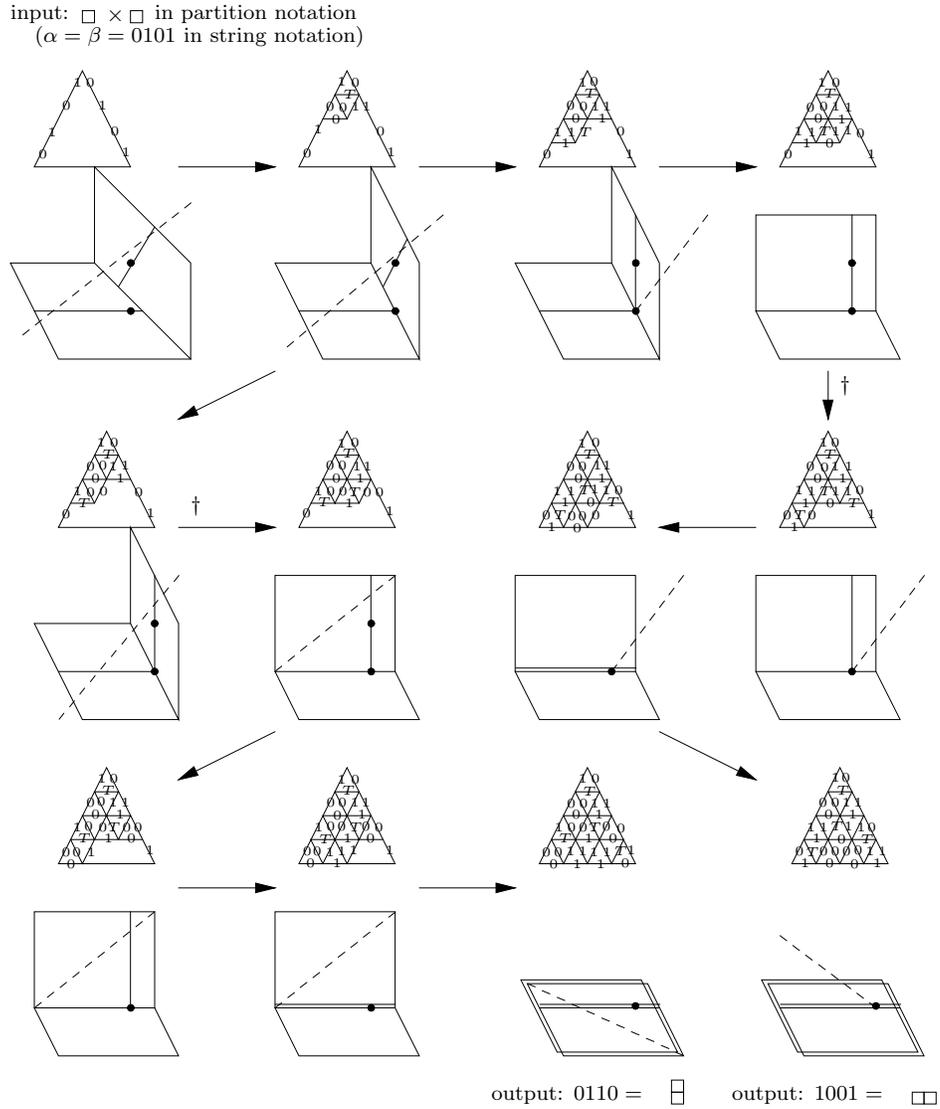
\begin{figure}[htbp]
\begin{center}
\setlength{\unitlength}{0.00083333in}
\begingroup\makeatletter\ifx\SetFigFont\undefined%
\gdef\SetFigFont#1#2#3#4#5{%
  \reset@font\fontsize{#1}{#2pt}%
  \fontfamily{#3}\fontseries{#4}\fontshape{#5}%
  \selectfont}%
\fi\endgroup%
{\renewcommand{\dashlinestretch}{30}
\begin{picture}(5799,6897)(0,-10)
\put(5422,1592){\makebox(0,0)[lb]{{\SetFigFont{5}{6.0}{\rmdefault}{\mddefault}{\updefault}$1$}}}
\put(5262,4986){\blacken\ellipse{40}{40}}
\put(5262,4986){\ellipse{40}{40}}
\path(4662,4986)(4812,4686)(5562,4686)
	(5412,4986)(4662,4986)
\path(4662,4986)(4662,5586)(5412,5586)(5412,4986)
\path(5262,5586)(5262,4986)
\put(762,4986){\blacken\ellipse{40}{40}}
\put(762,4986){\ellipse{40}{40}}
\put(762,5286){\blacken\ellipse{40}{40}}
\put(762,5286){\ellipse{40}{40}}
\put(2412,4986){\blacken\ellipse{40}{40}}
\put(2412,4986){\ellipse{40}{40}}
\put(2412,5286){\blacken\ellipse{40}{40}}
\put(2412,5286){\ellipse{40}{40}}
\put(3912,4986){\blacken\ellipse{40}{40}}
\put(3912,4986){\ellipse{40}{40}}
\put(3912,5286){\blacken\ellipse{40}{40}}
\put(3912,5286){\ellipse{40}{40}}
\path(537,5286)(12,5286)(312,4686)(1137,4686)
\path(537,5286)(1137,4686)(1137,5286)
	(537,5886)(537,5286)
\path(162,4986)(837,4986)
\path(912,5511)(687,5136)
\path(2262,5286)(1662,5286)(1962,4686)(2562,4686)
\path(2262,5286)(2562,4686)(2562,5286)
	(2262,5886)(2262,5286)
\path(1812,4986)(2412,4986)
\path(2487,5436)(2337,5136)
\path(3762,5286)(3162,5286)(3462,4686)(4062,4686)
\path(3762,5286)(4062,4686)(4062,5286)
	(3762,5886)(3762,5286)
\path(3312,4986)(3912,4986)
\path(3912,5586)(3912,4986)
\put(2262,3036){\blacken\ellipse{40}{40}}
\put(2262,3036){\ellipse{40}{40}}
\put(2262,2736){\blacken\ellipse{40}{40}}
\put(2262,2736){\ellipse{40}{40}}
\path(1662,2736)(1812,2436)(2562,2436)
	(2412,2736)(1662,2736)
\path(1662,2736)(1662,3336)(2412,3336)(2412,2736)
\path(2262,3336)(2262,2736)
\put(5262,2736){\blacken\ellipse{40}{40}}
\put(5262,2736){\ellipse{40}{40}}
\path(4662,2736)(4812,2436)(5562,2436)
	(5412,2736)(4662,2736)
\path(4662,2736)(4662,3336)(5412,3336)(5412,2736)
\path(5262,3336)(5262,2736)
\put(3762,2736){\blacken\ellipse{40}{40}}
\put(3762,2736){\ellipse{40}{40}}
\path(3162,2736)(3312,2436)(4062,2436)
	(3912,2736)(3162,2736)
\path(3165,2758)(3915,2758)
\path(3162,2736)(3162,3336)(3912,3336)(3912,2736)
\put(912,2736){\blacken\ellipse{40}{40}}
\put(912,2736){\ellipse{40}{40}}
\put(912,3036){\blacken\ellipse{40}{40}}
\put(912,3036){\ellipse{40}{40}}
\path(762,3036)(162,3036)(462,2436)(1062,2436)
\path(762,3036)(1062,2436)(1062,3036)
	(762,3636)(762,3036)
\path(312,2736)(912,2736)
\path(912,3336)(912,2736)
\put(762,636){\blacken\ellipse{40}{40}}
\put(762,636){\ellipse{40}{40}}
\path(162,636)(312,336)(1062,336)
	(912,636)(162,636)
\path(162,636)(162,1236)(912,1236)(912,636)
\path(762,1236)(762,636)
\put(2262,636){\blacken\ellipse{40}{40}}
\put(2262,636){\ellipse{40}{40}}
\path(1662,636)(1812,336)(2562,336)
	(2412,636)(1662,636)
\path(1665,658)(2415,658)
\path(1662,636)(1662,1236)(2412,1236)(2412,636)
\put(3910,648){\blacken\ellipse{40}{40}}
\put(3910,648){\ellipse{40}{40}}
\path(3237,786)(3462,336)(4212,336)
	(3987,786)(3237,786)
\path(3315,658)(4065,658)
\path(3312,636)(4062,636)
\path(3198,811)(3423,361)(4173,361)
	(3948,811)(3198,811)
\put(5410,648){\blacken\ellipse{40}{40}}
\put(5410,648){\ellipse{40}{40}}
\path(4737,786)(4962,336)(5712,336)
	(5487,786)(4737,786)
\path(4815,658)(5565,658)
\path(4812,636)(5562,636)
\path(4698,811)(4923,361)(5673,361)
	(5448,811)(4698,811)
\path(762,6861)(837,6861)(837,6786)
	(762,6786)(762,6861)
\path(4137,186)(4212,186)(4212,111)
	(4137,111)(4137,186)
\path(4137,111)(4212,111)(4212,36)
	(4137,36)(4137,111)
\dashline{60.000}(1137,5661)(87,4836)
\dashline{60.000}(1737,4761)(2712,5586)
\path(1062,5886)(1662,5886)
\blacken\path(1542.000,5856.000)(1662.000,5886.000)(1542.000,5916.000)(1542.000,5856.000)
\path(2562,5886)(3162,5886)
\blacken\path(3042.000,5856.000)(3162.000,5886.000)(3042.000,5916.000)(3042.000,5856.000)
\path(4062,5886)(4662,5886)
\blacken\path(4542.000,5856.000)(4662.000,5886.000)(4542.000,5916.000)(4542.000,5856.000)
\path(1662,4611)(1062,4311)
\blacken\path(1155.915,4391.498)(1062.000,4311.000)(1182.748,4337.833)(1155.915,4391.498)
\path(5112,4611)(5112,4311)
\blacken\path(5082.000,4431.000)(5112.000,4311.000)(5142.000,4431.000)(5082.000,4431.000)
\path(1662,2361)(1062,2061)
\blacken\path(1155.915,2141.498)(1062.000,2061.000)(1182.748,2087.833)(1155.915,2141.498)
\path(4062,2361)(4662,2061)
\blacken\path(4541.252,2087.833)(4662.000,2061.000)(4568.085,2141.498)(4541.252,2087.833)
\path(1062,1386)(1662,1386)
\blacken\path(1542.000,1356.000)(1662.000,1386.000)(1542.000,1416.000)(1542.000,1356.000)
\path(2562,1386)(3162,1386)
\blacken\path(3042.000,1356.000)(3162.000,1386.000)(3042.000,1416.000)(3042.000,1356.000)
\path(1062,3636)(1662,3636)
\blacken\path(1542.000,3606.000)(1662.000,3636.000)(1542.000,3666.000)(1542.000,3606.000)
\path(4662,3636)(4062,3636)
\blacken\path(4182.000,3666.000)(4062.000,3636.000)(4182.000,3606.000)(4182.000,3666.000)
\dashline{60.000}(1062,3336)(312,2436)
\dashline{60.000}(912,1236)(162,636)
\dashline{60.000}(2412,3336)(1662,2736)
\dashline{60.000}(1662,636)(2412,1236)
\dashline{60.000}(3237,786)(4212,336)
\dashline{60.000}(5262,2736)(5712,3336)
\dashline{60.000}(3762,2736)(4212,3336)
\dashline{60.000}(5412,636)(4812,1086)
\path(5712,111)(5787,111)(5787,36)
	(5712,36)(5712,111)
\path(5637,111)(5712,111)(5712,36)
	(5637,36)(5637,111)
\dashline{60.000}(3912,4986)(4362,5586)
\path(2037,6336)(2187,6336)
\path(3537,6336)(3687,6336)(3612,6186)(3537,6336)
\path(3462,6186)(3612,6186)
\path(5262,6186)(4962,6186)(5037,6036)(4887,6036)
\path(5037,6036)(5112,6186)
\path(5187,3786)(4887,3786)
\path(5262,1986)(5112,1986)(5337,1536)
	(5412,1686)(4962,1686)(5037,1536)(5262,1986)
\path(5337,1836)(5037,1836)(5187,1536)(5337,1836)
\path(2187,3786)(2037,4086)(2187,4086)
	(2037,3786)(1962,3936)(2262,3936)
	(2187,3786)(1887,3786)
\path(537,3786)(462,3936)(762,3936)
\path(387,3786)(537,3786)(687,4086)
	(537,4086)(612,3936)
\path(2187,6336)(2112,6186)(2037,6336)
\path(2112,6186)(1962,6186)
\path(3462,6186)(3537,6036)
\path(5262,6186)(5187,6036)(5037,6036)
\path(5112,6186)(5187,6036)
\path(5337,3786)(5187,3786)
\path(5037,3786)(4962,3636)(4887,3786)
\path(3537,3786)(3612,3636)(3687,3786)
\path(537,1686)(462,1836)(762,1836)
	(687,1686)(537,1986)(687,1986)
	(462,1536)(387,1686)(837,1686)
\path(2337,1686)(1887,1686)(1962,1536)
	(2187,1986)(2037,1986)(2187,1686)
	(2262,1836)(1962,1836)(2112,1536)(2187,1686)
\path(3762,1836)(3462,1836)(3612,1536)(3762,1836)
\path(3762,6186)(3612,6186)(3537,6036)(3387,6036)
\path(5037,4086)(5187,4086)(5037,3786)
	(4962,3936)(5262,3936)(5187,3786)(5037,4086)
\path(462,6861)(537,6861)(537,6786)
	(462,6786)(462,6861)
\path(162,5886)(462,6486)(762,5886)(162,5886)
\path(1812,5886)(2112,6486)(2412,5886)(1812,5886)
\path(3312,5886)(3612,6486)(3912,5886)(3312,5886)
\path(4812,5886)(5112,6486)(5412,5886)(4812,5886)
\path(4812,3636)(5112,4236)(5412,3636)(4812,3636)
\path(3312,3636)(3612,4236)(3912,3636)(3312,3636)
\path(4887,1536)(5187,2136)(5487,1536)(4887,1536)
\path(3312,1536)(3612,2136)(3912,1536)(3312,1536)
\path(312,1536)(612,2136)(912,1536)(312,1536)
\path(1812,1536)(2112,2136)(2412,1536)(1812,1536)
\path(1812,3636)(2112,4236)(2412,3636)(1812,3636)
\path(312,3636)(612,4236)(912,3636)(312,3636)
\path(3837,3786)(3387,3786)(3462,3636)
	(3687,4086)(3537,4086)(3687,3786)
	(3762,3936)(3462,3936)(3537,3786)
\path(3687,1986)(3537,1986)(3762,1536)
	(3837,1686)(3387,1686)(3462,1536)(3687,1986)
\path(5037,6336)(5187,6336)(5112,6186)(5037,6336)
\put(12,6786){\makebox(0,0)[lb]{{\SetFigFont{8}{9.6}{\rmdefault}{\mddefault}{\updefault}input:  }}}
\put(612,6786){\makebox(0,0)[lb]{{\SetFigFont{8}{9.6}{\rmdefault}{\mddefault}{\updefault}$\times$}}}
\put(3012,36){\makebox(0,0)[lb]{{\SetFigFont{8}{9.6}{\rmdefault}{\mddefault}{\updefault}output:  $0110 =$}}}
\put(4512,36){\makebox(0,0)[lb]{{\SetFigFont{8}{9.6}{\rmdefault}{\mddefault}{\updefault}output:  $1001=$}}}
\put(1137,3711){\makebox(0,0)[lb]{{\SetFigFont{8}{9.6}{\rmdefault}{\mddefault}{\updefault}$\dagger$}}}
\put(5187,4461){\makebox(0,0)[lb]{{\SetFigFont{8}{9.6}{\rmdefault}{\mddefault}{\updefault}$\dagger$}}}
\put(162,6636){\makebox(0,0)[lb]{{\SetFigFont{8}{9.6}{\rmdefault}{\mddefault}{\updefault}($\alpha=\beta=0101$ in string notation)}}}
\put(912,6786){\makebox(0,0)[lb]{{\SetFigFont{8}{9.6}{\rmdefault}{\mddefault}{\updefault}in partition notation}}}
\put(403,6385){\makebox(0,0)[lb]{{\SetFigFont{5}{6.0}{\rmdefault}{\mddefault}{\updefault}$1$}}}
\put(482,6385){\makebox(0,0)[lb]{{\SetFigFont{5}{6.0}{\rmdefault}{\mddefault}{\updefault}$0$}}}
\put(557,6242){\makebox(0,0)[lb]{{\SetFigFont{5}{6.0}{\rmdefault}{\mddefault}{\updefault}$1$}}}
\put(334,6241){\makebox(0,0)[lb]{{\SetFigFont{5}{6.0}{\rmdefault}{\mddefault}{\updefault}$0$}}}
\put(243,6079){\makebox(0,0)[lb]{{\SetFigFont{5}{6.0}{\rmdefault}{\mddefault}{\updefault}$1$}}}
\put(635,6082){\makebox(0,0)[lb]{{\SetFigFont{5}{6.0}{\rmdefault}{\mddefault}{\updefault}$0$}}}
\put(187,5942){\makebox(0,0)[lb]{{\SetFigFont{5}{6.0}{\rmdefault}{\mddefault}{\updefault}$0$}}}
\put(704,5955){\makebox(0,0)[lb]{{\SetFigFont{5}{6.0}{\rmdefault}{\mddefault}{\updefault}$1$}}}
\put(2135,6389){\makebox(0,0)[lb]{{\SetFigFont{5}{6.0}{\rmdefault}{\mddefault}{\updefault}$0$}}}
\put(1834,5948){\makebox(0,0)[lb]{{\SetFigFont{5}{6.0}{\rmdefault}{\mddefault}{\updefault}$0$}}}
\put(1903,6094){\makebox(0,0)[lb]{{\SetFigFont{5}{6.0}{\rmdefault}{\mddefault}{\updefault}$1$}}}
\put(2207,6236){\makebox(0,0)[lb]{{\SetFigFont{5}{6.0}{\rmdefault}{\mddefault}{\updefault}$1$}}}
\put(2292,6082){\makebox(0,0)[lb]{{\SetFigFont{5}{6.0}{\rmdefault}{\mddefault}{\updefault}$0$}}}
\put(2360,5951){\makebox(0,0)[lb]{{\SetFigFont{5}{6.0}{\rmdefault}{\mddefault}{\updefault}$1$}}}
\put(2137,6235){\makebox(0,0)[lb]{{\SetFigFont{5}{6.0}{\rmdefault}{\mddefault}{\updefault}$1$}}}
\put(2059,6232){\makebox(0,0)[lb]{{\SetFigFont{5}{6.0}{\rmdefault}{\mddefault}{\updefault}$0$}}}
\put(2049,6391){\makebox(0,0)[lb]{{\SetFigFont{5}{6.0}{\rmdefault}{\mddefault}{\updefault}$1$}}}
\put(1978,6239){\makebox(0,0)[lb]{{\SetFigFont{5}{6.0}{\rmdefault}{\mddefault}{\updefault}$0$}}}
\put(3331,5938){\makebox(0,0)[lb]{{\SetFigFont{5}{6.0}{\rmdefault}{\mddefault}{\updefault}$0$}}}
\put(3847,5941){\makebox(0,0)[lb]{{\SetFigFont{5}{6.0}{\rmdefault}{\mddefault}{\updefault}$1$}}}
\put(3448,6007){\makebox(0,0)[lb]{{\SetFigFont{5}{6.0}{\rmdefault}{\mddefault}{\updefault}$1$}}}
\put(3560,6079){\makebox(0,0)[lb]{{\SetFigFont{5}{6.0}{\rmdefault}{\mddefault}{\updefault}$T$}}}
\put(3484,6078){\makebox(0,0)[lb]{{\SetFigFont{5}{6.0}{\rmdefault}{\mddefault}{\updefault}$1$}}}
\put(3677,6163){\makebox(0,0)[lb]{{\SetFigFont{5}{6.0}{\rmdefault}{\mddefault}{\updefault}$1$}}}
\put(3632,6238){\makebox(0,0)[lb]{{\SetFigFont{5}{6.0}{\rmdefault}{\mddefault}{\updefault}$1$}}}
\put(3559,6235){\makebox(0,0)[lb]{{\SetFigFont{5}{6.0}{\rmdefault}{\mddefault}{\updefault}$0$}}}
\put(3399,6079){\makebox(0,0)[lb]{{\SetFigFont{5}{6.0}{\rmdefault}{\mddefault}{\updefault}$1$}}}
\put(3478,6242){\makebox(0,0)[lb]{{\SetFigFont{5}{6.0}{\rmdefault}{\mddefault}{\updefault}$0$}}}
\put(3556,6389){\makebox(0,0)[lb]{{\SetFigFont{5}{6.0}{\rmdefault}{\mddefault}{\updefault}$1$}}}
\put(3707,6238){\makebox(0,0)[lb]{{\SetFigFont{5}{6.0}{\rmdefault}{\mddefault}{\updefault}$1$}}}
\put(3789,6082){\makebox(0,0)[lb]{{\SetFigFont{5}{6.0}{\rmdefault}{\mddefault}{\updefault}$0$}}}
\put(4835,5941){\makebox(0,0)[lb]{{\SetFigFont{5}{6.0}{\rmdefault}{\mddefault}{\updefault}$0$}}}
\put(5354,5938){\makebox(0,0)[lb]{{\SetFigFont{5}{6.0}{\rmdefault}{\mddefault}{\updefault}$1$}}}
\put(5288,6073){\makebox(0,0)[lb]{{\SetFigFont{5}{6.0}{\rmdefault}{\mddefault}{\updefault}$0$}}}
\put(4952,6010){\makebox(0,0)[lb]{{\SetFigFont{5}{6.0}{\rmdefault}{\mddefault}{\updefault}$1$}}}
\put(5095,6007){\makebox(0,0)[lb]{{\SetFigFont{5}{6.0}{\rmdefault}{\mddefault}{\updefault}$0$}}}
\put(5210,6088){\makebox(0,0)[lb]{{\SetFigFont{5}{6.0}{\rmdefault}{\mddefault}{\updefault}$1$}}}
\put(5135,6082){\makebox(0,0)[lb]{{\SetFigFont{5}{6.0}{\rmdefault}{\mddefault}{\updefault}$1$}}}
\put(5053,6082){\makebox(0,0)[lb]{{\SetFigFont{5}{6.0}{\rmdefault}{\mddefault}{\updefault}$T$}}}
\put(4985,6079){\makebox(0,0)[lb]{{\SetFigFont{5}{6.0}{\rmdefault}{\mddefault}{\updefault}$1$}}}
\put(4899,6082){\makebox(0,0)[lb]{{\SetFigFont{5}{6.0}{\rmdefault}{\mddefault}{\updefault}$1$}}}
\put(5024,6163){\makebox(0,0)[lb]{{\SetFigFont{5}{6.0}{\rmdefault}{\mddefault}{\updefault}$0$}}}
\put(5211,6229){\makebox(0,0)[lb]{{\SetFigFont{5}{6.0}{\rmdefault}{\mddefault}{\updefault}$1$}}}
\put(5132,6235){\makebox(0,0)[lb]{{\SetFigFont{5}{6.0}{\rmdefault}{\mddefault}{\updefault}$1$}}}
\put(4974,6229){\makebox(0,0)[lb]{{\SetFigFont{5}{6.0}{\rmdefault}{\mddefault}{\updefault}$0$}}}
\put(5132,6389){\makebox(0,0)[lb]{{\SetFigFont{5}{6.0}{\rmdefault}{\mddefault}{\updefault}$0$}}}
\put(5046,6385){\makebox(0,0)[lb]{{\SetFigFont{5}{6.0}{\rmdefault}{\mddefault}{\updefault}$1$}}}
\put(5354,3691){\makebox(0,0)[lb]{{\SetFigFont{5}{6.0}{\rmdefault}{\mddefault}{\updefault}$1$}}}
\put(4874,3609){\makebox(0,0)[lb]{{\SetFigFont{5}{6.0}{\rmdefault}{\mddefault}{\updefault}$1$}}}
\put(4828,3688){\makebox(0,0)[lb]{{\SetFigFont{5}{6.0}{\rmdefault}{\mddefault}{\updefault}$0$}}}
\put(4985,3685){\makebox(0,0)[lb]{{\SetFigFont{5}{6.0}{\rmdefault}{\mddefault}{\updefault}$0$}}}
\put(4946,3760){\makebox(0,0)[lb]{{\SetFigFont{5}{6.0}{\rmdefault}{\mddefault}{\updefault}$1$}}}
\put(5099,3760){\makebox(0,0)[lb]{{\SetFigFont{5}{6.0}{\rmdefault}{\mddefault}{\updefault}$0$}}}
\put(5210,3838){\makebox(0,0)[lb]{{\SetFigFont{5}{6.0}{\rmdefault}{\mddefault}{\updefault}$1$}}}
\put(5282,3839){\makebox(0,0)[lb]{{\SetFigFont{5}{6.0}{\rmdefault}{\mddefault}{\updefault}$0$}}}
\put(5134,3838){\makebox(0,0)[lb]{{\SetFigFont{5}{6.0}{\rmdefault}{\mddefault}{\updefault}$1$}}}
\put(5056,3838){\makebox(0,0)[lb]{{\SetFigFont{5}{6.0}{\rmdefault}{\mddefault}{\updefault}$T$}}}
\put(4982,3835){\makebox(0,0)[lb]{{\SetFigFont{5}{6.0}{\rmdefault}{\mddefault}{\updefault}$1$}}}
\put(4903,3841){\makebox(0,0)[lb]{{\SetFigFont{5}{6.0}{\rmdefault}{\mddefault}{\updefault}$1$}}}
\put(5020,3916){\makebox(0,0)[lb]{{\SetFigFont{5}{6.0}{\rmdefault}{\mddefault}{\updefault}$0$}}}
\put(5176,3917){\makebox(0,0)[lb]{{\SetFigFont{5}{6.0}{\rmdefault}{\mddefault}{\updefault}$1$}}}
\put(5213,3989){\makebox(0,0)[lb]{{\SetFigFont{5}{6.0}{\rmdefault}{\mddefault}{\updefault}$1$}}}
\put(5128,3988){\makebox(0,0)[lb]{{\SetFigFont{5}{6.0}{\rmdefault}{\mddefault}{\updefault}$1$}}}
\put(5063,3985){\makebox(0,0)[lb]{{\SetFigFont{5}{6.0}{\rmdefault}{\mddefault}{\updefault}$0$}}}
\put(3857,3689){\makebox(0,0)[lb]{{\SetFigFont{5}{6.0}{\rmdefault}{\mddefault}{\updefault}$1$}}}
\put(3520,3606){\makebox(0,0)[lb]{{\SetFigFont{5}{6.0}{\rmdefault}{\mddefault}{\updefault}$0$}}}
\put(3635,3685){\makebox(0,0)[lb]{{\SetFigFont{5}{6.0}{\rmdefault}{\mddefault}{\updefault}$0$}}}
\put(3485,3681){\makebox(0,0)[lb]{{\SetFigFont{5}{6.0}{\rmdefault}{\mddefault}{\updefault}$0$}}}
\put(3327,3688){\makebox(0,0)[lb]{{\SetFigFont{5}{6.0}{\rmdefault}{\mddefault}{\updefault}$0$}}}
\put(3446,3764){\makebox(0,0)[lb]{{\SetFigFont{5}{6.0}{\rmdefault}{\mddefault}{\updefault}$1$}}}
\put(3605,3766){\makebox(0,0)[lb]{{\SetFigFont{5}{6.0}{\rmdefault}{\mddefault}{\updefault}$0$}}}
\put(3482,3835){\makebox(0,0)[lb]{{\SetFigFont{5}{6.0}{\rmdefault}{\mddefault}{\updefault}$1$}}}
\put(3628,3848){\makebox(0,0)[lb]{{\SetFigFont{5}{6.0}{\rmdefault}{\mddefault}{\updefault}$1$}}}
\put(3703,3844){\makebox(0,0)[lb]{{\SetFigFont{5}{6.0}{\rmdefault}{\mddefault}{\updefault}$1$}}}
\put(3778,3842){\makebox(0,0)[lb]{{\SetFigFont{5}{6.0}{\rmdefault}{\mddefault}{\updefault}$0$}}}
\put(3674,3914){\makebox(0,0)[lb]{{\SetFigFont{5}{6.0}{\rmdefault}{\mddefault}{\updefault}$1$}}}
\put(3707,3994){\makebox(0,0)[lb]{{\SetFigFont{5}{6.0}{\rmdefault}{\mddefault}{\updefault}$1$}}}
\put(3647,3995){\makebox(0,0)[lb]{{\SetFigFont{5}{6.0}{\rmdefault}{\mddefault}{\updefault}$1$}}}
\put(3527,3916){\makebox(0,0)[lb]{{\SetFigFont{5}{6.0}{\rmdefault}{\mddefault}{\updefault}$0$}}}
\put(3396,3838){\makebox(0,0)[lb]{{\SetFigFont{5}{6.0}{\rmdefault}{\mddefault}{\updefault}$1$}}}
\put(3478,3992){\makebox(0,0)[lb]{{\SetFigFont{5}{6.0}{\rmdefault}{\mddefault}{\updefault}$0$}}}
\put(3635,4138){\makebox(0,0)[lb]{{\SetFigFont{5}{6.0}{\rmdefault}{\mddefault}{\updefault}$0$}}}
\put(3552,4135){\makebox(0,0)[lb]{{\SetFigFont{5}{6.0}{\rmdefault}{\mddefault}{\updefault}$1$}}}
\put(1831,3701){\makebox(0,0)[lb]{{\SetFigFont{5}{6.0}{\rmdefault}{\mddefault}{\updefault}$0$}}}
\put(2350,3708){\makebox(0,0)[lb]{{\SetFigFont{5}{6.0}{\rmdefault}{\mddefault}{\updefault}$1$}}}
\put(2210,3832){\makebox(0,0)[lb]{{\SetFigFont{5}{6.0}{\rmdefault}{\mddefault}{\updefault}$0$}}}
\put(2059,3838){\makebox(0,0)[lb]{{\SetFigFont{5}{6.0}{\rmdefault}{\mddefault}{\updefault}$0$}}}
\put(1981,3835){\makebox(0,0)[lb]{{\SetFigFont{5}{6.0}{\rmdefault}{\mddefault}{\updefault}$0$}}}
\put(1906,3838){\makebox(0,0)[lb]{{\SetFigFont{5}{6.0}{\rmdefault}{\mddefault}{\updefault}$1$}}}
\put(2286,3832){\makebox(0,0)[lb]{{\SetFigFont{5}{6.0}{\rmdefault}{\mddefault}{\updefault}$0$}}}
\put(1978,3995){\makebox(0,0)[lb]{{\SetFigFont{5}{6.0}{\rmdefault}{\mddefault}{\updefault}$0$}}}
\put(2144,3992){\makebox(0,0)[lb]{{\SetFigFont{5}{6.0}{\rmdefault}{\mddefault}{\updefault}$1$}}}
\put(2213,3992){\makebox(0,0)[lb]{{\SetFigFont{5}{6.0}{\rmdefault}{\mddefault}{\updefault}$1$}}}
\put(2053,3995){\makebox(0,0)[lb]{{\SetFigFont{5}{6.0}{\rmdefault}{\mddefault}{\updefault}$0$}}}
\put(2132,4138){\makebox(0,0)[lb]{{\SetFigFont{5}{6.0}{\rmdefault}{\mddefault}{\updefault}$0$}}}
\put(2049,4135){\makebox(0,0)[lb]{{\SetFigFont{5}{6.0}{\rmdefault}{\mddefault}{\updefault}$1$}}}
\put(857,3698){\makebox(0,0)[lb]{{\SetFigFont{5}{6.0}{\rmdefault}{\mddefault}{\updefault}$1$}}}
\put(785,3832){\makebox(0,0)[lb]{{\SetFigFont{5}{6.0}{\rmdefault}{\mddefault}{\updefault}$0$}}}
\put(331,3695){\makebox(0,0)[lb]{{\SetFigFont{5}{6.0}{\rmdefault}{\mddefault}{\updefault}$0$}}}
\put(563,3838){\makebox(0,0)[lb]{{\SetFigFont{5}{6.0}{\rmdefault}{\mddefault}{\updefault}$0$}}}
\put(485,3832){\makebox(0,0)[lb]{{\SetFigFont{5}{6.0}{\rmdefault}{\mddefault}{\updefault}$0$}}}
\put(403,3838){\makebox(0,0)[lb]{{\SetFigFont{5}{6.0}{\rmdefault}{\mddefault}{\updefault}$1$}}}
\put(632,4138){\makebox(0,0)[lb]{{\SetFigFont{5}{6.0}{\rmdefault}{\mddefault}{\updefault}$0$}}}
\put(560,3995){\makebox(0,0)[lb]{{\SetFigFont{5}{6.0}{\rmdefault}{\mddefault}{\updefault}$0$}}}
\put(625,3991){\makebox(0,0)[lb]{{\SetFigFont{5}{6.0}{\rmdefault}{\mddefault}{\updefault}$1$}}}
\put(701,3994){\makebox(0,0)[lb]{{\SetFigFont{5}{6.0}{\rmdefault}{\mddefault}{\updefault}$1$}}}
\put(484,3992){\makebox(0,0)[lb]{{\SetFigFont{5}{6.0}{\rmdefault}{\mddefault}{\updefault}$0$}}}
\put(671,3913){\makebox(0,0)[lb]{{\SetFigFont{5}{6.0}{\rmdefault}{\mddefault}{\updefault}$1$}}}
\put(331,1586){\makebox(0,0)[lb]{{\SetFigFont{5}{6.0}{\rmdefault}{\mddefault}{\updefault}$0$}}}
\put(488,1582){\makebox(0,0)[lb]{{\SetFigFont{5}{6.0}{\rmdefault}{\mddefault}{\updefault}$1$}}}
\put(595,1669){\makebox(0,0)[lb]{{\SetFigFont{5}{6.0}{\rmdefault}{\mddefault}{\updefault}$1$}}}
\put(739,1666){\makebox(0,0)[lb]{{\SetFigFont{5}{6.0}{\rmdefault}{\mddefault}{\updefault}$0$}}}
\put(854,1595){\makebox(0,0)[lb]{{\SetFigFont{5}{6.0}{\rmdefault}{\mddefault}{\updefault}$1$}}}
\put(781,1738){\makebox(0,0)[lb]{{\SetFigFont{5}{6.0}{\rmdefault}{\mddefault}{\updefault}$0$}}}
\put(709,1741){\makebox(0,0)[lb]{{\SetFigFont{5}{6.0}{\rmdefault}{\mddefault}{\updefault}$0$}}}
\put(560,1741){\makebox(0,0)[lb]{{\SetFigFont{5}{6.0}{\rmdefault}{\mddefault}{\updefault}$0$}}}
\put(400,1741){\makebox(0,0)[lb]{{\SetFigFont{5}{6.0}{\rmdefault}{\mddefault}{\updefault}$1$}}}
\put(484,1895){\makebox(0,0)[lb]{{\SetFigFont{5}{6.0}{\rmdefault}{\mddefault}{\updefault}$0$}}}
\put(625,1894){\makebox(0,0)[lb]{{\SetFigFont{5}{6.0}{\rmdefault}{\mddefault}{\updefault}$1$}}}
\put(700,1894){\makebox(0,0)[lb]{{\SetFigFont{5}{6.0}{\rmdefault}{\mddefault}{\updefault}$1$}}}
\put(553,2042){\makebox(0,0)[lb]{{\SetFigFont{5}{6.0}{\rmdefault}{\mddefault}{\updefault}$1$}}}
\put(629,2042){\makebox(0,0)[lb]{{\SetFigFont{5}{6.0}{\rmdefault}{\mddefault}{\updefault}$0$}}}
\put(2021,1513){\makebox(0,0)[lb]{{\SetFigFont{5}{6.0}{\rmdefault}{\mddefault}{\updefault}$1$}}}
\put(1834,1588){\makebox(0,0)[lb]{{\SetFigFont{5}{6.0}{\rmdefault}{\mddefault}{\updefault}$0$}}}
\put(1978,1585){\makebox(0,0)[lb]{{\SetFigFont{5}{6.0}{\rmdefault}{\mddefault}{\updefault}$1$}}}
\put(2050,1588){\makebox(0,0)[lb]{{\SetFigFont{5}{6.0}{\rmdefault}{\mddefault}{\updefault}$1$}}}
\put(2125,1595){\makebox(0,0)[lb]{{\SetFigFont{5}{6.0}{\rmdefault}{\mddefault}{\updefault}$1$}}}
\put(2350,1598){\makebox(0,0)[lb]{{\SetFigFont{5}{6.0}{\rmdefault}{\mddefault}{\updefault}$1$}}}
\put(1909,1738){\makebox(0,0)[lb]{{\SetFigFont{5}{6.0}{\rmdefault}{\mddefault}{\updefault}$1$}}}
\put(1982,1742){\makebox(0,0)[lb]{{\SetFigFont{5}{6.0}{\rmdefault}{\mddefault}{\updefault}$0$}}}
\put(2050,1742){\makebox(0,0)[lb]{{\SetFigFont{5}{6.0}{\rmdefault}{\mddefault}{\updefault}$0$}}}
\put(2125,1741){\makebox(0,0)[lb]{{\SetFigFont{5}{6.0}{\rmdefault}{\mddefault}{\updefault}$T$}}}
\put(2200,1742){\makebox(0,0)[lb]{{\SetFigFont{5}{6.0}{\rmdefault}{\mddefault}{\updefault}$0$}}}
\put(2281,1741){\makebox(0,0)[lb]{{\SetFigFont{5}{6.0}{\rmdefault}{\mddefault}{\updefault}$0$}}}
\put(2177,1817){\makebox(0,0)[lb]{{\SetFigFont{5}{6.0}{\rmdefault}{\mddefault}{\updefault}$1$}}}
\put(2024,1814){\makebox(0,0)[lb]{{\SetFigFont{5}{6.0}{\rmdefault}{\mddefault}{\updefault}$0$}}}
\put(1981,1892){\makebox(0,0)[lb]{{\SetFigFont{5}{6.0}{\rmdefault}{\mddefault}{\updefault}$0$}}}
\put(2132,1895){\makebox(0,0)[lb]{{\SetFigFont{5}{6.0}{\rmdefault}{\mddefault}{\updefault}$1$}}}
\put(2200,1895){\makebox(0,0)[lb]{{\SetFigFont{5}{6.0}{\rmdefault}{\mddefault}{\updefault}$1$}}}
\put(2053,1892){\makebox(0,0)[lb]{{\SetFigFont{5}{6.0}{\rmdefault}{\mddefault}{\updefault}$0$}}}
\put(2135,2039){\makebox(0,0)[lb]{{\SetFigFont{5}{6.0}{\rmdefault}{\mddefault}{\updefault}$0$}}}
\put(2050,2039){\makebox(0,0)[lb]{{\SetFigFont{5}{6.0}{\rmdefault}{\mddefault}{\updefault}$1$}}}
\put(3373,1510){\makebox(0,0)[lb]{{\SetFigFont{5}{6.0}{\rmdefault}{\mddefault}{\updefault}$0$}}}
\put(3520,1506){\makebox(0,0)[lb]{{\SetFigFont{5}{6.0}{\rmdefault}{\mddefault}{\updefault}$1$}}}
\put(3667,1507){\makebox(0,0)[lb]{{\SetFigFont{5}{6.0}{\rmdefault}{\mddefault}{\updefault}$1$}}}
\put(3820,1512){\makebox(0,0)[lb]{{\SetFigFont{5}{6.0}{\rmdefault}{\mddefault}{\updefault}$0$}}}
\put(3327,1585){\makebox(0,0)[lb]{{\SetFigFont{5}{6.0}{\rmdefault}{\mddefault}{\updefault}$0$}}}
\put(3403,1585){\makebox(0,0)[lb]{{\SetFigFont{5}{6.0}{\rmdefault}{\mddefault}{\updefault}$0$}}}
\put(3478,1592){\makebox(0,0)[lb]{{\SetFigFont{5}{6.0}{\rmdefault}{\mddefault}{\updefault}$1$}}}
\put(3778,1589){\makebox(0,0)[lb]{{\SetFigFont{5}{6.0}{\rmdefault}{\mddefault}{\updefault}$T$}}}
\put(3857,1595){\makebox(0,0)[lb]{{\SetFigFont{5}{6.0}{\rmdefault}{\mddefault}{\updefault}$1$}}}
\put(3706,1585){\makebox(0,0)[lb]{{\SetFigFont{5}{6.0}{\rmdefault}{\mddefault}{\updefault}$1$}}}
\put(3625,1585){\makebox(0,0)[lb]{{\SetFigFont{5}{6.0}{\rmdefault}{\mddefault}{\updefault}$1$}}}
\put(3556,1592){\makebox(0,0)[lb]{{\SetFigFont{5}{6.0}{\rmdefault}{\mddefault}{\updefault}$1$}}}
\put(3742,1666){\makebox(0,0)[lb]{{\SetFigFont{5}{6.0}{\rmdefault}{\mddefault}{\updefault}$0$}}}
\put(3791,1735){\makebox(0,0)[lb]{{\SetFigFont{5}{6.0}{\rmdefault}{\mddefault}{\updefault}$0$}}}
\put(3703,1742){\makebox(0,0)[lb]{{\SetFigFont{5}{6.0}{\rmdefault}{\mddefault}{\updefault}$0$}}}
\put(3557,1738){\makebox(0,0)[lb]{{\SetFigFont{5}{6.0}{\rmdefault}{\mddefault}{\updefault}$0$}}}
\put(3484,1738){\makebox(0,0)[lb]{{\SetFigFont{5}{6.0}{\rmdefault}{\mddefault}{\updefault}$0$}}}
\put(3406,1739){\makebox(0,0)[lb]{{\SetFigFont{5}{6.0}{\rmdefault}{\mddefault}{\updefault}$1$}}}
\put(3706,1886){\makebox(0,0)[lb]{{\SetFigFont{5}{6.0}{\rmdefault}{\mddefault}{\updefault}$1$}}}
\put(3628,1888){\makebox(0,0)[lb]{{\SetFigFont{5}{6.0}{\rmdefault}{\mddefault}{\updefault}$1$}}}
\put(3559,1885){\makebox(0,0)[lb]{{\SetFigFont{5}{6.0}{\rmdefault}{\mddefault}{\updefault}$0$}}}
\put(3481,1885){\makebox(0,0)[lb]{{\SetFigFont{5}{6.0}{\rmdefault}{\mddefault}{\updefault}$0$}}}
\put(3556,2041){\makebox(0,0)[lb]{{\SetFigFont{5}{6.0}{\rmdefault}{\mddefault}{\updefault}$1$}}}
\put(3632,2041){\makebox(0,0)[lb]{{\SetFigFont{5}{6.0}{\rmdefault}{\mddefault}{\updefault}$0$}}}
\put(4942,1513){\makebox(0,0)[lb]{{\SetFigFont{5}{6.0}{\rmdefault}{\mddefault}{\updefault}$1$}}}
\put(5095,1513){\makebox(0,0)[lb]{{\SetFigFont{5}{6.0}{\rmdefault}{\mddefault}{\updefault}$0$}}}
\put(5249,1509){\makebox(0,0)[lb]{{\SetFigFont{5}{6.0}{\rmdefault}{\mddefault}{\updefault}$0$}}}
\put(5392,1513){\makebox(0,0)[lb]{{\SetFigFont{5}{6.0}{\rmdefault}{\mddefault}{\updefault}$1$}}}
\put(4906,1589){\makebox(0,0)[lb]{{\SetFigFont{5}{6.0}{\rmdefault}{\mddefault}{\updefault}$0$}}}
\put(5056,1589){\makebox(0,0)[lb]{{\SetFigFont{5}{6.0}{\rmdefault}{\mddefault}{\updefault}$0$}}}
\put(5128,1592){\makebox(0,0)[lb]{{\SetFigFont{5}{6.0}{\rmdefault}{\mddefault}{\updefault}$0$}}}
\put(5200,1592){\makebox(0,0)[lb]{{\SetFigFont{5}{6.0}{\rmdefault}{\mddefault}{\updefault}$0$}}}
\put(5272,1594){\makebox(0,0)[lb]{{\SetFigFont{5}{6.0}{\rmdefault}{\mddefault}{\updefault}$0$}}}
\put(5346,1592){\makebox(0,0)[lb]{{\SetFigFont{5}{6.0}{\rmdefault}{\mddefault}{\updefault}$1$}}}
\put(5357,1738){\makebox(0,0)[lb]{{\SetFigFont{5}{6.0}{\rmdefault}{\mddefault}{\updefault}$0$}}}
\put(5291,1739){\makebox(0,0)[lb]{{\SetFigFont{5}{6.0}{\rmdefault}{\mddefault}{\updefault}$1$}}}
\put(5212,1742){\makebox(0,0)[lb]{{\SetFigFont{5}{6.0}{\rmdefault}{\mddefault}{\updefault}$1$}}}
\put(5056,1738){\makebox(0,0)[lb]{{\SetFigFont{5}{6.0}{\rmdefault}{\mddefault}{\updefault}$1$}}}
\put(4978,1741){\makebox(0,0)[lb]{{\SetFigFont{5}{6.0}{\rmdefault}{\mddefault}{\updefault}$1$}}}
\put(5053,1895){\makebox(0,0)[lb]{{\SetFigFont{5}{6.0}{\rmdefault}{\mddefault}{\updefault}$0$}}}
\put(5128,1891){\makebox(0,0)[lb]{{\SetFigFont{5}{6.0}{\rmdefault}{\mddefault}{\updefault}$0$}}}
\put(5200,1892){\makebox(0,0)[lb]{{\SetFigFont{5}{6.0}{\rmdefault}{\mddefault}{\updefault}$1$}}}
\put(5272,1891){\makebox(0,0)[lb]{{\SetFigFont{5}{6.0}{\rmdefault}{\mddefault}{\updefault}$1$}}}
\put(5131,2045){\makebox(0,0)[lb]{{\SetFigFont{5}{6.0}{\rmdefault}{\mddefault}{\updefault}$1$}}}
\put(3625,6392){\makebox(0,0)[lb]{{\SetFigFont{5}{6.0}{\rmdefault}{\mddefault}{\updefault}$0$}}}
\put(3588,6321){\makebox(0,0)[lb]{{\SetFigFont{5}{6.0}{\rmdefault}{\mddefault}{\updefault}$T$}}}
\put(2088,6314){\makebox(0,0)[lb]{{\SetFigFont{5}{6.0}{\rmdefault}{\mddefault}{\updefault}$T$}}}
\put(2015,6160){\makebox(0,0)[lb]{{\SetFigFont{5}{6.0}{\rmdefault}{\mddefault}{\updefault}$0$}}}
\put(3511,6166){\makebox(0,0)[lb]{{\SetFigFont{5}{6.0}{\rmdefault}{\mddefault}{\updefault}$0$}}}
\put(5083,6316){\makebox(0,0)[lb]{{\SetFigFont{5}{6.0}{\rmdefault}{\mddefault}{\updefault}$T$}}}
\put(5056,6235){\makebox(0,0)[lb]{{\SetFigFont{5}{6.0}{\rmdefault}{\mddefault}{\updefault}$0$}}}
\put(5170,6160){\makebox(0,0)[lb]{{\SetFigFont{5}{6.0}{\rmdefault}{\mddefault}{\updefault}$1$}}}
\put(5128,4138){\makebox(0,0)[lb]{{\SetFigFont{5}{6.0}{\rmdefault}{\mddefault}{\updefault}$0$}}}
\put(5057,4138){\makebox(0,0)[lb]{{\SetFigFont{5}{6.0}{\rmdefault}{\mddefault}{\updefault}$1$}}}
\put(5082,4067){\makebox(0,0)[lb]{{\SetFigFont{5}{6.0}{\rmdefault}{\mddefault}{\updefault}$T$}}}
\put(4974,3988){\makebox(0,0)[lb]{{\SetFigFont{5}{6.0}{\rmdefault}{\mddefault}{\updefault}$0$}}}
\put(5232,3760){\makebox(0,0)[lb]{{\SetFigFont{5}{6.0}{\rmdefault}{\mddefault}{\updefault}$T$}}}
\put(3730,3766){\makebox(0,0)[lb]{{\SetFigFont{5}{6.0}{\rmdefault}{\mddefault}{\updefault}$T$}}}
\put(4896,3688){\makebox(0,0)[lb]{{\SetFigFont{5}{6.0}{\rmdefault}{\mddefault}{\updefault}$T$}}}
\put(3582,4069){\makebox(0,0)[lb]{{\SetFigFont{5}{6.0}{\rmdefault}{\mddefault}{\updefault}$T$}}}
\put(3553,3991){\makebox(0,0)[lb]{{\SetFigFont{5}{6.0}{\rmdefault}{\mddefault}{\updefault}$0$}}}
\put(3557,3841){\makebox(0,0)[lb]{{\SetFigFont{5}{6.0}{\rmdefault}{\mddefault}{\updefault}$T$}}}
\put(3400,3685){\makebox(0,0)[lb]{{\SetFigFont{5}{6.0}{\rmdefault}{\mddefault}{\updefault}$T$}}}
\put(3360,3609){\makebox(0,0)[lb]{{\SetFigFont{5}{6.0}{\rmdefault}{\mddefault}{\updefault}$1$}}}
\put(2082,4070){\makebox(0,0)[lb]{{\SetFigFont{5}{6.0}{\rmdefault}{\mddefault}{\updefault}$T$}}}
\put(2018,3917){\makebox(0,0)[lb]{{\SetFigFont{5}{6.0}{\rmdefault}{\mddefault}{\updefault}$0$}}}
\put(2167,3913){\makebox(0,0)[lb]{{\SetFigFont{5}{6.0}{\rmdefault}{\mddefault}{\updefault}$1$}}}
\put(2119,3832){\makebox(0,0)[lb]{{\SetFigFont{5}{6.0}{\rmdefault}{\mddefault}{\updefault}$T$}}}
\put(1929,3764){\makebox(0,0)[lb]{{\SetFigFont{5}{6.0}{\rmdefault}{\mddefault}{\updefault}$T$}}}
\put(2098,3763){\makebox(0,0)[lb]{{\SetFigFont{5}{6.0}{\rmdefault}{\mddefault}{\updefault}$1$}}}
\put(550,4139){\makebox(0,0)[lb]{{\SetFigFont{5}{6.0}{\rmdefault}{\mddefault}{\updefault}$1$}}}
\put(583,4067){\makebox(0,0)[lb]{{\SetFigFont{5}{6.0}{\rmdefault}{\mddefault}{\updefault}$T$}}}
\put(514,3917){\makebox(0,0)[lb]{{\SetFigFont{5}{6.0}{\rmdefault}{\mddefault}{\updefault}$0$}}}
\put(430,3763){\makebox(0,0)[lb]{{\SetFigFont{5}{6.0}{\rmdefault}{\mddefault}{\updefault}$T$}}}
\put(3553,3682){\makebox(0,0)[lb]{{\SetFigFont{5}{6.0}{\rmdefault}{\mddefault}{\updefault}$0$}}}
\put(582,1969){\makebox(0,0)[lb]{{\SetFigFont{5}{6.0}{\rmdefault}{\mddefault}{\updefault}$T$}}}
\put(549,1898){\makebox(0,0)[lb]{{\SetFigFont{5}{6.0}{\rmdefault}{\mddefault}{\updefault}$0$}}}
\put(514,1817){\makebox(0,0)[lb]{{\SetFigFont{5}{6.0}{\rmdefault}{\mddefault}{\updefault}$0$}}}
\put(668,1817){\makebox(0,0)[lb]{{\SetFigFont{5}{6.0}{\rmdefault}{\mddefault}{\updefault}$1$}}}
\put(472,1745){\makebox(0,0)[lb]{{\SetFigFont{5}{6.0}{\rmdefault}{\mddefault}{\updefault}$0$}}}
\put(618,1741){\makebox(0,0)[lb]{{\SetFigFont{5}{6.0}{\rmdefault}{\mddefault}{\updefault}$T$}}}
\put(433,1670){\makebox(0,0)[lb]{{\SetFigFont{5}{6.0}{\rmdefault}{\mddefault}{\updefault}$T$}}}
\put(396,1585){\makebox(0,0)[lb]{{\SetFigFont{5}{6.0}{\rmdefault}{\mddefault}{\updefault}$0$}}}
\put(361,1510){\makebox(0,0)[lb]{{\SetFigFont{5}{6.0}{\rmdefault}{\mddefault}{\updefault}$0$}}}
\put(2079,1964){\makebox(0,0)[lb]{{\SetFigFont{5}{6.0}{\rmdefault}{\mddefault}{\updefault}$T$}}}
\put(1936,1663){\makebox(0,0)[lb]{{\SetFigFont{5}{6.0}{\rmdefault}{\mddefault}{\updefault}$T$}}}
\put(2088,1666){\makebox(0,0)[lb]{{\SetFigFont{5}{6.0}{\rmdefault}{\mddefault}{\updefault}$1$}}}
\put(2231,1666){\makebox(0,0)[lb]{{\SetFigFont{5}{6.0}{\rmdefault}{\mddefault}{\updefault}$0$}}}
\put(1897,1589){\makebox(0,0)[lb]{{\SetFigFont{5}{6.0}{\rmdefault}{\mddefault}{\updefault}$0$}}}
\put(1857,1513){\makebox(0,0)[lb]{{\SetFigFont{5}{6.0}{\rmdefault}{\mddefault}{\updefault}$0$}}}
\put(3583,1963){\makebox(0,0)[lb]{{\SetFigFont{5}{6.0}{\rmdefault}{\mddefault}{\updefault}$T$}}}
\put(3657,1816){\makebox(0,0)[lb]{{\SetFigFont{5}{6.0}{\rmdefault}{\mddefault}{\updefault}$1$}}}
\put(3514,1817){\makebox(0,0)[lb]{{\SetFigFont{5}{6.0}{\rmdefault}{\mddefault}{\updefault}$0$}}}
\put(3616,1741){\makebox(0,0)[lb]{{\SetFigFont{5}{6.0}{\rmdefault}{\mddefault}{\updefault}$T$}}}
\put(3432,1666){\makebox(0,0)[lb]{{\SetFigFont{5}{6.0}{\rmdefault}{\mddefault}{\updefault}$T$}}}
\put(3592,1667){\makebox(0,0)[lb]{{\SetFigFont{5}{6.0}{\rmdefault}{\mddefault}{\updefault}$1$}}}
\put(5195,2045){\makebox(0,0)[lb]{{\SetFigFont{5}{6.0}{\rmdefault}{\mddefault}{\updefault}$0$}}}
\put(5164,1964){\makebox(0,0)[lb]{{\SetFigFont{5}{6.0}{\rmdefault}{\mddefault}{\updefault}$T$}}}
\put(5239,1816){\makebox(0,0)[lb]{{\SetFigFont{5}{6.0}{\rmdefault}{\mddefault}{\updefault}$1$}}}
\put(5088,1820){\makebox(0,0)[lb]{{\SetFigFont{5}{6.0}{\rmdefault}{\mddefault}{\updefault}$0$}}}
\put(5126,1742){\makebox(0,0)[lb]{{\SetFigFont{5}{6.0}{\rmdefault}{\mddefault}{\updefault}$T$}}}
\put(5014,1667){\makebox(0,0)[lb]{{\SetFigFont{5}{6.0}{\rmdefault}{\mddefault}{\updefault}$1$}}}
\put(5168,1663){\makebox(0,0)[lb]{{\SetFigFont{5}{6.0}{\rmdefault}{\mddefault}{\updefault}$0$}}}
\put(5315,1663){\makebox(0,0)[lb]{{\SetFigFont{5}{6.0}{\rmdefault}{\mddefault}{\updefault}$T$}}}
\put(4973,1586){\makebox(0,0)[lb]{{\SetFigFont{5}{6.0}{\rmdefault}{\mddefault}{\updefault}$T$}}}
\put(5262,5286){\blacken\ellipse{40}{40}}
\put(5262,5286){\ellipse{40}{40}}
\end{picture}
}

\end{center}
\caption{A motivating example of the degenerations
of the geometric Littlewood-Richardson rule,
illustrating $\si_{1}^2 = \si_{1,1} + \si_{2}$,
or in string notation,
$\si_{0101}^2 = \si_{0110} + \si_{1001}$ (equ.\ \eqref{basic}).  The puzzles and daggers will be explained shortly.  In the puzzles, $T$ is written instead of $10$ for convenience. \lremind{keyexample}
}
\label{keyexample}
\end{figure}

We continue the example, using this degeneration to deform an
intersection of two Schubert cycles.  (The reader may recognize this
as a variant of Example~\ref{fourlines}.)  This example, understood
well enough, leads to the general rule for the Grassmannian in
cohomology, and the many other rules (conjectural and otherwise)
stated in Part~\ref{partone}.  Let $\al = \be$ be the unique Schubert
divisor class (the class $\Box = (1)$ in terms of partitions), so
$\overline{\Omega}_{\al} = \overline{\Omega}_{\be}$ both correspond to
the set of lines in $\proj^3$ meeting a fixed line.  We will degenerate
the locus of lines meeting two general (skew) lines into a union of
Schubert varieties.  This is depicted in Figure~\ref{keyexample}.

Consider the locus of lines meeting both the fixed line $\proj F_2$
and the moving line $\proj M_2$.  This is clearly a two-dimensional
cycle (isomorphic to $\proj F_2 \times \proj M_2 \cong \proj^1 \times
\proj^1$).  After the first degeneration, the moving line has not
moved, so the locus remains the same.  In the second degeneration, the
moving line has degenerated to meet the fixed line at a point $\proj
F_1$.  Then there are two irreducible families of lines meeting both
the fixed and moving lines.  First, the lines could contain the point
$\proj F_1$ (depicted in the upper right panel).  Second, the lines
could lie in the plane spanned by $\proj F_2$ and $\proj M_2$ (the
center left panel).  In the first case, we continue the degenerations,
always interpreting our locus as those lines containing the point
$\proj F_1$.  At the end, we interpret this as the Schubert variety
$\overline{\Omega}_{(2)}(M_{\bullet}) = \overline{\Omega}_{(2)}(F_{\bullet})$.
After the next degeneration in the second case, the locus of lines in
question can be restated as those lines contained in the plane $\proj
M_3$; and this description continues until then end, where we
interpret this as the Schubert variety
$\overline{\Omega}_{(1,1)}(M_{\bullet}) = \overline{\Omega}_{(1,1)}(F_{\bullet})$.
We have thus illustrated the identity\lremind{basic}
\begin{equation}
\label{basic}
\si_{1} \cdot
\si_{1} =
\si_{2} +
\si_{1,1} .
\end{equation}

\subsection{The degenerations in general}
\label{degenerationorder}\lremind{degenerationorder}
In general, the degeneration corresponds to a path in the
Bruhat order of $S_n$.  The path corresponds to partial factorizations
from the right of the longest word $w_0$:
$$
w_0 = e_{n-1} \cdots e_2 e_1  \quad \cdots \quad e_{n-1} e_{n-2} e_{n-3} 
\quad e_{n-1} e_{n-2} \quad e_{n-1}.
$$
We denote the $\binom n 2 + 1$ permutations by
$d_0 = w_0$, \dots, 
 $d_{\binom n 2 - 2} = e_{n-2} e_1$,
$d_{\binom n 2 - 1} = e_{n-1}$,
$d_{\binom n 2} = e$.

At each stage, we will consider cycles of the form
$\overline{\Omega_{\alpha} (M_{\bullet}) \cap
  \Omega_{\beta}(F_{\bullet})}$.  At the start, where $M_{\bullet}$ and
$F_{\bullet}$ are transverse, we have the intersection of Schubert
cycles that we seek to degenerate.  At the end, where $M_{\bullet} =
F_{\bullet}$, we have Schubert cycles with respect to a single flag.
At each stage, this cycle $\overline{\Omega_{\alpha} (M_{\bullet}) \cap
  \Omega_{\beta}(F_{\bullet})}$ (with $M_{\bullet}$ related to
$F_{\bullet}$ by $d_i$, $i< \binom n 2$) will degenerate to one or two
cycles of the same form, except that $M_{\bullet}$ and $F_{\bullet}$
are related by $d_{i+1}$.  These one or two cycles will each appear
with multiplicity $1$.

In \cite{vakil:checkers}, these intermediate cycles are
described geometrically cleanly using {\em checkerboards}.  For the
generalizations we will consider, we will describe them here in terms
of Knutson and Tao's {\em puzzles} \cite{knutson:puzzle}, which we
describe now.  The bijection between checkerboard games and puzzles is
given in \cite[App.~ A]{vakil:checkers}.  

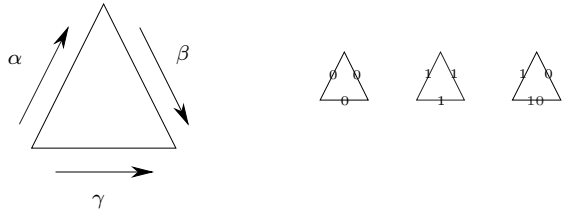
\begin{figure}[htbp]
\begin{center}
\setlength{\unitlength}{0.00083333in}
\begingroup\makeatletter\ifx\SetFigFont\undefined%
\gdef\SetFigFont#1#2#3#4#5{%
  \reset@font\fontsize{#1}{#2pt}%
  \fontfamily{#3}\fontseries{#4}\fontshape{#5}%
  \selectfont}%
\fi\endgroup%
{\renewcommand{\dashlinestretch}{30}
\begin{picture}(3462,1333)(0,-10)
\put(3234,674){\makebox(0,0)[lb]{{\SetFigFont{5}{6.0}{\rmdefault}{\mddefault}{\updefault}$10$}}}
\path(2550,706)(2700,1006)(2850,706)(2550,706)
\path(75,556)(375,1156)
\blacken\path(348.167,1035.252)(375.000,1156.000)(294.502,1062.085)(337.434,1080.868)(348.167,1035.252)
\path(825,1156)(1125,556)
\blacken\path(1044.502,649.915)(1125.000,556.000)(1098.167,676.748)(1087.434,631.132)(1044.502,649.915)
\path(300,256)(900,256)
\blacken\path(780.000,226.000)(900.000,256.000)(780.000,286.000)(816.000,256.000)(780.000,226.000)
\path(3150,706)(3300,1006)(3450,706)(3150,706)
\path(1950,706)(2100,1006)(2250,706)(1950,706)
\put(1050,931){\makebox(0,0)[lb]{{\SetFigFont{8}{9.6}{\rmdefault}{\mddefault}{\updefault}$\beta$}}}
\put(0,931){\makebox(0,0)[lb]{{\SetFigFont{8}{9.6}{\rmdefault}{\mddefault}{\updefault}$\alpha$}}}
\put(525,31){\makebox(0,0)[lb]{{\SetFigFont{8}{9.6}{\rmdefault}{\mddefault}{\updefault}$\gamma$}}}
\put(2075,672){\makebox(0,0)[lb]{{\SetFigFont{5}{6.0}{\rmdefault}{\mddefault}{\updefault}$0$}}}
\put(2756,841){\makebox(0,0)[lb]{{\SetFigFont{5}{6.0}{\rmdefault}{\mddefault}{\updefault}$1$}}}
\put(2595,841){\makebox(0,0)[lb]{{\SetFigFont{5}{6.0}{\rmdefault}{\mddefault}{\updefault}$1$}}}
\put(2001,836){\makebox(0,0)[lb]{{\SetFigFont{5}{6.0}{\rmdefault}{\mddefault}{\updefault}$0$}}}
\put(3186,846){\makebox(0,0)[lb]{{\SetFigFont{5}{6.0}{\rmdefault}{\mddefault}{\updefault}$1$}}}
\put(3344,842){\makebox(0,0)[lb]{{\SetFigFont{5}{6.0}{\rmdefault}{\mddefault}{\updefault}$0$}}}
\put(2675,677){\makebox(0,0)[lb]{{\SetFigFont{5}{6.0}{\rmdefault}{\mddefault}{\updefault}$1$}}}
\put(2150,835){\makebox(0,0)[lb]{{\SetFigFont{5}{6.0}{\rmdefault}{\mddefault}{\updefault}$0$}}}
\path(150,406)(600,1306)(1050,406)(150,406)
\end{picture}
}

\end{center}
\caption{The puzzles of Knutson and Tao (the pieces may be
  rotated but not reflected).\lremind{puzzleD}
}
\label{puzzleD}
\end{figure}

{\bf Littlewood-Richardson rule (puzzle version).}  The
Littlewood-Richardson coefficient $c^{\gamma}_{\al, \be}$ may be
computed as follows.  We write $\al$, $\be$, and $\ga$ in terms of
strings of $n$ $0$'s and $1$'s, of which $k$ are $1$'s.  We write
the digits of $\al$, $\be$, and $\ga$ along the sides of a lattice
equilateral triangle of side length $n$, as shown in
Figure~\ref{puzzleD}. Then $c^{\gamma}_{\al, \be}$ is the number of ways
of filling in this puzzle with unit triangular pieces of the sort
shown in Figure~\ref{puzzleD}, where two pieces can share an edge only
if the edges have the same label.  

\begin{figure}[htbp]
\begin{center}
\setlength{\unitlength}{0.00083333in}
\begingroup\makeatletter\ifx\SetFigFont\undefined%
\gdef\SetFigFont#1#2#3#4#5{%
  \reset@font\fontsize{#1}{#2pt}%
  \fontfamily{#3}\fontseries{#4}\fontshape{#5}%
  \selectfont}%
\fi\endgroup%
{\renewcommand{\dashlinestretch}{30}
\begin{picture}(535,436)(0,-10)
\path(162,380)(462,380)(312,80)
\path(162,380)(12,80)(312,80)
\put(52,189){\makebox(0,0)[lb]{{\SetFigFont{8}{9.6}{\rmdefault}{\mddefault}{\updefault}$0$}}}
\put(351,187){\makebox(0,0)[lb]{{\SetFigFont{8}{9.6}{\rmdefault}{\mddefault}{\updefault}$0$}}}
\put(126,31){\makebox(0,0)[lb]{{\SetFigFont{8}{9.6}{\rmdefault}{\mddefault}{\updefault}$1$}}}
\put(255,337){\makebox(0,0)[lb]{{\SetFigFont{8}{9.6}{\rmdefault}{\mddefault}{\updefault}$1$}}}
\end{picture}
}

\end{center}
\caption{The original rhombus piece of Knutson and Tao (obtained by gluing
two of the third piece of Fig.~\ref{puzzleD} together).
It may be rotated but not reflected.\lremind{rhombus}
}
\label{rhombus}
\end{figure}
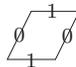

Examples are shown in Figure~\ref{keyexample}.  As an enlightening
exercise, the reader may enjoy showing that
$c^{101010}_{010101,010101}=2$.  (In \cite{knutson:puzzle}, the third
piece is different: as such a piece can only be glued along the $10$
edge along another such piece, one may as well consider ``rhombi''
consisting of two such pieces glued along $10$-edges, see
Fig.~\ref{rhombus}.  However, we shall see that these triangular
pieces are more suitable for generalization.  But the rhombi of
Knutson and Tao lead to less cluttered diagrams, so we occasionally
ignore the $10$ edges in figures.)

\subsection{The geometric Littlewood-Richardson rule in the guise of
partially completed puzzles}
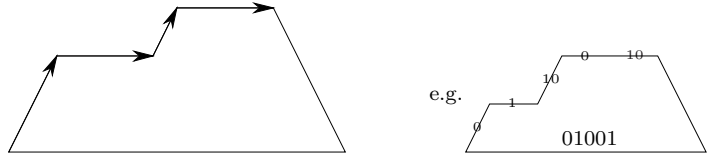
\begin{figure}[htbp]
\begin{center}
\setlength{\unitlength}{0.00083333in}
\begingroup\makeatletter\ifx\SetFigFont\undefined%
\gdef\SetFigFont#1#2#3#4#5{%
  \reset@font\fontsize{#1}{#2pt}%
  \fontfamily{#3}\fontseries{#4}\fontshape{#5}%
  \selectfont}%
\fi\endgroup%
{\renewcommand{\dashlinestretch}{30}
\begin{picture}(4374,969)(0,-10)
\put(3862,593){\makebox(0,0)[lb]{{\SetFigFont{5}{6.0}{\rmdefault}{\mddefault}{\updefault}$10$}}}
\path(312,612)(912,612)
\blacken\path(792.000,582.000)(912.000,612.000)(792.000,642.000)(828.000,612.000)(792.000,582.000)
\path(1062,912)(1662,912)
\blacken\path(1542.000,882.000)(1662.000,912.000)(1542.000,942.000)(1578.000,912.000)(1542.000,882.000)
\path(12,12)(312,612)
\blacken\path(285.167,491.252)(312.000,612.000)(231.502,518.085)(274.434,536.868)(285.167,491.252)
\path(912,612)(1062,912)
\blacken\path(1035.167,791.252)(1062.000,912.000)(981.502,818.085)(1024.434,836.868)(1035.167,791.252)
\path(2862,12)(3012,312)(3312,312)
	(3462,612)(4062,612)(4362,12)(2862,12)
\put(3462,57){\makebox(0,0)[lb]{{\SetFigFont{8}{9.6}{\rmdefault}{\mddefault}{\updefault}$01001$}}}
\put(2637,312){\makebox(0,0)[lb]{{\SetFigFont{8}{9.6}{\rmdefault}{\mddefault}{\updefault}e.g.}}}
\put(2906,143){\makebox(0,0)[lb]{{\SetFigFont{5}{6.0}{\rmdefault}{\mddefault}{\updefault}$0$}}}
\put(3124,293){\makebox(0,0)[lb]{{\SetFigFont{5}{6.0}{\rmdefault}{\mddefault}{\updefault}$1$}}}
\put(3338,443){\makebox(0,0)[lb]{{\SetFigFont{5}{6.0}{\rmdefault}{\mddefault}{\updefault}$10$}}}
\put(3580,585){\makebox(0,0)[lb]{{\SetFigFont{5}{6.0}{\rmdefault}{\mddefault}{\updefault}$0$}}}
\path(12,12)(312,612)(912,612)
	(1062,912)(1662,912)(2112,12)(12,12)
\end{picture}
}

\end{center}
\caption{Reading $\al$ from a partially filled puzzle
(horizontal $10$'s are read as $1$, leading edge $10$ is read as $0$).
The string  $\al$ describes the $k$-plane's intersection with the moving flag $M_{\bullet}$.
\lremind{readalpha}
}
\label{readalpha}
\end{figure}

\begin{figure}[htbp]
\begin{center}
\setlength{\unitlength}{0.00083333in}
\begingroup\makeatletter\ifx\SetFigFont\undefined%
\gdef\SetFigFont#1#2#3#4#5{%
  \reset@font\fontsize{#1}{#2pt}%
  \fontfamily{#3}\fontseries{#4}\fontshape{#5}%
  \selectfont}%
\fi\endgroup%
{\renewcommand{\dashlinestretch}{30}
\begin{picture}(4374,969)(0,-10)
\put(4262,149){\makebox(0,0)[lb]{{\SetFigFont{5}{6.0}{\rmdefault}{\mddefault}{\updefault}$1$}}}
\path(312,612)(912,612)
\blacken\path(792.000,582.000)(912.000,612.000)(792.000,642.000)(828.000,612.000)(792.000,582.000)
\path(1062,912)(1662,912)
\blacken\path(1542.000,882.000)(1662.000,912.000)(1542.000,942.000)(1578.000,912.000)(1542.000,882.000)
\path(1662,912)(2112,12)
\blacken\path(2031.502,105.915)(2112.000,12.000)(2085.167,132.748)(2074.434,87.132)(2031.502,105.915)
\path(2862,12)(3012,312)(3312,312)
	(3462,612)(4062,612)(4362,12)(2862,12)
\put(3462,57){\makebox(0,0)[lb]{{\SetFigFont{8}{9.6}{\rmdefault}{\mddefault}{\updefault}$00011$}}}
\put(2637,312){\makebox(0,0)[lb]{{\SetFigFont{8}{9.6}{\rmdefault}{\mddefault}{\updefault}e.g.}}}
\put(3113,288){\makebox(0,0)[lb]{{\SetFigFont{5}{6.0}{\rmdefault}{\mddefault}{\updefault}$10$}}}
\put(3344,443){\makebox(0,0)[lb]{{\SetFigFont{5}{6.0}{\rmdefault}{\mddefault}{\updefault}$10$}}}
\put(3588,586){\makebox(0,0)[lb]{{\SetFigFont{5}{6.0}{\rmdefault}{\mddefault}{\updefault}$0$}}}
\put(3861,587){\makebox(0,0)[lb]{{\SetFigFont{5}{6.0}{\rmdefault}{\mddefault}{\updefault}$1$}}}
\put(4118,430){\makebox(0,0)[lb]{{\SetFigFont{5}{6.0}{\rmdefault}{\mddefault}{\updefault}$1$}}}
\path(12,12)(312,612)(912,612)
	(1062,912)(1662,912)(2112,12)(12,12)
\end{picture}
}

\end{center}
\caption{Reading $\be$ from a partially filled puzzle
($10$'s are read as $0$, and leading edge $10$ turns the first $1$ afterward to $0$).  The string $\be$ describes the $k$-plane's intersection with the fixed flag $F_{\bullet}$.\lremind{readbeta}}
\label{readbeta}
\end{figure}
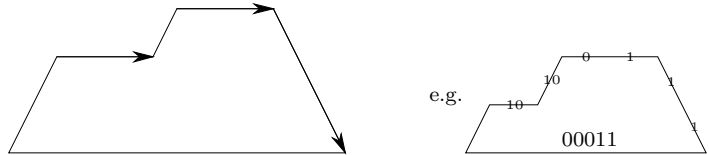

\lremind{pp}\label{pp}The geometric Littlewood-Richardson rule may be
understood in terms of partially completed puzzles.  The description
here is slightly cleaner than, although equivalent to, that of
\cite[\S A.2]{vakil:checkers}.  To each term in the degeneration
order, we associate the shape of a partially filled puzzle, as shown
in Figure~\ref{degenorder}.  Notice that there is at most one
southwest-to-northeast edge in the middle of the puzzle; call this the
{\em leading edge}.  The cycles we will consider at step $d_i$ are as
follows: We label the northwestern, northern, and northeastern edges
of the partially filled puzzle, where the left and right edges of the
large triangle have labels $0$ or $1$, and the remaining northern
border edges are labeled $0$, $1$, or $10$.  To such a labeled
partially filled puzzle, we define two strings of $0$'s and $1$'s (of
which $k$ are $1$'s), denoted $\al$ and $\be$ as follows.  For $\al$,
we read off the list of $0$'s and $1$'s on the edges in the order
shown in Figure~\ref{readalpha}, where horizontal $10$ edges are
interpreted as $1$, and a $10$ leading edge is interpreted as a $0$.
For $\be$, we read off the list of $0$'s and $1$'s as shown in
Figure~\ref{readbeta}, where we first interpret the horizontal $10$
edges as $0$'s, and then if the leading edge is $10$, we turn the
first $1$ in $\be$ after the leading edge to $0$.

Note  that $\al$ and $\be$ are not independent.
This is clearest after the last step, where the above description 
ensures that $\al=\be$.

\begin{figure}[htbp]
\begin{center}
\setlength{\unitlength}{0.00083333in}
\begingroup\makeatletter\ifx\SetFigFont\undefined%
\gdef\SetFigFont#1#2#3#4#5{%
  \reset@font\fontsize{#1}{#2pt}%
  \fontfamily{#3}\fontseries{#4}\fontshape{#5}%
  \selectfont}%
\fi\endgroup%
{\renewcommand{\dashlinestretch}{30}
\begin{picture}(3924,2008)(0,-10)
\put(3807,609){\makebox(0,0)[lb]{{\SetFigFont{5}{6.0}{\rmdefault}{\mddefault}{\updefault}$1$}}}
\path(912,1081)(1062,781)(1212,1081)
\blacken\path(231.502,274.915)(312.000,181.000)(285.167,301.748)(274.434,256.132)(231.502,274.915)
\dashline{60.000}(312,181)(87,631)
\blacken\path(531.502,274.915)(612.000,181.000)(585.167,301.748)(574.434,256.132)(531.502,274.915)
\dashline{60.000}(612,181)(312,781)
\blacken\path(831.502,274.915)(912.000,181.000)(885.167,301.748)(874.434,256.132)(831.502,274.915)
\dashline{60.000}(912,181)(537,931)
\blacken\path(1131.502,274.915)(1212.000,181.000)(1185.167,301.748)(1174.434,256.132)(1131.502,274.915)
\dashline{60.000}(1212,181)(762,1081)
\blacken\path(1431.502,274.915)(1512.000,181.000)(1485.167,301.748)(1474.434,256.132)(1431.502,274.915)
\dashline{60.000}(1512,181)(1137,931)
\blacken\path(3638.833,301.748)(3612.000,181.000)(3692.498,274.915)(3649.566,256.132)(3638.833,301.748)
\dashline{60.000}(3612,181)(3837,631)
\blacken\path(3338.833,301.748)(3312.000,181.000)(3392.498,274.915)(3349.566,256.132)(3338.833,301.748)
\dashline{60.000}(3312,181)(3687,931)
\blacken\path(3038.833,301.748)(3012.000,181.000)(3092.498,274.915)(3049.566,256.132)(3038.833,301.748)
\dashline{60.000}(3012,181)(3462,1081)
\blacken\path(2738.833,301.748)(2712.000,181.000)(2792.498,274.915)(2749.566,256.132)(2738.833,301.748)
\dashline{60.000}(2712,181)(3162,1081)
\blacken\path(2438.833,301.748)(2412.000,181.000)(2492.498,274.915)(2449.566,256.132)(2438.833,301.748)
\dashline{60.000}(2412,181)(2637,631)
\path(2562,781)(2712,481)(2862,781)
\dashline{60.000}(2937,931)(3012,181)
\blacken\path(2970.208,297.419)(3012.000,181.000)(3029.911,303.390)(3003.642,264.583)(2970.208,297.419)
\path(162,781)(762,1981)(1212,1081)
\path(2562,781)(3162,1981)(3612,1081)
\path(2412,481)(2562,781)(2862,781)
	(3012,1081)(3612,1081)(3912,481)(2412,481)
\put(312,31){\makebox(0,0)[lb]{{\SetFigFont{8}{9.6}{\rmdefault}{\mddefault}{\updefault}$0$}}}
\put(612,31){\makebox(0,0)[lb]{{\SetFigFont{8}{9.6}{\rmdefault}{\mddefault}{\updefault}$1$}}}
\put(912,31){\makebox(0,0)[lb]{{\SetFigFont{8}{9.6}{\rmdefault}{\mddefault}{\updefault}$0$}}}
\put(1212,31){\makebox(0,0)[lb]{{\SetFigFont{8}{9.6}{\rmdefault}{\mddefault}{\updefault}$0$}}}
\put(1512,31){\makebox(0,0)[lb]{{\SetFigFont{8}{9.6}{\rmdefault}{\mddefault}{\updefault}$1$}}}
\put(2337,31){\makebox(0,0)[lb]{{\SetFigFont{8}{9.6}{\rmdefault}{\mddefault}{\updefault}$0$}}}
\put(2637,31){\makebox(0,0)[lb]{{\SetFigFont{8}{9.6}{\rmdefault}{\mddefault}{\updefault}$0$}}}
\put(2937,31){\makebox(0,0)[lb]{{\SetFigFont{8}{9.6}{\rmdefault}{\mddefault}{\updefault}$0$}}}
\put(3237,31){\makebox(0,0)[lb]{{\SetFigFont{8}{9.6}{\rmdefault}{\mddefault}{\updefault}$1$}}}
\put(3537,31){\makebox(0,0)[lb]{{\SetFigFont{8}{9.6}{\rmdefault}{\mddefault}{\updefault}$1$}}}
\put(686,1381){\makebox(0,0)[lb]{{\SetFigFont{8}{9.6}{\rmdefault}{\mddefault}{\updefault}$\al$}}}
\put(3085,1381){\makebox(0,0)[lb]{{\SetFigFont{8}{9.6}{\rmdefault}{\mddefault}{\updefault}$\be$}}}
\put(63,620){\makebox(0,0)[lb]{{\SetFigFont{5}{6.0}{\rmdefault}{\mddefault}{\updefault}$0$}}}
\put(289,756){\makebox(0,0)[lb]{{\SetFigFont{5}{6.0}{\rmdefault}{\mddefault}{\updefault}$1$}}}
\put(497,908){\makebox(0,0)[lb]{{\SetFigFont{5}{6.0}{\rmdefault}{\mddefault}{\updefault}$10$}}}
\put(734,1056){\makebox(0,0)[lb]{{\SetFigFont{5}{6.0}{\rmdefault}{\mddefault}{\updefault}$0$}}}
\put(954,915){\makebox(0,0)[lb]{{\SetFigFont{5}{6.0}{\rmdefault}{\mddefault}{\updefault}$0$}}}
\put(1109,908){\makebox(0,0)[lb]{{\SetFigFont{5}{6.0}{\rmdefault}{\mddefault}{\updefault}$1$}}}
\put(1016,1058){\makebox(0,0)[lb]{{\SetFigFont{5}{6.0}{\rmdefault}{\mddefault}{\updefault}$10$}}}
\put(2620,613){\makebox(0,0)[lb]{{\SetFigFont{5}{6.0}{\rmdefault}{\mddefault}{\updefault}$0$}}}
\put(2753,620){\makebox(0,0)[lb]{{\SetFigFont{5}{6.0}{\rmdefault}{\mddefault}{\updefault}$1$}}}
\put(2671,753){\makebox(0,0)[lb]{{\SetFigFont{5}{6.0}{\rmdefault}{\mddefault}{\updefault}$10$}}}
\put(2884,897){\makebox(0,0)[lb]{{\SetFigFont{5}{6.0}{\rmdefault}{\mddefault}{\updefault}$10$}}}
\put(3132,1052){\makebox(0,0)[lb]{{\SetFigFont{5}{6.0}{\rmdefault}{\mddefault}{\updefault}$0$}}}
\put(3427,1053){\makebox(0,0)[lb]{{\SetFigFont{5}{6.0}{\rmdefault}{\mddefault}{\updefault}$1$}}}
\put(3664,908){\makebox(0,0)[lb]{{\SetFigFont{5}{6.0}{\rmdefault}{\mddefault}{\updefault}$1$}}}
\path(12,481)(162,781)(462,781)
	(612,1081)(1212,1081)(1512,481)(12,481)
\end{picture}
}

\end{center}
\caption{Alternate (equivalent) description of how to read $\al$ and $\be$ (cf.\ Figures~\ref{readalpha} and~\ref{readbeta}).\lremind{newab}}
\label{newab}
\end{figure}
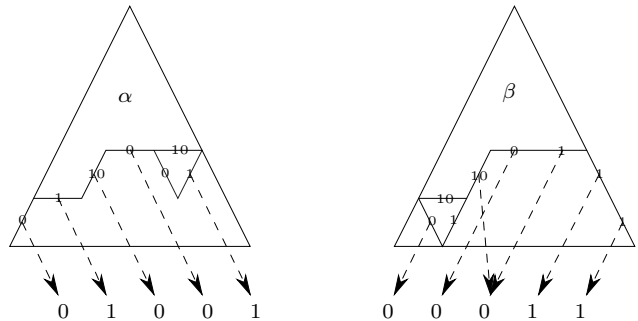

The reader may prefer the following equivalent description of $\al$
and $\be$, depicted in Figure~\ref{newab}.  To each horizontal edge
$10$, temporarily glue on a copy of the last triangle of
Figure~\ref{puzzleD} (with sides $0/1/10$ read clockwise), so now all
exposed edges are $0$'s and $1$'s, except possibly the leading edge.
Then $\al$ is the string visible looking northwest form the bottom
edge of the puzzle, except a leading edge $01$ is read as $0$; and
$\be$ is the string visible looking northeast from the bottom edge of
the puzzle, except that a leading edge $01$ turns the next $1$ into
$0$.

To such a partially filled puzzle, we associate a subvariety
of the Grassmannian as follows.
(This construction will be central to the
generalizations in every section of Part~\ref{partone}!)
  We fix two flags $F_{\bullet}$
and $M_{\bullet}$ in the relative position determined by the partially
completed puzzle, and we consider \lremind{puzzlevariety}
\begin{equation}\label{puzzlevariety}
\overline{ \Omega_{\al}(M_{\bullet}) \cap \Omega_{\be}(F_{\bullet})}.
\end{equation}
Call such a subvariety of the Grassmannian a {\em puzzle variety}.
The reader is
strongly encouraged to see this definition in practice by examining a
couple of panels of Figure~\ref{keyexample}.

Then the geometric Littlewood-Richardson rule may now be described
quickly as follows.  Suppose we have a partially filled puzzle,
corresponding to two flags $F_{\bullet}$ and $M_{\bullet}$ in given relative
position, and a puzzle variety (those $k$-planes meeting the two
reference flags as described in \eqref{puzzlevariety}).  Then if the
moving flag $M_{\bullet}$ is degenerated (to the next step in the
degeneration order, Figure~\ref{degenorder}), the cycle in the
Grassmannian degenerates to all possible ways of adding more
triangular pieces to the puzzle to create a ``next-larger'' partially
completed puzzle.  (Note that the edges on the bottom of the puzzle
must be labeled $0$ or $1$, and may not be labeled $10$.)

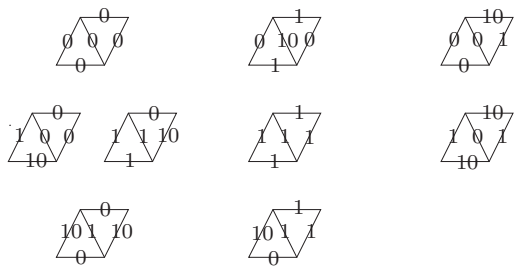
\begin{figure}[htbp]
\begin{center}
\setlength{\unitlength}{0.00083333in}
\begingroup\makeatletter\ifx\SetFigFont\undefined%
\gdef\SetFigFont#1#2#3#4#5{%
  \reset@font\fontsize{#1}{#2pt}%
  \fontfamily{#3}\fontseries{#4}\fontshape{#5}%
  \selectfont}%
\fi\endgroup%
{\renewcommand{\dashlinestretch}{30}
\begin{picture}(3250,1647)(0,-10)
\put(1634,31){\makebox(0,0)[lb]{{\SetFigFont{8}{9.6}{\rmdefault}{\mddefault}{\updefault}$0$}}}
\path(312,77)(462,377)(762,377)
	(612,77)(312,77)
\path(1512,1277)(1662,1577)(1962,1577)
	(1812,1277)(1512,1277)
\path(1512,677)(1662,977)(1962,977)
	(1812,677)(1512,677)
\path(1512,77)(1662,377)(1962,377)
	(1812,77)(1512,77)
\path(2712,677)(2862,977)(3162,977)
	(3012,677)(2712,677)
\path(2712,1277)(2862,1577)(3162,1577)
	(3012,1277)(2712,1277)
\path(462,1577)(612,1277)
\drawline(20,902)(20,902)
\path(312,1277)(462,1577)(762,1577)
	(612,1277)(312,1277)
\path(1662,1577)(1812,1277)
\path(2862,977)(3012,677)
\path(1812,677)(1662,977)
\path(1662,377)(1812,77)
\path(612,77)(462,377)
\path(312,677)(162,977)
\path(3012,1277)(2862,1577)
\path(612,677)(762,977)(1062,977)
	(912,677)(612,677)
\path(912,677)(762,977)
\put(430,1233){\makebox(0,0)[lb]{{\SetFigFont{8}{9.6}{\rmdefault}{\mddefault}{\updefault}$0$}}}
\put(506,1390){\makebox(0,0)[lb]{{\SetFigFont{8}{9.6}{\rmdefault}{\mddefault}{\updefault}$0$}}}
\put(662,1390){\makebox(0,0)[lb]{{\SetFigFont{8}{9.6}{\rmdefault}{\mddefault}{\updefault}$0$}}}
\put(348,1383){\makebox(0,0)[lb]{{\SetFigFont{8}{9.6}{\rmdefault}{\mddefault}{\updefault}$0$}}}
\put(576,1548){\makebox(0,0)[lb]{{\SetFigFont{8}{9.6}{\rmdefault}{\mddefault}{\updefault}$0$}}}
\put(1793,1537){\makebox(0,0)[lb]{{\SetFigFont{8}{9.6}{\rmdefault}{\mddefault}{\updefault}$1$}}}
\put(1546,1383){\makebox(0,0)[lb]{{\SetFigFont{8}{9.6}{\rmdefault}{\mddefault}{\updefault}$0$}}}
\put(1686,1391){\makebox(0,0)[lb]{{\SetFigFont{8}{9.6}{\rmdefault}{\mddefault}{\updefault}$10$}}}
\put(1858,1393){\makebox(0,0)[lb]{{\SetFigFont{8}{9.6}{\rmdefault}{\mddefault}{\updefault}$0$}}}
\put(1644,1234){\makebox(0,0)[lb]{{\SetFigFont{8}{9.6}{\rmdefault}{\mddefault}{\updefault}$1$}}}
\put(2759,1396){\makebox(0,0)[lb]{{\SetFigFont{8}{9.6}{\rmdefault}{\mddefault}{\updefault}$0$}}}
\put(2907,1397){\makebox(0,0)[lb]{{\SetFigFont{8}{9.6}{\rmdefault}{\mddefault}{\updefault}$0$}}}
\put(3066,1393){\makebox(0,0)[lb]{{\SetFigFont{8}{9.6}{\rmdefault}{\mddefault}{\updefault}$1$}}}
\put(2964,1537){\makebox(0,0)[lb]{{\SetFigFont{8}{9.6}{\rmdefault}{\mddefault}{\updefault}$10$}}}
\put(2821,1233){\makebox(0,0)[lb]{{\SetFigFont{8}{9.6}{\rmdefault}{\mddefault}{\updefault}$0$}}}
\put(287,939){\makebox(0,0)[lb]{{\SetFigFont{8}{9.6}{\rmdefault}{\mddefault}{\updefault}$0$}}}
\put(49,793){\makebox(0,0)[lb]{{\SetFigFont{8}{9.6}{\rmdefault}{\mddefault}{\updefault}$1$}}}
\put(207,790){\makebox(0,0)[lb]{{\SetFigFont{8}{9.6}{\rmdefault}{\mddefault}{\updefault}$0$}}}
\put(357,790){\makebox(0,0)[lb]{{\SetFigFont{8}{9.6}{\rmdefault}{\mddefault}{\updefault}$0$}}}
\put(115,639){\makebox(0,0)[lb]{{\SetFigFont{8}{9.6}{\rmdefault}{\mddefault}{\updefault}$10$}}}
\put(649,788){\makebox(0,0)[lb]{{\SetFigFont{8}{9.6}{\rmdefault}{\mddefault}{\updefault}$1$}}}
\put(737,639){\makebox(0,0)[lb]{{\SetFigFont{8}{9.6}{\rmdefault}{\mddefault}{\updefault}$1$}}}
\put(823,788){\makebox(0,0)[lb]{{\SetFigFont{8}{9.6}{\rmdefault}{\mddefault}{\updefault}$1$}}}
\put(941,798){\makebox(0,0)[lb]{{\SetFigFont{8}{9.6}{\rmdefault}{\mddefault}{\updefault}$10$}}}
\put(887,935){\makebox(0,0)[lb]{{\SetFigFont{8}{9.6}{\rmdefault}{\mddefault}{\updefault}$0$}}}
\put(1795,938){\makebox(0,0)[lb]{{\SetFigFont{8}{9.6}{\rmdefault}{\mddefault}{\updefault}$1$}}}
\put(1559,788){\makebox(0,0)[lb]{{\SetFigFont{8}{9.6}{\rmdefault}{\mddefault}{\updefault}$1$}}}
\put(1704,789){\makebox(0,0)[lb]{{\SetFigFont{8}{9.6}{\rmdefault}{\mddefault}{\updefault}$1$}}}
\put(1857,785){\makebox(0,0)[lb]{{\SetFigFont{8}{9.6}{\rmdefault}{\mddefault}{\updefault}$1$}}}
\put(1639,635){\makebox(0,0)[lb]{{\SetFigFont{8}{9.6}{\rmdefault}{\mddefault}{\updefault}$1$}}}
\put(2964,939){\makebox(0,0)[lb]{{\SetFigFont{8}{9.6}{\rmdefault}{\mddefault}{\updefault}$10$}}}
\put(2753,788){\makebox(0,0)[lb]{{\SetFigFont{8}{9.6}{\rmdefault}{\mddefault}{\updefault}$1$}}}
\put(2907,788){\makebox(0,0)[lb]{{\SetFigFont{8}{9.6}{\rmdefault}{\mddefault}{\updefault}$0$}}}
\put(3062,789){\makebox(0,0)[lb]{{\SetFigFont{8}{9.6}{\rmdefault}{\mddefault}{\updefault}$1$}}}
\put(2805,631){\makebox(0,0)[lb]{{\SetFigFont{8}{9.6}{\rmdefault}{\mddefault}{\updefault}$10$}}}
\put(335,195){\makebox(0,0)[lb]{{\SetFigFont{8}{9.6}{\rmdefault}{\mddefault}{\updefault}$10$}}}
\put(501,195){\makebox(0,0)[lb]{{\SetFigFont{8}{9.6}{\rmdefault}{\mddefault}{\updefault}$1$}}}
\put(651,200){\makebox(0,0)[lb]{{\SetFigFont{8}{9.6}{\rmdefault}{\mddefault}{\updefault}$10$}}}
\put(583,333){\makebox(0,0)[lb]{{\SetFigFont{8}{9.6}{\rmdefault}{\mddefault}{\updefault}$0$}}}
\put(431,35){\makebox(0,0)[lb]{{\SetFigFont{8}{9.6}{\rmdefault}{\mddefault}{\updefault}$0$}}}
\put(1790,348){\makebox(0,0)[lb]{{\SetFigFont{8}{9.6}{\rmdefault}{\mddefault}{\updefault}$1$}}}
\put(1867,195){\makebox(0,0)[lb]{{\SetFigFont{8}{9.6}{\rmdefault}{\mddefault}{\updefault}$1$}}}
\put(1704,195){\makebox(0,0)[lb]{{\SetFigFont{8}{9.6}{\rmdefault}{\mddefault}{\updefault}$1$}}}
\put(1529,189){\makebox(0,0)[lb]{{\SetFigFont{8}{9.6}{\rmdefault}{\mddefault}{\updefault}$10$}}}
\path(12,677)(162,977)(462,977)
	(312,677)(12,677)
\end{picture}
}
\end{center}
\caption{The puzzle interpretation of the geometric Littlewood-Richardson rule.\lremind{pairs}}
\label{pairs}
\end{figure}

This involves adding two adjacent triangles (forming a rhombus), and
possibly a third triangle to ``complete a row''.  There is no choice
for the third triangle: if we know two sides of a puzzle piece, then
we know the third.  The reader may readily verify that there are nine
possible choices for the pair of pieces, shown suggestively in
Figure~\ref{pairs}.  Once the northern edge and the northwestern edge
are specified, there is usually one choice for the rhombus, and there
may be two.  For those familiar with the checkerboard description of
this rule, they correspond to the ``stay'' and ``swap'' options
respectively.

\subsection{Connection to tableaux}
\begin{figure}[htbp]
\begin{center}
\setlength{\unitlength}{0.00083333in}
\begingroup\makeatletter\ifx\SetFigFont\undefined%
\gdef\SetFigFont#1#2#3#4#5{%
  \reset@font\fontsize{#1}{#2pt}%
  \fontfamily{#3}\fontseries{#4}\fontshape{#5}%
  \selectfont}%
\fi\endgroup%
{\renewcommand{\dashlinestretch}{30}
\begin{picture}(546,434)(0,-10)
\put(284,335){\makebox(0,0)[lb]{{\SetFigFont{8}{9.6}{\rmdefault}{\mddefault}{\updefault}$1$}}}
\path(162,371)(312,71)
\put(147,186){\makebox(0,0)[lb]{{\SetFigFont{8}{9.6}{\rmdefault}{\mddefault}{\updefault}$10$}}}
\path(12,71)(162,371)(462,371)
	(312,71)(12,71)
\put(362,186){\makebox(0,0)[lb]{{\SetFigFont{8}{9.6}{\rmdefault}{\mddefault}{\updefault}$0$}}}
\put(53,183){\makebox(0,0)[lb]{{\SetFigFont{8}{9.6}{\rmdefault}{\mddefault}{\updefault}$0$}}}
\put(135,31){\makebox(0,0)[lb]{{\SetFigFont{8}{9.6}{\rmdefault}{\mddefault}{\updefault}$1$}}}
\end{picture}
}

\end{center}
\caption{These puzzle pieces (in this orientation!) correspond to the entries in the tableau
Littlewood-Richardson rule.\lremind{puztab}}
\label{puztab}
\end{figure}
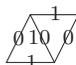

The Littlewood-Richardson rule is 
traditionally
given in terms of tableaux.
The bijection to tableaux is straightforward:  whenever the
$k$-plane changes its relationship to the moving flag $M_{\bullet}$ (e.g.\ the
moves marked $\dagger$ in Figure~\ref{keyexample}),
a number is put into the tableau.  See \cite[\S 3.1]{vakil:checkers}
for the precise statement of which number is placed, and in which row.
An elegant  bijection from puzzles to tableaux, due to Terry Tao, is
given in \cite[Fig.~11]{vakil:checkers}; whenever the puzzle piece(s)
of Figure~\ref{puztab} appear in the puzzle, a number is placed
in the tableau.

While discussing tableaux, we should mention K.~Purbhoo's beautiful
``mosaics'', with which he not only bijects puzzles and tableaux, but
proves the commutativity and associativity of each \cite{purbhoo}.

\section{The $K$-theory (or Grothendieck ring) of the Grassmannian}
\lremind{K}\label{K}The Grothendieck groups or $K$-theory of the
Grassmannian is generated by the classes of the structure sheaves of
the Schubert cells.  Buch gave a Littlewood-Richardson rule in
$K$-theory in \cite{buch:K}.
Buch's rule states that $K$-theory Littlewood-Richardson coefficients
enumerate certain ``set-valued tableaux''.  There is also a sign: the
sign of $c^{\ga}_{\al, \be}$ is $(-1)^{\codim \Omega_{\ga} - \codim
  \Omega_{\al} - \codim \Omega_{\be} }$, depending on the difference of the
codimension of $\ga$ from its ``expected codimension''.

As an aside, we note that Buch conjectured positivity --- with this
sign convention --- in the Grothendieck ring of a flag variety in
general.  A beautiful proof was given by M.~Brion in
\cite{brion:positivity}.

\begin{figure}[htbp]
\begin{center}
\setlength{\unitlength}{0.00083333in}
\begingroup\makeatletter\ifx\SetFigFont\undefined%
\gdef\SetFigFont#1#2#3#4#5{%
  \reset@font\fontsize{#1}{#2pt}%
  \fontfamily{#3}\fontseries{#4}\fontshape{#5}%
  \selectfont}%
\fi\endgroup%
{\renewcommand{\dashlinestretch}{30}
\begin{picture}(683,671)(0,-10)
\put(360,121){\makebox(0,0)[lb]{{\SetFigFont{8}{9.6}{\rmdefault}{\mddefault}{\updefault}$1$}}}
\path(162,312)(312,612)(462,312)(162,312)
\put(133,572){\makebox(0,0)[lb]{{\SetFigFont{8}{9.6}{\rmdefault}{\mddefault}{\updefault}$0$}}}
\put(437,569){\makebox(0,0)[lb]{{\SetFigFont{8}{9.6}{\rmdefault}{\mddefault}{\updefault}$1$}}}
\put(499,412){\makebox(0,0)[lb]{{\SetFigFont{8}{9.6}{\rmdefault}{\mddefault}{\updefault}$0$}}}
\put(63,417){\makebox(0,0)[lb]{{\SetFigFont{8}{9.6}{\rmdefault}{\mddefault}{\updefault}$1$}}}
\put(203,119){\makebox(0,0)[lb]{{\SetFigFont{8}{9.6}{\rmdefault}{\mddefault}{\updefault}$0$}}}
\path(12,612)(612,612)(312,12)(12,612)
\end{picture}
}

\end{center}
\caption{The new $K$-theory puzzle piece.  It may {\em not} be rotated.
(The central edges might be plausibly labeled $10$, although this
has no mathematical effect.)\lremind{puzzleK}}
\label{puzzleK}
\end{figure}
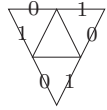

Buch's generalization of the tableau rule to
$K$-theory yields a generalization to $K$-theory of the geometric
checker rule \cite[Thm.~3.4]{vakil:checkers} and the puzzle rule
\cite[Thm.~3.6]{vakil:checkers}.  The puzzle rule is as follows: there
is a new puzzle piece, shown in Figure~\ref{puzzleK}, and the
$K$-theory Littlewood-Richardson coefficient is the number of puzzles
completed with the usual pieces and the $K$-theory piece, where each
puzzle is counted with sign corresponding to the number of $K$-theory
pieces.  As an example, the reader can verify that in $K$-theory,
\lremind{basicK}
\begin{equation}
\label{basicK}
\si_{0101}^2 = \si_{0110} + \si_{1001} - \si_{1010}.
\end{equation}
in string notation, extending \eqref{basic}.

The puzzle statement is quite clean.  The checker statement (not given 
here) suggests a geometric conjecture \cite[Conj.~3.5]{vakil:checkers}
(by the second author with Buch, with additional comments by
Knutson), which is also quite clean.  We will now do something
seemingly perverse: give a {\em less clean} puzzle statement.  The advantage
is that the puzzle rule will have an interpretation as a {\em
  conjectural $K$-theoretic geometric Littlewood-Richardson rule.}
(When we extend our discussion to equivariant cohomology, \S \ref{HT},
and equivariant $K$-theory, \S \ref{KT}, the checker description will
break down, but the geometric puzzle description will still work.)

\begin{figure}[htbp]
\begin{center}
\setlength{\unitlength}{0.00083333in}
\begingroup\makeatletter\ifx\SetFigFont\undefined%
\gdef\SetFigFont#1#2#3#4#5{%
  \reset@font\fontsize{#1}{#2pt}%
  \fontfamily{#3}\fontseries{#4}\fontshape{#5}%
  \selectfont}%
\fi\endgroup%
{\renewcommand{\dashlinestretch}{30}
\begin{picture}(2348,432)(0,-10)
\put(1259,185){\makebox(0,0)[lb]{{\SetFigFont{8}{9.6}{\rmdefault}{\mddefault}{\updefault}$10K$}}}
\path(912,71)(1062,371)(1362,371)
	(1212,71)(912,71)
\path(1812,71)(1962,371)(2262,371)
	(2112,71)(1812,71)
\path(162,371)(312,71)
\path(1062,371)(1212,71)
\path(1962,371)(2112,71)
\put(284,333){\makebox(0,0)[lb]{{\SetFigFont{8}{9.6}{\rmdefault}{\mddefault}{\updefault}$0$}}}
\put(205,188){\makebox(0,0)[lb]{{\SetFigFont{8}{9.6}{\rmdefault}{\mddefault}{\updefault}$1$}}}
\put(150,31){\makebox(0,0)[lb]{{\SetFigFont{8}{9.6}{\rmdefault}{\mddefault}{\updefault}$1$}}}
\put(1035,31){\makebox(0,0)[lb]{{\SetFigFont{8}{9.6}{\rmdefault}{\mddefault}{\updefault}$0$}}}
\put(1900,33){\makebox(0,0)[lb]{{\SetFigFont{8}{9.6}{\rmdefault}{\mddefault}{\updefault}$10$}}}
\put(2164,181){\makebox(0,0)[lb]{{\SetFigFont{8}{9.6}{\rmdefault}{\mddefault}{\updefault}$0$}}}
\put(1986,181){\makebox(0,0)[lb]{{\SetFigFont{8}{9.6}{\rmdefault}{\mddefault}{\updefault}$10$}}}
\put(1693,181){\makebox(0,0)[lb]{{\SetFigFont{8}{9.6}{\rmdefault}{\mddefault}{\updefault}$10K$}}}
\put(1186,328){\makebox(0,0)[lb]{{\SetFigFont{8}{9.6}{\rmdefault}{\mddefault}{\updefault}$0$}}}
\put(2079,329){\makebox(0,0)[lb]{{\SetFigFont{8}{9.6}{\rmdefault}{\mddefault}{\updefault}$1$}}}
\put(57,182){\makebox(0,0)[lb]{{\SetFigFont{8}{9.6}{\rmdefault}{\mddefault}{\updefault}$1$}}}
\put(358,189){\makebox(0,0)[lb]{{\SetFigFont{8}{9.6}{\rmdefault}{\mddefault}{\updefault}$10K$}}}
\put(797,180){\makebox(0,0)[lb]{{\SetFigFont{8}{9.6}{\rmdefault}{\mddefault}{\updefault}$10K$}}}
\put(1119,180){\makebox(0,0)[lb]{{\SetFigFont{8}{9.6}{\rmdefault}{\mddefault}{\updefault}$1$}}}
\path(12,71)(162,371)(462,371)
	(312,71)(12,71)
\end{picture}
}
\end{center}
\caption{Alternate $K$-theory puzzle pieces.  They may {\em not} be rotated.
The first piece contributes a sign of $-1$.  The labels on the internal
edges don't matter.\lremind{geopuzzleK}}
\label{geopuzzleK}
\end{figure}
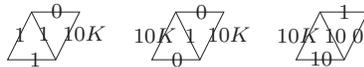

The alternate puzzle pieces are shown in Figure~\ref{geopuzzleK}.
Note that we have a new flavor of edge (labeled $10K$), which may only
be oriented southwest/northeast.  The interested reader will readily verify
that the puzzles with the alternate $K$-puzzle pieces of
Figure~\ref{geopuzzleK} are in straightforward bijection with the
puzzles using the original $K$-piece of Figure~\ref{puzzleK}.  

Corresponding to a partially-completed $K$-puzzle, there is a {\em
  $K$-theory puzzle variety}, which we describe by defining $\al$ and
$\be$ and then using the puzzle variety definition
of \eqref{puzzlevariety}.  To define $\be$, use the same recipe as before
(Fig.\ ~\ref{readbeta}), treating $10K$ as $10$.  To define $\al$
(Fig.\ ~\ref{readalpha}), treat $10K$ as $10$, except the first $1$
after that position is exchanged with the $0$ immediately to its left.

Once again, when completing a puzzle from top to bottom, there is
usually only one choice at each stage.  However, where there used to
be two choices, there are now three.  This leads to a conjectural
geometric interpretation of Buch's rule.  Recall that if a variety
degenerates into two pieces, say $X$ degenerates into $X_1 \cup X_2$,
then in homology, $[X] = [X_1] + [X_2]$, but in $K$-theory
$$[\oh_X] = [\oh_{X_1}] + [\oh_{X_2}] - [\oh_{X_1 \cap X_2}],$$
where the intersection is the scheme-theoretic intersection.
Then we have a {\em $K$-theoretic geometric Littlewood-Richardson
  rule}: when the puzzle variety breaks in the geometric
Littlewood-Richardson rule (the center left panel of
Fig.~\ref{pairs}), into $X_1$ and $X_2$, say, and $X_3$ is the puzzle
variety corresponding to the third (new, $K$-theoretic) option,
then:\lremind{Kconj}

\begin{conjecture}
  $X_3$ is a $GL(n)$-translate of the scheme-theoretic intersection $X_1 \cap
  X_2$.
\label{Kconj}
\end{conjecture}

This is straightforward to check set-theoretically.
The reader may verify that Conjecture~\ref{Kconj} is equivalent to
the $K$-theoretic conjecture \cite[Conj.~3.5]{vakil:checkers}, and hence implies Buch's combinatorial
$K$-theoretic Littlewood-Richardson rule.

The reader may prefer to instead examine the example of
Figure~\ref{keyexample}.  In the single instance when the variety
breaks into two pieces (the second panel in the top row breaks into
the union of the third panel in the top row and the first panel in the
second row), the scheme-theoretic intersection corresponds to the
locus of lines passing through the fixed point $\proj F_1$ (the
condition of the second panel in the top row), and lying in the plane
spanned by the moving line $\proj M_1$ and the fixed line $\proj F_1$
(the condition of the first panel of the second row).  This is clearly
a translate of the locus of lines passing through the fixed point
$\proj F_1$ and lying in the moving plane $\proj M_3$, which is the
puzzle variety predicted by Conjecture~\ref{Kconj}.

\section{The equivariant cohomology of the Grassmannian}
\lremind{HT}\label{HT}Suppose $T$ is the natural $n$-dimensional torus
acting on $\CC^n$.  Choose an order of the $T$-fixed basis $\vv_1$,
\dots, $\vv_n$, and let the $T$-equivariant cohomology of a point be
$\Z[y_1, \dots, y_n]$, where $y_i$ corresponds to $\vv_i$.  We next
generalize our geometric construction, at least conjecturally, to
$T$-equivariant cohomology.  The equivariant Schubert classes are
defined as equivariant cohomology classes corresponding to Schubert
varieties with respect to the fixed flag
$$
F_{\bullet} = \{ \{ 0 \} \subset \langle \vv_1 \rangle \subset \langle \vv_1, \vv_2
\rangle \subset \cdots \subset \langle \vv_1, \vv_2, \dots, \vv_n
\rangle = \C^n \} .$$ A beautiful argument showing positivity in
equivariant cohomology was given  by W.~Graham in \cite{graham:equivariant},
confirming a conjecture of D.~Peterson.

\begin{figure}[htbp]
\begin{center}
\setlength{\unitlength}{0.00083333in}
\begingroup\makeatletter\ifx\SetFigFont\undefined%
\gdef\SetFigFont#1#2#3#4#5{%
  \reset@font\fontsize{#1}{#2pt}%
  \fontfamily{#3}\fontseries{#4}\fontshape{#5}%
  \selectfont}%
\fi\endgroup%
{\renewcommand{\dashlinestretch}{30}
\begin{picture}(405,639)(0,-10)
\texture{88555555 55000000 555555 55000000 555555 55000000 555555 55000000 
	555555 55000000 555555 55000000 555555 55000000 555555 55000000 
	555555 55000000 555555 55000000 555555 55000000 555555 55000000 
	555555 55000000 555555 55000000 555555 55000000 555555 55000000 }
\shade\path(12,312)(162,612)(312,312)
	(162,12)(12,312)
\path(12,312)(162,612)(312,312)
	(162,12)(12,312)
\put(200,127){\makebox(0,0)[lb]{{\SetFigFont{8}{9.6}{\rmdefault}{\mddefault}{\updefault}$0$}}}
\put(55,125){\makebox(0,0)[lb]{{\SetFigFont{8}{9.6}{\rmdefault}{\mddefault}{\updefault}$1$}}}
\put(221,406){\makebox(0,0)[lb]{{\SetFigFont{8}{9.6}{\rmdefault}{\mddefault}{\updefault}$1$}}}
\put(59,405){\makebox(0,0)[lb]{{\SetFigFont{8}{9.6}{\rmdefault}{\mddefault}{\updefault}$0$}}}
\end{picture}
}

\end{center}
\caption{The new equivariant puzzle piece  of Knutson and Tao.
It may {\em not} be
  rotated or reflected).\lremind{eqpuzzle}}
\label{eqpuzzle}
\end{figure}
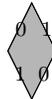

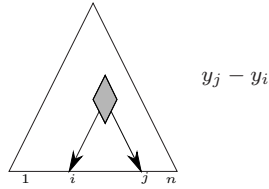
\begin{figure}[htbp]
\begin{center}
\setlength{\unitlength}{0.00083333in}
\begingroup\makeatletter\ifx\SetFigFont\undefined%
\gdef\SetFigFont#1#2#3#4#5{%
  \reset@font\fontsize{#1}{#2pt}%
  \fontfamily{#3}\fontseries{#4}\fontshape{#5}%
  \selectfont}%
\fi\endgroup%
{\renewcommand{\dashlinestretch}{30}
\begin{picture}(1688,1179)(0,-10)
\put(1212,627){\makebox(0,0)[lb]{{\SetFigFont{8}{9.6}{\rmdefault}{\mddefault}{\updefault}$y_j-y_i$}}}
\path(537,1152)(1062,102)(12,102)(537,1152)
\path(612,552)(387,102)
\blacken\path(413.833,222.748)(387.000,102.000)(467.498,195.915)(424.566,177.132)(413.833,222.748)
\path(612,552)(837,102)
\blacken\path(756.502,195.915)(837.000,102.000)(810.167,222.748)(799.434,177.132)(756.502,195.915)
\put(87,27){\makebox(0,0)[lb]{{\SetFigFont{5}{6.0}{\rmdefault}{\mddefault}{\updefault}$1$}}}
\put(387,27){\makebox(0,0)[lb]{{\SetFigFont{5}{6.0}{\rmdefault}{\mddefault}{\updefault}$i$}}}
\put(837,27){\makebox(0,0)[lb]{{\SetFigFont{5}{6.0}{\rmdefault}{\mddefault}{\updefault}$j$}}}
\put(987,27){\makebox(0,0)[lb]{{\SetFigFont{5}{6.0}{\rmdefault}{\mddefault}{\updefault}$n$}}}
\texture{88555555 55000000 555555 55000000 555555 55000000 555555 55000000 
	555555 55000000 555555 55000000 555555 55000000 555555 55000000 
	555555 55000000 555555 55000000 555555 55000000 555555 55000000 
	555555 55000000 555555 55000000 555555 55000000 555555 55000000 }
\shade\path(537,552)(612,702)(687,552)
	(612,402)(537,552)
\path(537,552)(612,702)(687,552)
	(612,402)(537,552)
\end{picture}
}

\end{center}
\caption{The weight contributed by an equivariant puzzle piece.  See
  the second panel of Figure~\ref{enex} for a specific
  example.\lremind{eqweight}}
\label{eqweight}
\end{figure}

Knutson and Tao gave an elegant equivariant Littlewood-Richardson rule
in terms of their puzzles.  Equivariant puzzles have the same pieces
as puzzles for ordinary cohomology (Fig.~\ref{puzzleD}), plus an additional
piece, shown in Figure~\ref{eqpuzzle}.  It may {\em not}
be rotated or reflected.  Each puzzle appears with a certain weight,
which is a product of contributions from the equivariant pieces in the
puzzle.  The contribution of the piece shown in Figure~\ref{eqweight},
which projects to position $i$ in the southwest direction and to
position $j$ in the southeast direction, is $y_j-y_i$. 

\begin{figure}[htbp]
\begin{center}
\setlength{\unitlength}{0.00083333in}
\begingroup\makeatletter\ifx\SetFigFont\undefined%
\gdef\SetFigFont#1#2#3#4#5{%
  \reset@font\fontsize{#1}{#2pt}%
  \fontfamily{#3}\fontseries{#4}\fontshape{#5}%
  \selectfont}%
\fi\endgroup%
{\renewcommand{\dashlinestretch}{30}
\begin{picture}(2788,1017)(0,-10)
\put(135,48){\makebox(0,0)[lb]{{\SetFigFont{8}{9.6}{\rmdefault}{\mddefault}{\updefault}$1$}}}
\path(12,90)(462,990)(912,90)(12,90)
\path(162,390)(762,390)
\path(462,390)(312,90)
\path(1812,90)(2262,990)(2712,90)(1812,90)
\path(612,90)(762,390)
\path(1962,390)(2112,90)(2412,690)
	(2112,690)(2412,90)(2562,390)
\texture{88555555 55000000 555555 55000000 555555 55000000 555555 55000000 
	555555 55000000 555555 55000000 555555 55000000 555555 55000000 
	555555 55000000 555555 55000000 555555 55000000 555555 55000000 
	555555 55000000 555555 55000000 555555 55000000 555555 55000000 }
\shade\path(2262,390)(2412,690)(2562,390)
	(2412,90)(2262,390)
\path(2262,390)(2412,690)(2562,390)
	(2412,90)(2262,390)
\put(491,809){\makebox(0,0)[lb]{{\SetFigFont{8}{9.6}{\rmdefault}{\mddefault}{\updefault}$0$}}}
\put(238,350){\makebox(0,0)[lb]{{\SetFigFont{8}{9.6}{\rmdefault}{\mddefault}{\updefault}$1$}}}
\put(546,354){\makebox(0,0)[lb]{{\SetFigFont{8}{9.6}{\rmdefault}{\mddefault}{\updefault}$0$}}}
\put(181,493){\makebox(0,0)[lb]{{\SetFigFont{8}{9.6}{\rmdefault}{\mddefault}{\updefault}$1$}}}
\put(346,495){\makebox(0,0)[lb]{{\SetFigFont{8}{9.6}{\rmdefault}{\mddefault}{\updefault}$1$}}}
\put(439,648){\makebox(0,0)[lb]{{\SetFigFont{8}{9.6}{\rmdefault}{\mddefault}{\updefault}$0$}}}
\put(641,510){\makebox(0,0)[lb]{{\SetFigFont{8}{9.6}{\rmdefault}{\mddefault}{\updefault}$1$}}}
\put(808,205){\makebox(0,0)[lb]{{\SetFigFont{8}{9.6}{\rmdefault}{\mddefault}{\updefault}$0$}}}
\put(640,200){\makebox(0,0)[lb]{{\SetFigFont{8}{9.6}{\rmdefault}{\mddefault}{\updefault}$0$}}}
\put(731,43){\makebox(0,0)[lb]{{\SetFigFont{8}{9.6}{\rmdefault}{\mddefault}{\updefault}$0$}}}
\put(2604,201){\makebox(0,0)[lb]{{\SetFigFont{8}{9.6}{\rmdefault}{\mddefault}{\updefault}$0$}}}
\put(2450,205){\makebox(0,0)[lb]{{\SetFigFont{8}{9.6}{\rmdefault}{\mddefault}{\updefault}$0$}}}
\put(2305,203){\makebox(0,0)[lb]{{\SetFigFont{8}{9.6}{\rmdefault}{\mddefault}{\updefault}$1$}}}
\put(2230,43){\makebox(0,0)[lb]{{\SetFigFont{8}{9.6}{\rmdefault}{\mddefault}{\updefault}$1$}}}
\put(2530,39){\makebox(0,0)[lb]{{\SetFigFont{8}{9.6}{\rmdefault}{\mddefault}{\updefault}$0$}}}
\put(1939,31){\makebox(0,0)[lb]{{\SetFigFont{8}{9.6}{\rmdefault}{\mddefault}{\updefault}$0$}}}
\put(1850,189){\makebox(0,0)[lb]{{\SetFigFont{8}{9.6}{\rmdefault}{\mddefault}{\updefault}$0$}}}
\put(2005,187){\makebox(0,0)[lb]{{\SetFigFont{8}{9.6}{\rmdefault}{\mddefault}{\updefault}$0$}}}
\put(2163,191){\makebox(0,0)[lb]{{\SetFigFont{8}{9.6}{\rmdefault}{\mddefault}{\updefault}$1$}}}
\put(1990,471){\makebox(0,0)[lb]{{\SetFigFont{8}{9.6}{\rmdefault}{\mddefault}{\updefault}$1$}}}
\put(2164,465){\makebox(0,0)[lb]{{\SetFigFont{8}{9.6}{\rmdefault}{\mddefault}{\updefault}$0$}}}
\put(2471,484){\makebox(0,0)[lb]{{\SetFigFont{8}{9.6}{\rmdefault}{\mddefault}{\updefault}$1$}}}
\put(2309,483){\makebox(0,0)[lb]{{\SetFigFont{8}{9.6}{\rmdefault}{\mddefault}{\updefault}$0$}}}
\put(2235,643){\makebox(0,0)[lb]{{\SetFigFont{8}{9.6}{\rmdefault}{\mddefault}{\updefault}$0$}}}
\put(2167,813){\makebox(0,0)[lb]{{\SetFigFont{8}{9.6}{\rmdefault}{\mddefault}{\updefault}$0$}}}
\put(2304,811){\makebox(0,0)[lb]{{\SetFigFont{8}{9.6}{\rmdefault}{\mddefault}{\updefault}$0$}}}
\put(357,199){\makebox(0,0)[lb]{{\SetFigFont{8}{9.6}{\rmdefault}{\mddefault}{\updefault}$0$}}}
\put(65,210){\makebox(0,0)[lb]{{\SetFigFont{8}{9.6}{\rmdefault}{\mddefault}{\updefault}$0$}}}
\put(496,200){\makebox(0,0)[lb]{{\SetFigFont{8}{9.6}{\rmdefault}{\mddefault}{\updefault}$0$}}}
\put(357,806){\makebox(0,0)[lb]{{\SetFigFont{8}{9.6}{\rmdefault}{\mddefault}{\updefault}$0$}}}
\put(434,47){\makebox(0,0)[lb]{{\SetFigFont{8}{9.6}{\rmdefault}{\mddefault}{\updefault}$0$}}}
\path(612,690)(312,690)(612,90)
\end{picture}
}

\end{center}
\caption{Example of the two equivariant puzzles associated to a single
  equivariant Schubert problem $\si_{010} \cdot \si_{010}$.  The first
  already appears in ordinary cohomology.  The second appears with
  weight $y_3-y_2$ (notice the equivariant rhombus in the southeast of
  the figure).  The $10$ and $01$ edges are omitted for the sake of
  clarity. This verifies equ.\ \eqref{basicHT}:  $\si_{010}^2 = \si_{100} + (y_3-y_2) \si_{010}$. \lremind{enex}}
\label{enex}
\end{figure}
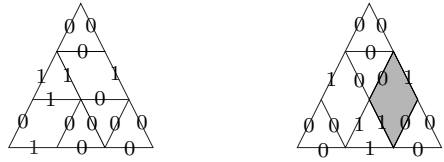

A specific example is given in Figure~\ref{enex}, discussed at length
in \S \ref{eqex}.  This example corresponds to two Schubert condition
corresponding to those points lying on the fixed line $\proj \langle
\vv_1, \vv_2 \rangle$.

We now give a conjectural geometric interpretation to this
combinatorial rule, due to Knutson and the second author, extending
the geometric Littlewood-Richardson rule in ordinary cohomology.  The
torus action ``prevents'' the equivariant cycle from degenerating (the
degeneration of \S \ref{hhh} is not equivariant), so we cannot use our
degeneration interpretation.  As in the example of Figure~\ref{enex},
these Schubert problems are excess intersection problems: the Schubert
cycles $\overline{\Omega}_{\alpha}(F_{\bullet})$ and
$\overline{\Omega}_{\beta}(F_{\bullet})$ in $G(k,n)$ do not intersect
properly (in the ``expected'' dimension).  So we use instead the
following trick.

The degeneration order may be interpreted in the flag variety $Fl(n)$
parametrizing the moving flag's relative position to the fixed flag.
It corresponds to a sequence of nested Schubert varieties in $Fl(n)$,
each a Cartier divisor on the previous one, where the first element is
all of $Fl(n)$, and the last is a point.  (This Cartier requirement is
a strong constraint on the degeneration order, equivalent to the fact
that each Schubert variety in the sequence  is smooth.)  We consider subvarieties of
$G(k,n) \times Fl(n)$, where $Fl(n)$ parametrizes the ``moving'' flag
$M_{\bullet}$, of the form \lremind{eqPV}
\begin{equation}
\pi_1^* \left( \overline{ \Omega_{\al}(M_{\bullet}) \cap
    \Omega_{\be}(F_{\bullet}) } \right) \cap \pi_2^*
\overline{\Omega}_{d_i}(F_{\bullet})\label{eqPV}
\end{equation}
where $d_i$ is in the degeneration
order (\S \ref{degenerationorder}).  In other words, $M_{\bullet}$ is required to be in given relative
position to $F_{\bullet}$ (or in a more degenerate position), and the
$k$-plane is required to be in given relative position to both
$M_{\bullet}$ and $F_{\bullet}$.  We will see that $\al$ and $\be$ may force
the fixed and moving flags to be in more degenerate position than that
required by $d_i$; an explicit example will be given in \S \ref{eqex}.

Call such a variety an {\em equivariant puzzle variety}.  Note that
the fiber of the equivariant puzzle variety \eqref{eqPV} over a point
of the (open) Schubert cell $\Omega_{d_i}(F_{\bullet})$ is a puzzle
variety \eqref{puzzlevariety}.  We will shortly associate such a
variety to a partially completed equivariant puzzle.  Note that if
$d_i$ is the final step in the degeneration order ($i= \binom n 2$),
the equivariant puzzle variety is a Schubert variety with respect to
the fixed flag $F_\bullet$.

Now pull back (to a given equivariant puzzle variety) the divisor
$\oO_{d_{i+1}}(F_{\bullet})$ on $\oO_{d_i}(F_{\bullet})$ corresponding
to the next Schubert variety in the degeneration order.  We conjecture
that either (i) this divisor does not contain our equivariant puzzle
variety (i.e.\ it pulls back to a Cartier divisor on the equivariant
puzzle variety), in which case it is (scheme-theoretically) the union
of one or two other equivariant puzzle varieties (corresponding to
$d_{i+1}$), or (ii) it contains the equivariant puzzle variety.  In
the latter case, the excess intersection problem (intersection with
the zero-section of a line bundle) is simplest sort of excess
intersection: we obtain a contribution of the equivariant first Chern
class of the line bundle, which is pure weight.  (This is where the
Cartier hypothesis is used.)

We now make this precise, and conclude with an example, which will
perhaps be most helpful to the reader.

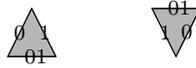
\begin{figure}[htbp]
\begin{center}
\setlength{\unitlength}{0.00083333in}
\begingroup\makeatletter\ifx\SetFigFont\undefined%
\gdef\SetFigFont#1#2#3#4#5{%
  \reset@font\fontsize{#1}{#2pt}%
  \fontfamily{#3}\fontseries{#4}\fontshape{#5}%
  \selectfont}%
\fi\endgroup%
{\renewcommand{\dashlinestretch}{30}
\begin{picture}(1280,435)(0,-10)
\texture{88555555 55000000 555555 55000000 555555 55000000 555555 55000000 
	555555 55000000 555555 55000000 555555 55000000 555555 55000000 
	555555 55000000 555555 55000000 555555 55000000 555555 55000000 
	555555 55000000 555555 55000000 555555 55000000 555555 55000000 }
\shade\path(1062,74)(912,374)(1212,374)(1062,74)
\path(1062,74)(912,374)(1212,374)(1062,74)
\shade\path(12,74)(162,374)(312,74)(12,74)
\path(12,74)(162,374)(312,74)(12,74)
\put(209,180){\makebox(0,0)[lb]{{\SetFigFont{8}{9.6}{\rmdefault}{\mddefault}{\updefault}$1$}}}
\put(117,31){\makebox(0,0)[lb]{{\SetFigFont{8}{9.6}{\rmdefault}{\mddefault}{\updefault}$01$}}}
\put(959,181){\makebox(0,0)[lb]{{\SetFigFont{8}{9.6}{\rmdefault}{\mddefault}{\updefault}$1$}}}
\put(1096,185){\makebox(0,0)[lb]{{\SetFigFont{8}{9.6}{\rmdefault}{\mddefault}{\updefault}$0$}}}
\put(1012,336){\makebox(0,0)[lb]{{\SetFigFont{8}{9.6}{\rmdefault}{\mddefault}{\updefault}$01$}}}
\put(52,180){\makebox(0,0)[lb]{{\SetFigFont{8}{9.6}{\rmdefault}{\mddefault}{\updefault}$0$}}}
\end{picture}
}

\end{center}
\caption{An alternate equivariant puzzle piece.  It may appear only in
  the orientations shown here; it may {\em not} be further rotated or
  reflected.  This is a trivial variation of Figure~\ref{eqpuzzle}.
\lremind{eqpuzzletwo}}
\label{eqpuzzletwo}
\end{figure}

Here is how to associate an equivariant puzzle variety to a partially
completed equivariant puzzle.  In order to parallel our earlier
description, we replace Knutson and Tao's piece of
Figure~\ref{eqpuzzle} with the two halves of Figure~\ref{eqpuzzletwo}.
(This is clearly a trivial variation of Fig.~\ref{eqpuzzle}.)  To each
partially completed puzzle (as in Fig.~\ref{degenorder}), we define
$\al$ and $\be$ as follows.  The recipe is the same as in
Figures~\ref{readalpha} and ~\ref{readbeta}, with the additional fact
that $01$ is read as $0$ for $\al$ and $1$ for $\be$.
(This may also be interpreted in the same way as Figure~\ref{newab}.)  Then the
conjectural equivariant geometric Littlewood-Richardson rule is as
follows.  To compute the product of two equivariant Schubert classes,
we fill a puzzle in the degeneration order  with equivariant puzzle pieces.  For each move,
we fill in the next two or three puzzle pieces.  There will always be
one or two choices.  If we place a copy of the second piece in
Figure~\ref{eqpuzzletwo} (the bottom half of Knutson-Tao's
piece), there is one choice, and this will correspond to excess
intersection.  The weight is as given in Figure~\ref{eqweight}.
Otherwise, there is no excess intersection, and the Cartier divisor
will be the scheme-theoretic union of the equivariant puzzle varieties
corresponding to the one or two choices.

Thus to compute a product of two equivariant Schubert classes, we
start with an empty puzzle with labeled northwest and northeast edges,
and complete it.  The partially completed puzzles we see en route
completely describe the geometry of the successive Cartier slicing.
Note that this rule is manifestly positive (as it is simply an
interpretation of the puzzle rule).  It is also clearly a
generalization of the geometric Littlewood-Richardson rule in ordinary
cohomology.  It is well-checked (for example for $n$ up to $5$, and
for many cases for larger $n$).

\subsection{A worked equivariant example}
\lremind{eqex}\label{eqex}We conclude our equivariant discussion with
a worked example, which we hope will illustrate what is happening
geometrically.  Consider the problem of Figure~\ref{enex},
corresponding to intersecting the locus of points meeting the fixed
line with itself.  We will verify that\lremind{basicHT}
\begin{equation}\label{basicHT}
\si_{010}^2 = \si_{100} + (y_3-y_2) \si_{010}.
\end{equation}
In the Grassmannian $G(1,3)$, also known as $\proj^2$, this is the
intersection of a codimension $1$ class with itself, so we expect a
dimension $0$ answer.  In ordinary cohomology (which we recover by
setting the equivariant parameters to zero), we expect to see the
class of a point (in ordinary cohomology, lines can deform, and the
intersection of two general lines is one point).

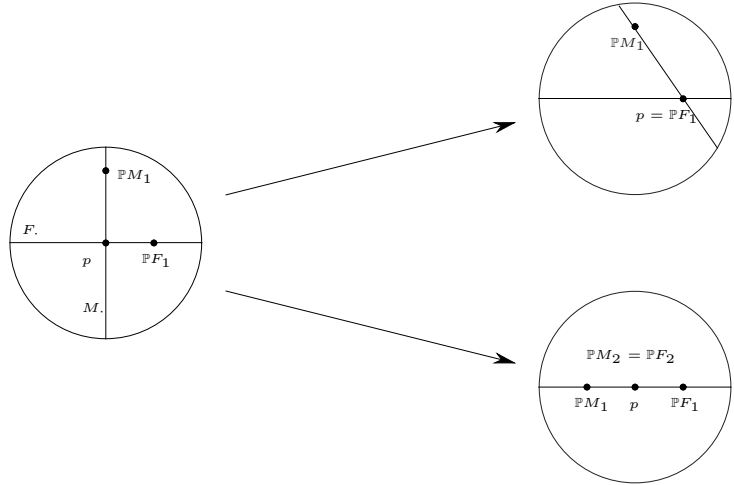
\begin{figure}[htbp]
\begin{center}
\setlength{\unitlength}{0.00083333in}
\begingroup\makeatletter\ifx\SetFigFont\undefined%
\gdef\SetFigFont#1#2#3#4#5{%
  \reset@font\fontsize{#1}{#2pt}%
  \fontfamily{#3}\fontseries{#4}\fontshape{#5}%
  \selectfont}%
\fi\endgroup%
{\renewcommand{\dashlinestretch}{30}
\begin{picture}(4524,3029)(0,-10)
\put(3762,2707){\makebox(0,0)[lb]{{\SetFigFont{5}{6.0}{\rmdefault}{\mddefault}{\updefault}$\proj M_1$}}}
\put(3912,607){\ellipse{1200}{1200}}
\put(3912,2407){\ellipse{1200}{1200}}
\put(912,1507){\blacken\ellipse{36}{36}}
\put(912,1507){\ellipse{36}{36}}
\put(4212,2407){\blacken\ellipse{36}{36}}
\put(4212,2407){\ellipse{36}{36}}
\put(4212,607){\blacken\ellipse{36}{36}}
\put(4212,607){\ellipse{36}{36}}
\put(612,1507){\blacken\ellipse{36}{36}}
\put(612,1507){\ellipse{36}{36}}
\put(3912,607){\blacken\ellipse{36}{36}}
\put(3912,607){\ellipse{36}{36}}
\put(3612,607){\blacken\ellipse{36}{36}}
\put(3612,607){\ellipse{36}{36}}
\put(612,1957){\blacken\ellipse{36}{36}}
\put(612,1957){\ellipse{36}{36}}
\put(3912,2857){\blacken\ellipse{36}{36}}
\put(3912,2857){\ellipse{36}{36}}
\path(1362,1207)(3162,757)
\blacken\path(3038.307,757.000)(3162.000,757.000)(3052.859,815.209)(3080.508,777.373)(3038.307,757.000)
\path(1362,1807)(3162,2257)
\blacken\path(3052.859,2198.791)(3162.000,2257.000)(3038.307,2257.000)(3080.508,2236.627)(3052.859,2198.791)
\path(12,1507)(1212,1507)
\path(3312,2407)(4512,2407)
\path(3312,607)(4512,607)
\path(612,2107)(612,907)
\path(4429,2099)(3815,2985)
\put(837,1357){\makebox(0,0)[lb]{{\SetFigFont{5}{6.0}{\rmdefault}{\mddefault}{\updefault}$\proj F_1$}}}
\put(4137,457){\makebox(0,0)[lb]{{\SetFigFont{5}{6.0}{\rmdefault}{\mddefault}{\updefault}$\proj F_1$}}}
\put(3612,757){\makebox(0,0)[lb]{{\SetFigFont{5}{6.0}{\rmdefault}{\mddefault}{\updefault}$\proj M_2 = \proj F_2$}}}
\put(687,1882){\makebox(0,0)[lb]{{\SetFigFont{5}{6.0}{\rmdefault}{\mddefault}{\updefault}$\proj M_1$}}}
\put(462,1357){\makebox(0,0)[lb]{{\SetFigFont{5}{6.0}{\rmdefault}{\mddefault}{\updefault}$p$}}}
\put(3912,2257){\makebox(0,0)[lb]{{\SetFigFont{5}{6.0}{\rmdefault}{\mddefault}{\updefault}$p=\proj F_1$}}}
\put(3876,457){\makebox(0,0)[lb]{{\SetFigFont{5}{6.0}{\rmdefault}{\mddefault}{\updefault}$p$}}}
\put(3537,457){\makebox(0,0)[lb]{{\SetFigFont{5}{6.0}{\rmdefault}{\mddefault}{\updefault}$\proj M_1$}}}
\put(462,1057){\makebox(0,0)[lb]{{\SetFigFont{5}{6.0}{\rmdefault}{\mddefault}{\updefault}$M_{\cdot}$}}}
\put(87,1546){\makebox(0,0)[lb]{{\SetFigFont{5}{6.0}{\rmdefault}{\mddefault}{\updefault}$F_{\cdot}$}}}
\put(612,1507){\ellipse{1200}{1200}}
\end{picture}
}

\end{center}
\caption{A worked example (\S \ref{eqex}) showing
a reducible Cartier slice.
\lremind{enexg}}
\label{enexg}
\end{figure}

Following the recipe, instead of working in $\proj^2$, we work in
$\proj^2 \times Fl(3) \cong \proj^2 \times \FF l (2)$, parametrizing
the moving line, the moving point (on the moving line), and the point
$p$ of $\proj^2$.  On this fivefold, we consider the threefold
corresponding to requiring the point $p$ to lie on both the moving
line $\proj M_2$ and on the fixed line $\proj F_2$.  The generic such
configuration is depicted on the left side of Figure~\ref{enexg}.

We now slice with our first divisor $D$, corresponding to requiring
the moving line $\proj M_2$ to pass through the fixed point $\proj
F_1$.  This Cartier divisor (a fourfold) has two components $D_1$
and $D_2$.

The first component $D_1$ is the geometrically clear one: the moving line
$\proj M_2$ rotates around the fixed point $\proj F_1$ (yielding one
dimension of moduli), $p$ is the fixed point $\proj F_1$, and the
moving point varies (yielding a second dimension of moduli).  This
corresponds to (the top two rows of) the first puzzle of
Figure~\ref{enex}.  The generic such configuration is depicted on the
top right side of Figure~\ref{enexg}.

But there is a second two-dimensional component $D_2$: the moving line
$\proj M_2$ equals the fixed line $\proj F_2$, $p$ varies on the
moving line (yielding one dimension of moduli), and the moving point
$\proj M_1$ varies on the moving line (yielding the second dimension
of moduli).  This component is not seen in the usual geometric
Littlewood-Richardson rule.  This corresponds to (the top two rows of) the
second puzzle of Figure~\ref{enex}.  The generic such configuration is
depicted on the bottom right side of Figure~\ref{enexg}.  Note that
the moving and fixed flags are forced to be in
more degenerate condition than required by the degeneration order:
the moving and fixed lines are forced to agree.

The reader may verify that there are no other components of $D$.

In order to see how the weight arises, we continue to follow the
second case.  In the next Cartier slice, the moving point is required
to agree with the fixed point ($\proj M_1 = \proj F_1$).  This is
clearly an irreducible divisor: the point $p$ may still vary on the
moving line $\proj M_2$, and the moving point is now fixed (despite
its name).  We then slice with the pullback of the next Cartier
divisor ($\oO_{d_3}(F_\bullet) \subset \oO_{d_2}(F_{\bullet})$),
requiring the moving line to agree with the fixed line.  However, this
is {\em not} a divisor on our equivariant puzzle variety: the
condition of the moving line agreeing with the fixed line is satisfied
by the entire equivariant puzzle variety.  Thus we get excess
intersection given by the (equivariant) first Chern class of this line
bundle corresponding to the divisor, which is readily checked to
be $y_3-y_2$.

\section{The equivariant $K$-theory of the Grassmannian}
\lremind{KT}\label{KT}The persevering reader who has read
Sections~\ref{K} and \ref{HT} will realize that this begs for an extension
to equivariant $K$-theory, and that this would ideally correspond to
the following: in a Cartier slice in the equivariant
discussion where the variety breaks into two pieces, the
scheme-theoretic intersection of the two pieces should be another
equivariant puzzle variety.  In this case we do not have pre-existing
puzzles to guide us.

Knutson and the second author indeed conjecture such a rule.  We first
describe the new equivariant $K$-theory puzzle pieces, which gives a
purely combinatorial conjectural rule in equivariant $K$-theory, and
then describe the geometry conjecturally associated to it.  Knutson
and the first author have verified this rule up to dimension $5$.  We
emphasize that although this may be interpreted as a purely
combinatorial rule, it was induced from many geometric examples.

The equivariant $K$-theory pieces include the original pieces from
equivariant cohomology --- those from ordinary cohomology shown in
Figure~\ref{puzzleD}, and the equivariant piece of Figure~\ref{eqpuzzle}.
The equivariant piece in the position shown in Figure~\ref{eqweight}
now contributes $1 - e^{y_i-y_j}$.  (The geometric reason
is that this is will be the Chern class
of a line bundle
in equivariant $K$-theory, not in equivariant cohomology.)

\begin{figure}[htbp]
\begin{center}
\setlength{\unitlength}{0.00083333in}
\begingroup\makeatletter\ifx\SetFigFont\undefined%
\gdef\SetFigFont#1#2#3#4#5{%
  \reset@font\fontsize{#1}{#2pt}%
  \fontfamily{#3}\fontseries{#4}\fontshape{#5}%
  \selectfont}%
\fi\endgroup%
{\renewcommand{\dashlinestretch}{30}
\begin{picture}(1458,670)(0,-10)
\put(1274,436){\makebox(0,0)[lb]{{\SetFigFont{8}{9.6}{\rmdefault}{\mddefault}{\updefault}$0$}}}
\path(1212,312)(1362,612)(1062,612)(1362,12)
\put(162,462){\makebox(0,0)[lb]{{\SetFigFont{5}{6.0}{\rmdefault}{\mddefault}{\updefault}$1$}}}
\put(12,162){\makebox(0,0)[lb]{{\SetFigFont{5}{6.0}{\rmdefault}{\mddefault}{\updefault}$0$}}}
\put(1062,387){\makebox(0,0)[lb]{{\SetFigFont{5}{6.0}{\rmdefault}{\mddefault}{\updefault}$1$}}}
\put(1212,87){\makebox(0,0)[lb]{{\SetFigFont{5}{6.0}{\rmdefault}{\mddefault}{\updefault}$0$}}}
\put(112,90){\makebox(0,0)[lb]{{\SetFigFont{5}{6.0}{\rmdefault}{\mddefault}{\updefault}$1$}}}
\put(265,392){\makebox(0,0)[lb]{{\SetFigFont{5}{6.0}{\rmdefault}{\mddefault}{\updefault}$0$}}}
\put(1159,462){\makebox(0,0)[lb]{{\SetFigFont{5}{6.0}{\rmdefault}{\mddefault}{\updefault}$0$}}}
\put(1307,162){\makebox(0,0)[lb]{{\SetFigFont{5}{6.0}{\rmdefault}{\mddefault}{\updefault}$1$}}}
\put(138,570){\makebox(0,0)[lb]{{\SetFigFont{8}{9.6}{\rmdefault}{\mddefault}{\updefault}$1$}}}
\put(49,427){\makebox(0,0)[lb]{{\SetFigFont{8}{9.6}{\rmdefault}{\mddefault}{\updefault}$1$}}}
\put(1183,571){\makebox(0,0)[lb]{{\SetFigFont{8}{9.6}{\rmdefault}{\mddefault}{\updefault}$0$}}}
\path(162,312)(12,612)(312,612)(12,12)
\end{picture}
}

\end{center}
\caption{The new  (conjectural) $K_T^*$ puzzle pieces.\lremind{KTpieces}}
\label{KTpieces}
\end{figure}
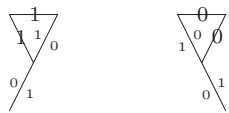

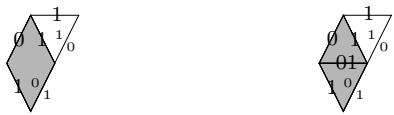
\begin{figure}[htbp]
\begin{center}
\setlength{\unitlength}{0.00083333in}
\begingroup\makeatletter\ifx\SetFigFont\undefined%
\gdef\SetFigFont#1#2#3#4#5{%
  \reset@font\fontsize{#1}{#2pt}%
  \fontfamily{#3}\fontseries{#4}\fontshape{#5}%
  \selectfont}%
\fi\endgroup%
{\renewcommand{\dashlinestretch}{30}
\begin{picture}(2429,680)(0,-10)
\texture{88555555 55000000 555555 55000000 555555 55000000 555555 55000000 
	555555 55000000 555555 55000000 555555 55000000 555555 55000000 
	555555 55000000 555555 55000000 555555 55000000 555555 55000000 
	555555 55000000 555555 55000000 555555 55000000 555555 55000000 }
\shade\path(1962,312)(2112,612)(2262,312)(1962,312)
\path(1962,312)(2112,612)(2262,312)(1962,312)
\shade\path(12,312)(162,612)(312,312)
	(162,12)(12,312)
\path(12,312)(162,612)(312,312)
	(162,12)(12,312)
\shade\path(2262,312)(2112,12)(1962,312)(2262,312)
\path(2262,312)(2112,12)(1962,312)(2262,312)
\path(312,312)(162,612)(462,612)
	(162,12)(12,312)(162,612)
\put(2337,387){\makebox(0,0)[lb]{{\SetFigFont{5}{6.0}{\rmdefault}{\mddefault}{\updefault}$0$}}}
\put(162,162){\makebox(0,0)[lb]{{\SetFigFont{5}{6.0}{\rmdefault}{\mddefault}{\updefault}$0$}}}
\put(387,387){\makebox(0,0)[lb]{{\SetFigFont{5}{6.0}{\rmdefault}{\mddefault}{\updefault}$0$}}}
\put(2112,162){\makebox(0,0)[lb]{{\SetFigFont{5}{6.0}{\rmdefault}{\mddefault}{\updefault}$0$}}}
\put(2262,462){\makebox(0,0)[lb]{{\SetFigFont{5}{6.0}{\rmdefault}{\mddefault}{\updefault}$1$}}}
\put(2187,87){\makebox(0,0)[lb]{{\SetFigFont{5}{6.0}{\rmdefault}{\mddefault}{\updefault}$1$}}}
\put(312,462){\makebox(0,0)[lb]{{\SetFigFont{5}{6.0}{\rmdefault}{\mddefault}{\updefault}$1$}}}
\put(237,87){\makebox(0,0)[lb]{{\SetFigFont{5}{6.0}{\rmdefault}{\mddefault}{\updefault}$1$}}}
\put(290,574){\makebox(0,0)[lb]{{\SetFigFont{8}{9.6}{\rmdefault}{\mddefault}{\updefault}$1$}}}
\put(198,423){\makebox(0,0)[lb]{{\SetFigFont{8}{9.6}{\rmdefault}{\mddefault}{\updefault}$1$}}}
\put(54,422){\makebox(0,0)[lb]{{\SetFigFont{8}{9.6}{\rmdefault}{\mddefault}{\updefault}$0$}}}
\put(50,126){\makebox(0,0)[lb]{{\SetFigFont{8}{9.6}{\rmdefault}{\mddefault}{\updefault}$1$}}}
\put(2236,581){\makebox(0,0)[lb]{{\SetFigFont{8}{9.6}{\rmdefault}{\mddefault}{\updefault}$1$}}}
\put(2151,423){\makebox(0,0)[lb]{{\SetFigFont{8}{9.6}{\rmdefault}{\mddefault}{\updefault}$1$}}}
\put(2009,427){\makebox(0,0)[lb]{{\SetFigFont{8}{9.6}{\rmdefault}{\mddefault}{\updefault}$0$}}}
\put(2065,276){\makebox(0,0)[lb]{{\SetFigFont{8}{9.6}{\rmdefault}{\mddefault}{\updefault}$01$}}}
\path(2262,312)(2112,612)(2412,612)(2112,12)
\put(2007,122){\makebox(0,0)[lb]{{\SetFigFont{8}{9.6}{\rmdefault}{\mddefault}{\updefault}$1$}}}
\path(2262,312)(1962,312)(2112,612)
\path(1962,312)(2112,12)
\end{picture}
}

\end{center}
\caption{Placing the first $K_T$-piece.\lremind{KTplaceone}
The left description is if one is using the equivariant
piece of Figure~\ref{eqpuzzle}, and the right is if one
is using the alternate equivariant piece of Figure~\ref{eqpuzzletwo}.}
\label{KTplaceone}
\end{figure}

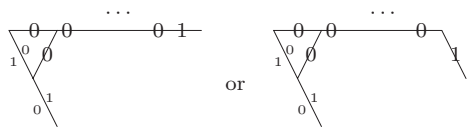
\begin{figure}[htbp]
\begin{center}
\setlength{\unitlength}{0.00083333in}
\begingroup\makeatletter\ifx\SetFigFont\undefined%
\gdef\SetFigFont#1#2#3#4#5{%
  \reset@font\fontsize{#1}{#2pt}%
  \fontfamily{#3}\fontseries{#4}\fontshape{#5}%
  \selectfont}%
\fi\endgroup%
{\renewcommand{\dashlinestretch}{30}
\begin{picture}(2948,786)(0,-10)
\put(2764,419){\makebox(0,0)[lb]{{\SetFigFont{8}{9.6}{\rmdefault}{\mddefault}{\updefault}$1$}}}
\path(312,612)(1062,612)
\path(1812,312)(1962,612)(1662,612)(1962,12)
\path(1962,612)(2712,612)
\path(1062,612)(1212,612)
\path(2712,612)(2862,312)
\put(87,462){\makebox(0,0)[lb]{{\SetFigFont{5}{6.0}{\rmdefault}{\mddefault}{\updefault}$0$}}}
\put(12,387){\makebox(0,0)[lb]{{\SetFigFont{5}{6.0}{\rmdefault}{\mddefault}{\updefault}$1$}}}
\put(237,162){\makebox(0,0)[lb]{{\SetFigFont{5}{6.0}{\rmdefault}{\mddefault}{\updefault}$1$}}}
\put(162,87){\makebox(0,0)[lb]{{\SetFigFont{5}{6.0}{\rmdefault}{\mddefault}{\updefault}$0$}}}
\put(612,687){\makebox(0,0)[lb]{{\SetFigFont{8}{9.6}{\rmdefault}{\mddefault}{\updefault}$\cdots$}}}
\put(1737,462){\makebox(0,0)[lb]{{\SetFigFont{5}{6.0}{\rmdefault}{\mddefault}{\updefault}$0$}}}
\put(1662,387){\makebox(0,0)[lb]{{\SetFigFont{5}{6.0}{\rmdefault}{\mddefault}{\updefault}$1$}}}
\put(1887,162){\makebox(0,0)[lb]{{\SetFigFont{5}{6.0}{\rmdefault}{\mddefault}{\updefault}$1$}}}
\put(1812,87){\makebox(0,0)[lb]{{\SetFigFont{5}{6.0}{\rmdefault}{\mddefault}{\updefault}$0$}}}
\put(2262,687){\makebox(0,0)[lb]{{\SetFigFont{8}{9.6}{\rmdefault}{\mddefault}{\updefault}$\cdots$}}}
\put(1362,237){\makebox(0,0)[lb]{{\SetFigFont{8}{9.6}{\rmdefault}{\mddefault}{\updefault}or}}}
\put(216,422){\makebox(0,0)[lb]{{\SetFigFont{8}{9.6}{\rmdefault}{\mddefault}{\updefault}$0$}}}
\put(135,570){\makebox(0,0)[lb]{{\SetFigFont{8}{9.6}{\rmdefault}{\mddefault}{\updefault}$0$}}}
\put(908,570){\makebox(0,0)[lb]{{\SetFigFont{8}{9.6}{\rmdefault}{\mddefault}{\updefault}$0$}}}
\put(1058,570){\makebox(0,0)[lb]{{\SetFigFont{8}{9.6}{\rmdefault}{\mddefault}{\updefault}$1$}}}
\put(1785,569){\makebox(0,0)[lb]{{\SetFigFont{8}{9.6}{\rmdefault}{\mddefault}{\updefault}$0$}}}
\put(1861,421){\makebox(0,0)[lb]{{\SetFigFont{8}{9.6}{\rmdefault}{\mddefault}{\updefault}$0$}}}
\put(1989,570){\makebox(0,0)[lb]{{\SetFigFont{8}{9.6}{\rmdefault}{\mddefault}{\updefault}$0$}}}
\put(340,570){\makebox(0,0)[lb]{{\SetFigFont{8}{9.6}{\rmdefault}{\mddefault}{\updefault}$0$}}}
\put(2549,570){\makebox(0,0)[lb]{{\SetFigFont{8}{9.6}{\rmdefault}{\mddefault}{\updefault}$0$}}}
\path(162,312)(312,612)(12,612)(312,12)
\end{picture}
}

\end{center}
\caption{Placing the second $K_T$-piece.\lremind{KTplace}}
\label{KTplace}
\end{figure}

We do {\em not} use the $K$-theory pieces of Figs.~11 or~12.  There
are instead two new pieces, shown in Fig.~\ref{KTpieces}, each
contributing a sign of $-1$.  Each $K_T$-piece has an unusual pair of
edges, each of which is $0$ on one side and $1$ on the other.  Readers
familiar with \cite{knutson:puzzle} will recognize these pairs of
edges as {\em gashes}.  These pieces must be placed in this
orientation (i.e.\ may not be reflected or rotated). Another new
feature is that there are constraints on where the two pieces may be
placed.  The first piece must be placed to the right of an equivariant
rhombus (shown in Figure~\ref{KTplaceone}).  The second piece may only
be placed (when completing the puzzle from top to bottom and left to
right as usual) if the edges to its right are a (possibly empty)
series of horizontal $0$'s followed by a $1$.  This is depicted in
Figure ~\ref{KTplace}.  This latter condition is a new, potentially
disturbing, nonlocality in filling a puzzle (although it is at least
local to the row).  However, it is dictated by the conjectural
geometry.

\begin{figure}[htbp]
\begin{center}
\setlength{\unitlength}{0.00083333in}
\begingroup\makeatletter\ifx\SetFigFont\undefined%
\gdef\SetFigFont#1#2#3#4#5{%
  \reset@font\fontsize{#1}{#2pt}%
  \fontfamily{#3}\fontseries{#4}\fontshape{#5}%
  \selectfont}%
\fi\endgroup%
{\renewcommand{\dashlinestretch}{30}
\begin{picture}(4725,2308)(0,-10)
\texture{88555555 55000000 555555 55000000 555555 55000000 555555 55000000 
	555555 55000000 555555 55000000 555555 55000000 555555 55000000 
	555555 55000000 555555 55000000 555555 55000000 555555 55000000 
	555555 55000000 555555 55000000 555555 55000000 555555 55000000 }
\shade\path(3450,1381)(3300,1681)(3450,1981)
	(3600,1681)(3450,1381)
\path(3450,1381)(3300,1681)(3450,1981)
	(3600,1681)(3450,1381)
\path(1350,1231)(1050,1831)(1350,1831)
	(1200,1531)(900,1531)
\blacken\path(1194.502,574.915)(1275.000,481.000)(1248.167,601.748)(1237.434,556.132)(1194.502,574.915)
\dashline{60.000}(1275,481)(825,1381)
\path(1575,1381)(1050,331)
\blacken\path(1076.833,451.748)(1050.000,331.000)(1130.498,424.915)(1087.566,406.132)(1076.833,451.748)
\dashline{60.000}(1050,1531)(1575,481)
\blacken\path(1494.502,574.915)(1575.000,481.000)(1548.167,601.748)(1537.434,556.132)(1494.502,574.915)
\dashline{60.000}(1275,1681)(1875,481)
\blacken\path(1794.502,574.915)(1875.000,481.000)(1848.167,601.748)(1837.434,556.132)(1794.502,574.915)
\dashline{60.000}(1050,1531)(450,331)
\blacken\path(476.833,451.748)(450.000,331.000)(530.498,424.915)(487.566,406.132)(476.833,451.748)
\dashline{60.000}(1425,1681)(750,331)
\blacken\path(776.833,451.748)(750.000,331.000)(830.498,424.915)(787.566,406.132)(776.833,451.748)
\path(675,331)(825,181)
\path(825,331)(675,181)
\path(3150,1381)(3600,2281)(4050,1381)
\path(3300,1681)(3450,1381)(3750,1981)
	(3450,1981)(3600,1681)
\dashline{60.000}(3225,1531)(3750,481)
\blacken\path(3669.502,574.915)(3750.000,481.000)(3723.167,601.748)(3712.434,556.132)(3669.502,574.915)
\dashline{60.000}(3525,1531)(4050,481)
\blacken\path(3969.502,574.915)(4050.000,481.000)(4023.167,601.748)(4012.434,556.132)(3969.502,574.915)
\dashline{60.000}(3675,1831)(4350,481)
\blacken\path(4269.502,574.915)(4350.000,481.000)(4323.167,601.748)(4312.434,556.132)(4269.502,574.915)
\dashline{60.000}(3375,1531)(2700,181)
\blacken\path(2726.833,301.748)(2700.000,181.000)(2780.498,274.915)(2737.566,256.132)(2726.833,301.748)
\dashline{60.000}(3825,1831)(3000,181)
\blacken\path(3026.833,301.748)(3000.000,181.000)(3080.498,274.915)(3037.566,256.132)(3026.833,301.748)
\dashline{60.000}(3975,1531)(3300,181)
\blacken\path(3326.833,301.748)(3300.000,181.000)(3380.498,274.915)(3337.566,256.132)(3326.833,301.748)
\path(750,1231)(1650,1231)
\path(3150,1381)(4050,1381)
\path(750,1231)(1200,2131)(1650,1231)
\put(1125,1681){\makebox(0,0)[lb]{{\SetFigFont{5}{6.0}{\rmdefault}{\mddefault}{\updefault}$0$}}}
\put(1050,1681){\makebox(0,0)[lb]{{\SetFigFont{5}{6.0}{\rmdefault}{\mddefault}{\updefault}$1$}}}
\put(1275,1381){\makebox(0,0)[lb]{{\SetFigFont{5}{6.0}{\rmdefault}{\mddefault}{\updefault}$1$}}}
\put(1200,1306){\makebox(0,0)[lb]{{\SetFigFont{5}{6.0}{\rmdefault}{\mddefault}{\updefault}$0$}}}
\put(1275,331){\makebox(0,0)[lb]{{\SetFigFont{8}{9.6}{\rmdefault}{\mddefault}{\updefault}$0$}}}
\put(1575,331){\makebox(0,0)[lb]{{\SetFigFont{8}{9.6}{\rmdefault}{\mddefault}{\updefault}$1$}}}
\put(1875,331){\makebox(0,0)[lb]{{\SetFigFont{8}{9.6}{\rmdefault}{\mddefault}{\updefault}$0$}}}
\put(375,181){\makebox(0,0)[lb]{{\SetFigFont{8}{9.6}{\rmdefault}{\mddefault}{\updefault}$1$}}}
\put(675,181){\makebox(0,0)[lb]{{\SetFigFont{8}{9.6}{\rmdefault}{\mddefault}{\updefault}$1$}}}
\put(975,181){\makebox(0,0)[lb]{{\SetFigFont{8}{9.6}{\rmdefault}{\mddefault}{\updefault}$0$}}}
\put(675,31){\makebox(0,0)[lb]{{\SetFigFont{8}{9.6}{\rmdefault}{\mddefault}{\updefault}$0$}}}
\put(2175,331){\makebox(0,0)[lb]{{\SetFigFont{8}{9.6}{\rmdefault}{\mddefault}{\updefault}$=\al$}}}
\put(0,181){\makebox(0,0)[lb]{{\SetFigFont{8}{9.6}{\rmdefault}{\mddefault}{\updefault}$\be=$}}}
\put(3675,1756){\makebox(0,0)[lb]{{\SetFigFont{5}{6.0}{\rmdefault}{\mddefault}{\updefault}$0$}}}
\put(3525,1456){\makebox(0,0)[lb]{{\SetFigFont{5}{6.0}{\rmdefault}{\mddefault}{\updefault}$1$}}}
\put(3600,1848){\makebox(0,0)[lb]{{\SetFigFont{5}{6.0}{\rmdefault}{\mddefault}{\updefault}$1$}}}
\put(3485,1593){\makebox(0,0)[lb]{{\SetFigFont{5}{6.0}{\rmdefault}{\mddefault}{\updefault}$0$}}}
\put(3750,331){\makebox(0,0)[lb]{{\SetFigFont{8}{9.6}{\rmdefault}{\mddefault}{\updefault}$1$}}}
\put(4050,331){\makebox(0,0)[lb]{{\SetFigFont{8}{9.6}{\rmdefault}{\mddefault}{\updefault}$0$}}}
\put(4350,331){\makebox(0,0)[lb]{{\SetFigFont{8}{9.6}{\rmdefault}{\mddefault}{\updefault}$0$}}}
\put(4725,331){\makebox(0,0)[lb]{{\SetFigFont{8}{9.6}{\rmdefault}{\mddefault}{\updefault}$=\al$}}}
\put(2625,31){\makebox(0,0)[lb]{{\SetFigFont{8}{9.6}{\rmdefault}{\mddefault}{\updefault}$1$}}}
\put(2925,31){\makebox(0,0)[lb]{{\SetFigFont{8}{9.6}{\rmdefault}{\mddefault}{\updefault}$0$}}}
\put(3225,31){\makebox(0,0)[lb]{{\SetFigFont{8}{9.6}{\rmdefault}{\mddefault}{\updefault}$1$}}}
\put(2250,31){\makebox(0,0)[lb]{{\SetFigFont{8}{9.6}{\rmdefault}{\mddefault}{\updefault}$\be=$}}}
\put(792,1352){\makebox(0,0)[lb]{{\SetFigFont{8}{9.6}{\rmdefault}{\mddefault}{\updefault}$0$}}}
\put(1027,1494){\makebox(0,0)[lb]{{\SetFigFont{8}{9.6}{\rmdefault}{\mddefault}{\updefault}$1$}}}
\put(950,1644){\makebox(0,0)[lb]{{\SetFigFont{8}{9.6}{\rmdefault}{\mddefault}{\updefault}$1$}}}
\put(1107,1951){\makebox(0,0)[lb]{{\SetFigFont{8}{9.6}{\rmdefault}{\mddefault}{\updefault}$0$}}}
\put(1255,1951){\makebox(0,0)[lb]{{\SetFigFont{8}{9.6}{\rmdefault}{\mddefault}{\updefault}$0$}}}
\put(1176,1794){\makebox(0,0)[lb]{{\SetFigFont{8}{9.6}{\rmdefault}{\mddefault}{\updefault}$0$}}}
\put(1251,1659){\makebox(0,0)[lb]{{\SetFigFont{8}{9.6}{\rmdefault}{\mddefault}{\updefault}$0$}}}
\put(1391,1652){\makebox(0,0)[lb]{{\SetFigFont{8}{9.6}{\rmdefault}{\mddefault}{\updefault}$1$}}}
\put(1541,1359){\makebox(0,0)[lb]{{\SetFigFont{8}{9.6}{\rmdefault}{\mddefault}{\updefault}$0$}}}
\put(3192,1486){\makebox(0,0)[lb]{{\SetFigFont{8}{9.6}{\rmdefault}{\mddefault}{\updefault}$1$}}}
\put(3492,2092){\makebox(0,0)[lb]{{\SetFigFont{8}{9.6}{\rmdefault}{\mddefault}{\updefault}$1$}}}
\put(3655,2092){\makebox(0,0)[lb]{{\SetFigFont{8}{9.6}{\rmdefault}{\mddefault}{\updefault}$1$}}}
\put(3599,1959){\makebox(0,0)[lb]{{\SetFigFont{8}{9.6}{\rmdefault}{\mddefault}{\updefault}$1$}}}
\put(3501,1785){\makebox(0,0)[lb]{{\SetFigFont{8}{9.6}{\rmdefault}{\mddefault}{\updefault}$1$}}}
\put(3350,1784){\makebox(0,0)[lb]{{\SetFigFont{8}{9.6}{\rmdefault}{\mddefault}{\updefault}$0$}}}
\put(3783,1800){\makebox(0,0)[lb]{{\SetFigFont{8}{9.6}{\rmdefault}{\mddefault}{\updefault}$0$}}}
\put(3948,1492){\makebox(0,0)[lb]{{\SetFigFont{8}{9.6}{\rmdefault}{\mddefault}{\updefault}$1$}}}
\put(2111.413,208.717){\arc{2880.163}{3.0702}{4.0927}}
\blacken\path(641.082,224.952)(675.000,106.000)(701.050,226.919)(672.246,189.955)(641.082,224.952)
\put(3358,1495){\makebox(0,0)[lb]{{\SetFigFont{8}{9.6}{\rmdefault}{\mddefault}{\updefault}$1$}}}
\end{picture}
}

\end{center}
\caption{Reading $\al$ and $\be$
from partially filled $K_T$-puzzles.\lremind{alphabetaKT}}
\label{alphabetaKT}
\end{figure}
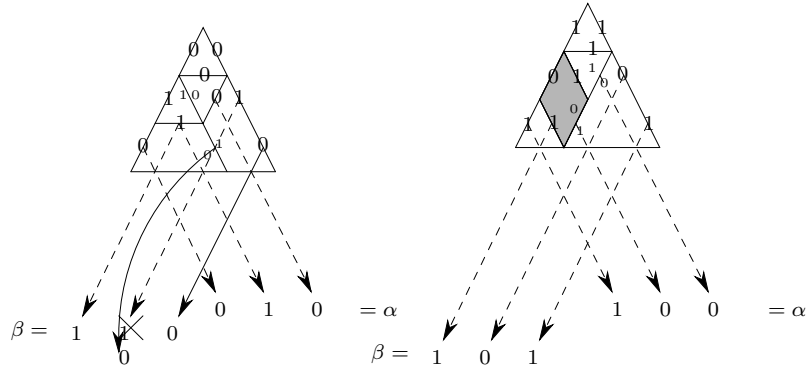

We now describe the conjectural geometry, by explaining how to
determine the equivariant puzzle variety corresponding to a partially
completed $K_T$-puzzle.
We first describe how to find $\al$ and $\be$ in words, and then
give examples which  may prove more enlightening.  We read off $\al$ and $\be$ from the edges as usual (see Figs.~ \ref{readalpha} and~\ref{readbeta}
respectively).  The $0$, $1$, $10$, and $01$ edges have the
same interpretation as before.    In addition:
\begin{itemize}
\item A gash dangling southwest, part of the first $K_T$-piece of Figure~\ref{KTpieces}, turns  the $0$ preceding it in $\al$ into a $1$.
(Because this piece must be placed immediately after an equivariant piece, as in Fig.~ \ref{KTplaceone}, the gash is {\em immediately} preceded by $0$
 in $\al$.)
\item  A gash dangling southeast, part of the second $K_T$-piece of Figure
~\ref{KTpieces}, turns the first $1$ following it in $\be$ into a $0$.
(There is guaranteed to be such a $1$ for each such gash by the placement rule of Fig.~\ref{KTplace}.)
\end{itemize}
Examples are given in Figure~\ref{alphabetaKT}.

\subsection{Extending Example~\ref{eqex} to equivariant $K$-theory}
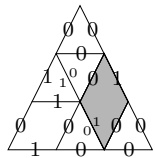
\begin{figure}[htbp]
\begin{center}
\setlength{\unitlength}{0.00083333in}
\begingroup\makeatletter\ifx\SetFigFont\undefined%
\gdef\SetFigFont#1#2#3#4#5{%
  \reset@font\fontsize{#1}{#2pt}%
  \fontfamily{#3}\fontseries{#4}\fontshape{#5}%
  \selectfont}%
\fi\endgroup%
{\renewcommand{\dashlinestretch}{30}
\begin{picture}(989,999)(0,-10)
\texture{88555555 55000000 555555 55000000 555555 55000000 555555 55000000 
	555555 55000000 555555 55000000 555555 55000000 555555 55000000 
	555555 55000000 555555 55000000 555555 55000000 555555 55000000 
	555555 55000000 555555 55000000 555555 55000000 555555 55000000 }
\shade\path(612,72)(462,372)(612,672)
	(762,372)(612,72)
\path(612,72)(462,372)(612,672)
	(762,372)(612,72)
\path(12,72)(462,972)(912,72)(12,72)
\path(162,372)(462,372)(312,72)
\path(612,72)(762,372)
\put(537,222){\makebox(0,0)[lb]{{\SetFigFont{5}{6.0}{\rmdefault}{\mddefault}{\updefault}$1$}}}
\put(326,463){\makebox(0,0)[lb]{{\SetFigFont{5}{6.0}{\rmdefault}{\mddefault}{\updefault}$1$}}}
\put(478,163){\makebox(0,0)[lb]{{\SetFigFont{5}{6.0}{\rmdefault}{\mddefault}{\updefault}$0$}}}
\put(391,526){\makebox(0,0)[lb]{{\SetFigFont{5}{6.0}{\rmdefault}{\mddefault}{\updefault}$0$}}}
\put(354,781){\makebox(0,0)[lb]{{\SetFigFont{8}{9.6}{\rmdefault}{\mddefault}{\updefault}$0$}}}
\put(505,781){\makebox(0,0)[lb]{{\SetFigFont{8}{9.6}{\rmdefault}{\mddefault}{\updefault}$0$}}}
\put(435,631){\makebox(0,0)[lb]{{\SetFigFont{8}{9.6}{\rmdefault}{\mddefault}{\updefault}$0$}}}
\put(225,487){\makebox(0,0)[lb]{{\SetFigFont{8}{9.6}{\rmdefault}{\mddefault}{\updefault}$1$}}}
\put(289,330){\makebox(0,0)[lb]{{\SetFigFont{8}{9.6}{\rmdefault}{\mddefault}{\updefault}$1$}}}
\put(56,181){\makebox(0,0)[lb]{{\SetFigFont{8}{9.6}{\rmdefault}{\mddefault}{\updefault}$0$}}}
\put(148,35){\makebox(0,0)[lb]{{\SetFigFont{8}{9.6}{\rmdefault}{\mddefault}{\updefault}$1$}}}
\put(356,179){\makebox(0,0)[lb]{{\SetFigFont{8}{9.6}{\rmdefault}{\mddefault}{\updefault}$0$}}}
\put(434,31){\makebox(0,0)[lb]{{\SetFigFont{8}{9.6}{\rmdefault}{\mddefault}{\updefault}$0$}}}
\put(733,31){\makebox(0,0)[lb]{{\SetFigFont{8}{9.6}{\rmdefault}{\mddefault}{\updefault}$0$}}}
\put(805,180){\makebox(0,0)[lb]{{\SetFigFont{8}{9.6}{\rmdefault}{\mddefault}{\updefault}$0$}}}
\put(654,181){\makebox(0,0)[lb]{{\SetFigFont{8}{9.6}{\rmdefault}{\mddefault}{\updefault}$0$}}}
\put(508,476){\makebox(0,0)[lb]{{\SetFigFont{8}{9.6}{\rmdefault}{\mddefault}{\updefault}$0$}}}
\path(462,372)(612,672)(312,672)(612,72)
\put(659,478){\makebox(0,0)[lb]{{\SetFigFont{8}{9.6}{\rmdefault}{\mddefault}{\updefault}$1$}}}
\end{picture}
}

\end{center}
\caption{The new $K_T$-puzzle arising in computing
$\si_{010}^2$. (The puzzles already arising in equivariant
cohomology are shown in Figure~\ref{enex}.)\lremind{swishy}}
\label{swishy}
\end{figure}

As an example, we continue the discussion of \S \ref{eqex}.  At the
key step of that example, we considered a Cartier divisor, and
observed that it had two components, which gave two equivariant cycles
that we analyzed further.  In equivariant $K$-theory, we are expecting
a third term, appearing with sign $-1$, corresponding to the
scheme-theoretic intersection of these first two pieces.  This leads
to the $K_T$-puzzle shown in Figure~\ref{swishy}.  This
corresponds to the fact that\lremind{basicKT}
\begin{eqnarray}
\label{basicKT}
\si_{010}^2 &=& \si_{100} + (1- e^{y_2-y_3}) \si_{010} - (1 - e^{y_2-y_3}) \si_{100} \\
&=&  e^{y_2-y_3} \si_{100} + (1 - e^{y_2-y_3}) \si_{010},
\nonumber
\end{eqnarray}
extending 
\eqref{basicHT}.
The first term on the right side of \eqref{basicKT}
corresponds to the left puzzle of Figure~\ref{enex},
the second term corresponds to the right side of Figure~\ref{enex},
and the third term corresponds to  the new $K_T$-puzzle of Figure~\ref{swishy}.  The partially completed puzzle corresponding to this puzzle is
the first panel of Figure~\ref{alphabetaKT}.

We now describe the geometry corresponding to this example.  Consider
the five-dimensional variety whose points generically correspond to the
configuration corresponding to the left side of Figure~\ref{enexg}.
As in the equivariant cohomology discussion, we consider the Cartier
divisor corresponding to requiring the moving line $\proj M_2$ to pass
through the fixed point $\proj F_1$.  This is reducible, and the
generic behaviors of the two irreducible components are shown on the
right side of Figure~\ref{enexg}.  These two components contributed to
the equivariant cohomology calculation.  We have a third term,
appearing with sign $-1$, corresponding to the scheme-theoretic
intersection of these two components.  This corresponds to those
configurations where the moving line $\proj M_2$ equals the fixed line
$\proj F_2$ (the condition of the lower-right panel of Fig.\
~\ref{enexg}) and also the point $p$ agrees with the fixed point
$\proj F_1$ (the condition of the upper-right panel of Fig.\
~\ref{enexg}).  This is indeed the equivariant puzzle variety
predicted by the conjecture.


The reader wishing to see the other $K_T$-piece in use
should compute (using the conjecture)
\begin{eqnarray*}
\si_{101}^2 &=& \si_{110} + (1 - e^{y_1-y_2}) \si_{101} - (1 - e^{y_1-y_2}) \si_{110} \\
&=& e^{y_1-y_2} \si_{110} + (1 - e^{y_1 - y_2}) \si_{101}.
\end{eqnarray*}
There are three puzzles, and the reader is encouraged to follow
the geometry of the successive Cartier slicing corresponding
to these puzzles.  The second panel of Figure~\ref{alphabetaKT}
will appear, contributing a $K_T$-term (not appearing in the $H_T$ calculation).

\section{A conjectural geometric Littlewood-Richardson rule
for the two-step flag variety}

\label{s:twostep}\lremind{s:twostep}We next use these same ideas to give a conjectural geometric
Littlewood-Richardson rule for two-step flag varieties (by the Knutson
and the second author), generalizing the geometric
Littlewood-Richardson rule for the Grassmannian.  Recent interest in
the two-step flag variety is likely due to the realization by A. Buch,
A. Kresch and H. Tamvakis \cite{bkt} that (i) the two-step problem is
intimately related to the quantum Grassmannian problem, and (ii) the
two-step problem appears to be simpler than the $m$-step
problem in general.

We first give Knutson's Littlewood-Richardson puzzle conjecture.

\begin{figure}[htbp]
\begin{center}
\setlength{\unitlength}{0.00083333in}
\begingroup\makeatletter\ifx\SetFigFont\undefined%
\gdef\SetFigFont#1#2#3#4#5{%
  \reset@font\fontsize{#1}{#2pt}%
  \fontfamily{#3}\fontseries{#4}\fontshape{#5}%
  \selectfont}%
\fi\endgroup%
{\renewcommand{\dashlinestretch}{30}
\begin{picture}(3851,434)(0,-10)
\put(3006,32){\makebox(0,0)[lb]{{\SetFigFont{8}{9.6}{\rmdefault}{\mddefault}{\updefault}$2(10)$}}}
\path(912,107)(1062,407)(1212,107)(912,107)
\path(1512,107)(1662,407)(1812,107)(1512,107)
\path(1962,107)(2112,407)(2262,107)(1962,107)
\path(2412,107)(2562,407)(2712,107)(2412,107)
\path(3012,107)(3162,407)(3312,107)(3012,107)
\path(3462,107)(3612,407)(3762,107)(3462,107)
\path(462,107)(612,407)(762,107)(462,107)
\put(63,223){\makebox(0,0)[lb]{{\SetFigFont{8}{9.6}{\rmdefault}{\mddefault}{\updefault}$0$}}}
\put(203,219){\makebox(0,0)[lb]{{\SetFigFont{8}{9.6}{\rmdefault}{\mddefault}{\updefault}$0$}}}
\put(129,59){\makebox(0,0)[lb]{{\SetFigFont{8}{9.6}{\rmdefault}{\mddefault}{\updefault}$0$}}}
\put(596,64){\makebox(0,0)[lb]{{\SetFigFont{8}{9.6}{\rmdefault}{\mddefault}{\updefault}$1$}}}
\put(524,223){\makebox(0,0)[lb]{{\SetFigFont{8}{9.6}{\rmdefault}{\mddefault}{\updefault}$1$}}}
\put(677,224){\makebox(0,0)[lb]{{\SetFigFont{8}{9.6}{\rmdefault}{\mddefault}{\updefault}$1$}}}
\put(963,223){\makebox(0,0)[lb]{{\SetFigFont{8}{9.6}{\rmdefault}{\mddefault}{\updefault}$2$}}}
\put(1113,223){\makebox(0,0)[lb]{{\SetFigFont{8}{9.6}{\rmdefault}{\mddefault}{\updefault}$2$}}}
\put(1034,64){\makebox(0,0)[lb]{{\SetFigFont{8}{9.6}{\rmdefault}{\mddefault}{\updefault}$2$}}}
\put(1563,223){\makebox(0,0)[lb]{{\SetFigFont{8}{9.6}{\rmdefault}{\mddefault}{\updefault}$2$}}}
\put(1712,226){\makebox(0,0)[lb]{{\SetFigFont{8}{9.6}{\rmdefault}{\mddefault}{\updefault}$0$}}}
\put(1610,65){\makebox(0,0)[lb]{{\SetFigFont{8}{9.6}{\rmdefault}{\mddefault}{\updefault}$20$}}}
\put(2024,217){\makebox(0,0)[lb]{{\SetFigFont{8}{9.6}{\rmdefault}{\mddefault}{\updefault}$1$}}}
\put(2161,225){\makebox(0,0)[lb]{{\SetFigFont{8}{9.6}{\rmdefault}{\mddefault}{\updefault}$0$}}}
\put(2071,59){\makebox(0,0)[lb]{{\SetFigFont{8}{9.6}{\rmdefault}{\mddefault}{\updefault}$10$}}}
\put(2458,217){\makebox(0,0)[lb]{{\SetFigFont{8}{9.6}{\rmdefault}{\mddefault}{\updefault}$2$}}}
\put(2612,223){\makebox(0,0)[lb]{{\SetFigFont{8}{9.6}{\rmdefault}{\mddefault}{\updefault}$1$}}}
\put(2510,59){\makebox(0,0)[lb]{{\SetFigFont{8}{9.6}{\rmdefault}{\mddefault}{\updefault}$21$}}}
\put(3057,229){\makebox(0,0)[lb]{{\SetFigFont{8}{9.6}{\rmdefault}{\mddefault}{\updefault}$2$}}}
\put(3171,228){\makebox(0,0)[lb]{{\SetFigFont{8}{9.6}{\rmdefault}{\mddefault}{\updefault}$10$}}}
\put(3472,217){\makebox(0,0)[lb]{{\SetFigFont{8}{9.6}{\rmdefault}{\mddefault}{\updefault}$21$}}}
\put(3664,223){\makebox(0,0)[lb]{{\SetFigFont{8}{9.6}{\rmdefault}{\mddefault}{\updefault}$0$}}}
\put(3452,31){\makebox(0,0)[lb]{{\SetFigFont{8}{9.6}{\rmdefault}{\mddefault}{\updefault}$(21)0$}}}
\path(12,107)(162,407)(312,107)(12,107)
\end{picture}
}

\end{center}
\caption{Knutson's $2$-step puzzle pieces, which may be rotated.  \lremind{twostep}}
\label{twostep}
\end{figure}
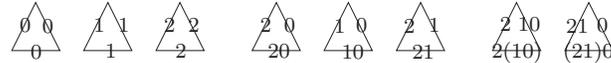

\begin{conjecture}[Knutson]  If
$\al$, $\be$, and $\ga$ are given in string notation,
the two-step Littlewood-Richardson number  $c^{\ga}_{\al, \be}$
is the number of puzzles with sides $\al$, $\be$, $\ga$ (written as in the
left panel of Fig.~\ref{puzzleD}), filled with puzzle pieces as given in 
Figure~\ref{twostep}.  \label{allentwostep}\lremind{allentwostep}
\end{conjecture}

This rule has been checked up to $n=16$ by Buch, Kresch, and
Tamvakis \cite{bkt}.  Given the number of cases to check, this
verification clearly required a great deal of ingenuity.
Conjecture~\ref{allentwostep} was originally stated as a conjecture by
Knutson for all partial flag varieties, but languished unpublished
once Knutson noted that this general version fails at $n=5$.
Buch, Kresch and Tamvakis noted that this generalization already
fails for three-step partial flag manifolds when $n=5$; Buch's patch
to the three-step conjecture is given in the next section.

The pieces are most cleanly described as follows.  The edges
correspond to binary trees with the nodes labeled by the integers $0$,
$1$, and $2$, such that the labels decrease strictly from left to right.
These trees can be represented by sequences of integers and
parentheses, as shown in Figure~\ref{twostep}.  The puzzle pieces 
consist of triangles labeled $x/x/x$ where $x \in \{ 0, 1, 2 \}$ as
well as triangles labeled $a/b/ab$ where $a$, $b$ and $ab$ are
acceptable edge-labels (read clockwise).  This description will be partially
conjecturally extended to the three-step case in the next section.

We now conjecture a geometric interpretation of partially filled
puzzles of this sort.  (This conjecture has been verified in a large
number of cases.) To any such partially filled puzzle, we associate a
{\em two-step puzzle variety}, extending the definition \eqref{puzzlevariety}
for Grassmannians,
of the form
$$
\overline{ \Omega_{\al}(M_\bullet) \cap \Omega_{\be}(F_\bullet)}.
$$
As before, we describe how to read $\al$ and $\be$ from the partially
filled puzzles, but now $\al$ and $\be$ will be strings of $0$'s, $1$'s,
and $2$'s.  It will be most convenient to describe the recipe using
the alternate visualization of Figure~\ref{newab}.  As in that figure,
on each horizontal edge labeled $10$, we temporarily glue on the
$1/0/10$ piece.  Similarly, for each horizontal edge labeled $20$,
$21$, $2(10)$, and $(21)0$, we temporarily glue on the $2/0/20$,
$2/1/21$, $2/10/2(10)$, and $21/0/(21)0$ piece (respectively).  Now
when attempting to read off $\al$ and $\be$, each edge visible will be
a $0$, $1$, or $2$ as desired, except for the following possibilities.
\begin{itemize}
\item As in the Grassmannian case, if the leading edge is labeled $10$, 
this counts for $0$ in $\al$, and in $\be$ turns the next $1$ to a $0$.
More generally, if a southwest/northeast edge is labeled
$ab$ ($a,b \in \{ 0,1,2 \}$, $a>b$), this counts for $b$ in $\al$, and in $\be$
turns the next $a$ to a $b$.
\item If the leading edge is $2(10)$, then this counts as $1$  in $\al$,
and in $\be$, the next $2$ is turned into a $1$, and the next $1$
(possibly the one just changed from a $2$) is turned into a $0$.
\item If the leading edge is $(21)0$, then this counts as $0$
in $\al$, and in $\be$ the next $2$ is turned into a $0$.  (The reader
may check that there is no $1$ before the $2$, so this could be
interpreted as the same statement  in the previous item:
the next $2$ is turned into a $1$, and the next $1$
is turned into a $0$.)
\item If there is a northwest-southeast edge labeled $10$
(this arises when temporarily gluing a $2/10/2(10)$-piece onto a horizontal
$2(10)$-edge), this counts as a $1$ in $\be$, and in $\al$, turns the next {\em earlier} $0$ into a $1$.
\end{itemize}

\section{Buch's conjectural combinatorial (non-geometric) rules in the
  three-step case, and for the two-step case in equivariant
  cohomology}

\label{s:anders}\lremind{s:anders}Buch has given combinatorial Littlewood-Richardson conjectures
in two additional cases:  three-step partial flag varieties in
ordinary cohomology, and two-step partial flag varieties in 
equivariant cohomology.  It is natural
given our earlier discussion to seek to understand the corresponding
conjectural geometry.  This may shed light on possible proofs.

We begin with Buch's three-step conjecture.  As usual, each
Littlewood-Richardson coefficient will count the number of puzzles
with sides corresponding to given strings (this time of $0$'s, $1$'s,
$2$'s, and $3$'s).  The triangular pieces are as follows.  Most of the
pieces will have edges that are analogous to the two-step case: they
correspond to binary trees with nodes labeled by $0$, $1$, $2$, or
$3$, with the labels strictly decreasing from left to the right.  We
write such trees as sequences of integers and parentheses.  The
complete list is $0$, $1$, $2$, $3$, $10$, $20$, $30$, $21$, $31$,
$32$, $(21)0$, $(31)0$, $(32)0$, $(32)1$, $2(10)$, $3(10)$, $3(20)$, $3(21)$, 
$((32)1)0$, $3((21)0)$, $(3(21))0$, $3(2(10))$, $(32)(10)$.  There are pieces
of the form $x/x/x$ where $x \in \{ 0, 1, 2, 3 \}$, and $a/b/ab$ where
$a$, $b$, and $ab$ are in the list above.  However, four extra pieces
are also required, shown in Figure~\ref{andersD}.  
Equivalently, an integer can be repeated in a tree, if the two copies are separated by exactly three parentheses.  
Buch has verified this
rule up to $n=9$. 

\begin{figure}[htbp]
\begin{center}
\setlength{\unitlength}{0.00083333in}
\begingroup\makeatletter\ifx\SetFigFont\undefined%
\gdef\SetFigFont#1#2#3#4#5{%
  \reset@font\fontsize{#1}{#2pt}%
  \fontfamily{#3}\fontseries{#4}\fontshape{#5}%
  \selectfont}%
\fi\endgroup%
{\renewcommand{\dashlinestretch}{30}
\begin{picture}(3353,406)(0,-10)
\put(2631,237){\makebox(0,0)[lb]{{\SetFigFont{8}{9.6}{\rmdefault}{\mddefault}{\updefault}$3$}}}
\path(1858,79)(2008,379)(2158,79)(1858,79)
\path(2608,79)(2758,379)(2908,79)(2608,79)
\path(358,79)(508,379)(658,79)(358,79)
\put(0,196){\makebox(0,0)[lb]{{\SetFigFont{8}{9.6}{\rmdefault}{\mddefault}{\updefault}$3(2(10))$}}}
\put(2032,195){\makebox(0,0)[lb]{{\SetFigFont{8}{9.6}{\rmdefault}{\mddefault}{\updefault}$(21)0$}}}
\put(958,31){\makebox(0,0)[lb]{{\SetFigFont{8}{9.6}{\rmdefault}{\mddefault}{\updefault}$(3(21))(10)$}}}
\put(1715,31){\makebox(0,0)[lb]{{\SetFigFont{8}{9.6}{\rmdefault}{\mddefault}{\updefault}$(32)((21)0)$}}}
\put(2495,31){\makebox(0,0)[lb]{{\SetFigFont{8}{9.6}{\rmdefault}{\mddefault}{\updefault}$3(((32)1)0)$}}}
\put(139,42){\makebox(0,0)[lb]{{\SetFigFont{8}{9.6}{\rmdefault}{\mddefault}{\updefault}$(3(2(10)))0$}}}
\put(912,180){\makebox(0,0)[lb]{{\SetFigFont{8}{9.6}{\rmdefault}{\mddefault}{\updefault}$3(21)$}}}
\put(2800,201){\makebox(0,0)[lb]{{\SetFigFont{8}{9.6}{\rmdefault}{\mddefault}{\updefault}$((32)1)0$}}}
\put(1786,231){\makebox(0,0)[lb]{{\SetFigFont{8}{9.6}{\rmdefault}{\mddefault}{\updefault}$32$}}}
\put(1296,224){\makebox(0,0)[lb]{{\SetFigFont{8}{9.6}{\rmdefault}{\mddefault}{\updefault}$10$}}}
\put(557,225){\makebox(0,0)[lb]{{\SetFigFont{8}{9.6}{\rmdefault}{\mddefault}{\updefault}$0$}}}
\path(1108,79)(1258,379)(1408,79)(1108,79)
\end{picture}
}

\end{center}
\caption{Buch's four unusual pieces in the $3$-step case. \lremind{andersD}}
\label{andersD}
\end{figure}

\begin{figure}[htbp]
\begin{center}
\setlength{\unitlength}{0.00083333in}
\begingroup\makeatletter\ifx\SetFigFont\undefined%
\gdef\SetFigFont#1#2#3#4#5{%
  \reset@font\fontsize{#1}{#2pt}%
  \fontfamily{#3}\fontseries{#4}\fontshape{#5}%
  \selectfont}%
\fi\endgroup%
{\renewcommand{\dashlinestretch}{30}
\begin{picture}(3443,639)(0,-10)
\texture{88555555 55000000 555555 55000000 555555 55000000 555555 55000000 
	555555 55000000 555555 55000000 555555 55000000 555555 55000000 
	555555 55000000 555555 55000000 555555 55000000 555555 55000000 
	555555 55000000 555555 55000000 555555 55000000 555555 55000000 }
\shade\path(2412,312)(2562,612)(2712,312)
	(2562,12)(2412,312)
\path(2412,312)(2562,612)(2712,312)
	(2562,12)(2412,312)
\shade\path(3012,312)(3162,612)(3312,312)
	(3162,12)(3012,312)
\path(3012,312)(3162,612)(3312,312)
	(3162,12)(3012,312)
\shade\path(612,312)(762,612)(912,312)
	(762,12)(612,312)
\path(612,312)(762,612)(912,312)
	(762,12)(612,312)
\shade\path(1212,312)(1362,612)(1512,312)
	(1362,12)(1212,312)
\path(1212,312)(1362,612)(1512,312)
	(1362,12)(1212,312)
\shade\path(1812,312)(1962,612)(2112,312)
	(1962,12)(1812,312)
\path(1812,312)(1962,612)(2112,312)
	(1962,12)(1812,312)
\shade\path(12,312)(162,612)(312,312)
	(162,12)(12,312)
\path(12,312)(162,612)(312,312)
	(162,12)(12,312)
\put(1409,428){\makebox(0,0)[lb]{{\SetFigFont{8}{9.6}{\rmdefault}{\mddefault}{\updefault}$2$}}}
\put(1256,126){\makebox(0,0)[lb]{{\SetFigFont{8}{9.6}{\rmdefault}{\mddefault}{\updefault}$2$}}}
\put(1443,131){\makebox(0,0)[lb]{{\SetFigFont{8}{9.6}{\rmdefault}{\mddefault}{\updefault}$1$}}}
\put(1826,428){\makebox(0,0)[lb]{{\SetFigFont{8}{9.6}{\rmdefault}{\mddefault}{\updefault}$10$}}}
\put(2014,423){\makebox(0,0)[lb]{{\SetFigFont{8}{9.6}{\rmdefault}{\mddefault}{\updefault}$2$}}}
\put(1856,126){\makebox(0,0)[lb]{{\SetFigFont{8}{9.6}{\rmdefault}{\mddefault}{\updefault}$2$}}}
\put(1998,132){\makebox(0,0)[lb]{{\SetFigFont{8}{9.6}{\rmdefault}{\mddefault}{\updefault}$10$}}}
\put(2451,428){\makebox(0,0)[lb]{{\SetFigFont{8}{9.6}{\rmdefault}{\mddefault}{\updefault}$0$}}}
\put(2592,429){\makebox(0,0)[lb]{{\SetFigFont{8}{9.6}{\rmdefault}{\mddefault}{\updefault}$21$}}}
\put(2434,120){\makebox(0,0)[lb]{{\SetFigFont{8}{9.6}{\rmdefault}{\mddefault}{\updefault}$21$}}}
\put(2609,126){\makebox(0,0)[lb]{{\SetFigFont{8}{9.6}{\rmdefault}{\mddefault}{\updefault}$0$}}}
\put(3026,439){\makebox(0,0)[lb]{{\SetFigFont{8}{9.6}{\rmdefault}{\mddefault}{\updefault}$10$}}}
\put(3193,435){\makebox(0,0)[lb]{{\SetFigFont{8}{9.6}{\rmdefault}{\mddefault}{\updefault}$21$}}}
\put(3032,126){\makebox(0,0)[lb]{{\SetFigFont{8}{9.6}{\rmdefault}{\mddefault}{\updefault}$21$}}}
\put(3197,126){\makebox(0,0)[lb]{{\SetFigFont{8}{9.6}{\rmdefault}{\mddefault}{\updefault}$10$}}}
\put(652,120){\makebox(0,0)[lb]{{\SetFigFont{8}{9.6}{\rmdefault}{\mddefault}{\updefault}$2$}}}
\put(226,428){\makebox(0,0)[lb]{{\SetFigFont{8}{9.6}{\rmdefault}{\mddefault}{\updefault}$1$}}}
\put(73,126){\makebox(0,0)[lb]{{\SetFigFont{8}{9.6}{\rmdefault}{\mddefault}{\updefault}$1$}}}
\put(215,126){\makebox(0,0)[lb]{{\SetFigFont{8}{9.6}{\rmdefault}{\mddefault}{\updefault}$0$}}}
\put(51,428){\makebox(0,0)[lb]{{\SetFigFont{8}{9.6}{\rmdefault}{\mddefault}{\updefault}$0$}}}
\put(651,428){\makebox(0,0)[lb]{{\SetFigFont{8}{9.6}{\rmdefault}{\mddefault}{\updefault}$0$}}}
\put(803,121){\makebox(0,0)[lb]{{\SetFigFont{8}{9.6}{\rmdefault}{\mddefault}{\updefault}$0$}}}
\put(809,428){\makebox(0,0)[lb]{{\SetFigFont{8}{9.6}{\rmdefault}{\mddefault}{\updefault}$2$}}}
\put(1273,429){\makebox(0,0)[lb]{{\SetFigFont{8}{9.6}{\rmdefault}{\mddefault}{\updefault}$1$}}}
\end{picture}
}

\end{center}
\caption{Buch's equivariant $2$-step pieces. 
They may not be rotated or reflected.  \lremind{andersB}}
\label{andersB}
\end{figure}
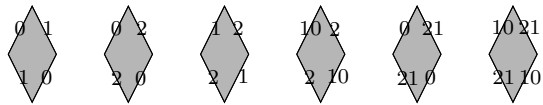

We next give Buch's equivariant $2$-step conjecture.  There are puzzle
pieces that are the same as the ordinary $2$-step conjecture
(Fig.~\ref{twostep}).  There are also six equivariant pieces, shown in
Figure~\ref{andersB}, generalizing Knutson and Tao's equivariant piece
for the Grassmannian (Fig.~\ref{eqpuzzle}).  Like the equivariant
Grassmannian piece, they may not be rotated or reflected.  Each
equivariant piece contributes a weight according to the same recipe as
the Grassmannian case (Fig.~\ref{eqweight}).  Buch has verified this
rule up to $n=7$. 

\section{A less explicit conjectural geometric Littlewood-Richardson
  rule for partial flag varieties in general}

  Despite the failure of Knutson's Conjecture~\ref{allentwostep} to
extend to all partial flag varieties, we are still led to a geometric
conjecture (Conj.~\ref{bigflag} below).  This is equivalent to
Conjecture~4.9 of \cite{vakil:checkers}.  It has been verified up to
$n=5$, and generalizes the geometric Littlewood-Richardson rule for
the Grassmannian, and the two-step geometric conjecture of
Section~\ref{s:twostep}.

Fix a variety $\FF l(a_1, \dots, a_m)$ of $m$-part partial flags in
$n$-space.  Let $S$ be the set of Schubert cells of $\FF l(a_1, \dots,
a_m)$ (i.e.\ certain strings of $0$'s, $1$'s, \dots, $m-1$'s), and let
$\{ d_0, \dots, d_{\binom n 2} \}$ be subset of the Schubert cells of
$\FF l(n)$ corresponding to the degeneration order (described in \S
\ref{degenerationorder}, and used throughout Part~\ref{partone}).  For
any element of $(\al, \be, d_i) \in S \times S \times \{ d_0, \dots,
d_{\binom n 2} \}$, fix two flags $M_{\bullet}$ and $F_{\bullet}$ in
relative position $d_i$ (there is a unique such pair up to
translation), and as before (e.g.\ \eqref{puzzlevariety}) define the {\em puzzle variety} $PV_{\al,
  \be, d_i}$
$$
\overline{ \Omega_{\al}(M_{\bullet}) \cap \Omega_{\be}(F_{\bullet})}.
$$
Note that  
$PV_{\al, \al, d_{\binom n 2}}$ is the Schubert variety $\oO_{\al}(F_{\bullet})$
and $PV_{\al, \be, d_0}$ 
is the intersection of general translates of the two Schubert varieties.
(Recall that $d_{\binom n 2}$ corresponds to $M_{\bullet} = F_{\bullet}$,
and $d_0$ corresponds to traverse $M_\bullet$ and $F_\bullet$.)

\begin{conjecture}
{\bf (\cite[Conj.~4.9]{vakil:checkers}, rephrased)}
There exists a subset $U$ of $S \times S \times \{d_0, \dots, d_{\binom n 2} \} = \{ (\al, \be, d_i) \}$
such that \label{bigflag}\lremind{bigflag}
\begin{itemize}
\item $(\al, \be, d_0) \in U$ for all $\al$ and $\be$.
\item If 
$(\al, \be, d_{\binom n 2}) \in U$, then $\al=\be$.
\item If $(\al, \be, d_i) \in U$ and $i< \binom n 2$ (i.e.\ $d_i$ is not the final
element of the degeneration order),
then upon degenerating $M_{\bullet}$ from relative position $d_i$
to $d_{i+1}$ with respect to $F_\bullet$, the puzzle variety
$PV_{\al, \be, d_i}$ degenerates to a union of other
puzzle varieties $PV_{\al', \be', d_{i+1}}$, where each $(\al', \be', d_{i+1})$
is in $U$, and 
each appears with multiplicity $1$.
\end{itemize}
\end{conjecture}

If this conjecture were true, and one had an explicit description of
which $\al'$ and $\be'$ arose at each step, one would have a
combinatorial Littlewood-Richardson with a geometric interpretation
(and motivation).  However, lacking even a conjectural explicit
description, this conjecture is admittedly vague and speculative.  The
main motivation for stating such a vague rule is that it suggests
where to look for a more precise rule, and also suggests how to
interpret more combinatorial conjectures.  It partly motivated the
conjectures of Knutson and the first author stated earlier.

One might also speculate that such a rule could also be extended to
equivariant $K$-theory, as in Section~\ref{KT}.

\part{MORE GENERAL RULES, MORE GENERAL DEGENERATION ORDERS}

\label{parttwo}
\section{The cohomology of flag varieties}\label{secflag}

In this section we will describe a different Littlewood-Richardson
rule for Grassmannians and explain how it generalizes to give a
Littlewood-Richardson rule for two-step flag varieties. These rules
will be in terms of combinatorial objects called Mondrian tableaux.
Mondrian tableaux are very efficient for encoding degenerations and
are very friendly to use. We encourage the reader to follow the
examples with graph paper and colored pencils in hand. 

\begin{figure}[htbp]
\begin{center}
\begin{picture}(0,0)%
\includegraphics{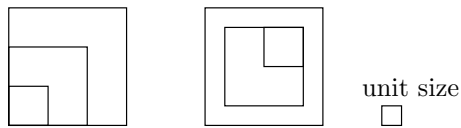}%
\end{picture}%
\setlength{\unitlength}{3947sp}%
\begingroup\makeatletter\ifx\SetFigFont\undefined%
\gdef\SetFigFont#1#2#3#4#5{%
  \reset@font\fontsize{#1}{#2pt}%
  \fontfamily{#3}\fontseries{#4}\fontshape{#5}%
  \selectfont}%
\fi\endgroup%
\begin{picture}(2775,763)(514,-523)
\put(2746,-326){\makebox(0,0)[lb]{\smash{{\SetFigFont{10}{12.0}{\rmdefault}{\mddefault}{\updefault}{\color[rgb]{0,0,0}unit size}%
}}}}
\end{picture}%

\end{center}
\caption{Two Mondrian tableaux associated to $\sigma_{2,1}$ in
  $G(3,6)$.}
\label{mondrian}
\end{figure}

A {\it Mondrian tableau} associated to a Schubert class
$\sigma_{\lambda_1, \cdots, \lambda_k}$ in $G(k,n)$ is a collection of
$k$ nested squares labeled by integers $1, \dots, k$, where the $j$th
square has size $n-k+j - \lambda_j$. The labels of the squares are
determined by the picture: assign $1$ to the smallest square; if a
square has label $i$, assign the label $i+1$ to the next larger
square.  Hence we will omit them when we depict Mondrian
tableaux. Figure \ref{mondrian} shows two examples of Mondrian
tableaux for $\sigma_{2,1}$ in $G(3,6)$.

In a Mondrian tableau a square of side-length $s$ denotes a vector
space of dimension $s$. We will denote squares in a Mondrian tableau
by capital letters in the math font (e.g., $A_i$) and the vector
spaces they represent by the corresponding letter in Roman font (e.g.,
$\mathrm{A}_i$). If a square $S_1$ is contained in a square $S_2$,
then the vector space $\mathrm{S}_1$ is a subspace of $\mathrm{S}_2$.
The reader should think of unit squares along the anti-diagonal of a
Mondrian tableau as a basis of the underlying vector space.  The
vector space represented by any square centered along the
anti-diagonal is the span of the basis elements it contains. A
Mondrian tableau associated to $\sigma_{\lambda}$ depicts the vector
spaces that have exceptional behavior for the $k$-planes parametrized
by the Schubert cycle. The $k$-planes are required to intersect the
vector space represented by the $i$th square in dimension at least
$i$.  

\begin{figure}[htbp]
\begin{center}
\epsfig{figure=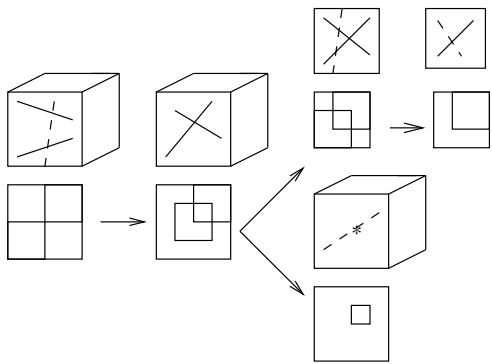}
\end{center}
\caption{The product $\sigma_1^2= \sigma_2 + \sigma_{1,1}$ in $G(2,4)=
  \GG(1,3)$: Mondrian tableaux and the projective geometry
  corresponding to them.}
\label{flipone}
\end{figure}

Before describing the rule in detail, we repeat the key example~\ref{akeyexample} of
Part \ref{partone} in terms of Mondrian tableaux. Figure \ref{flipone}
shows the calculation $\sigma_1^2 = \sigma_2 + \sigma_{1,1}$ in
$G(2,4) = \GG(1,3)$. The reader might want to refer to this example
while reading the rule.

The Mondrian tableaux rule follows the same basic strategy as the
rule of Part~\ref{partone}. We specialize the flags defining the Schubert varieties
until the intersection decomposes into Schubert varieties. However,
the Mondrian tableaux rule differs from the earlier rule in two
aspects.  First, the order of specialization is not pre-determined but
depends on the intersection problem. This added flexibility allows us
to avoid some geometric complications.  Second, the specializations do
not depend on a choice of basis. Although Schubert varieties are often
defined in terms of a fixed full-flag, one cannot canonically
associate a full-flag to a Schubert variety in $G(k,n)$. One can,
however, associate to it a canonical partial flag of at most
$k$-steps. The Mondrian tableaux rule will depend only on the
canonical partial flags.  

\begin{figure}[htbp]
\begin{center}
\begin{picture}(0,0)%
\includegraphics{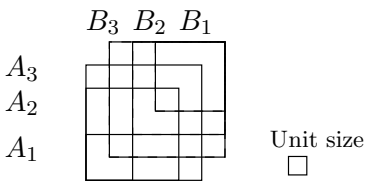}%
\end{picture}%
\setlength{\unitlength}{3947sp}%
\begingroup\makeatletter\ifx\SetFigFont\undefined%
\gdef\SetFigFont#1#2#3#4#5{%
  \reset@font\fontsize{#1}{#2pt}%
  \fontfamily{#3}\fontseries{#4}\fontshape{#5}%
  \selectfont}%
\fi\endgroup%
\begin{picture}(2248,1087)(376,-705)
\put(376,-547){\makebox(0,0)[lb]{\smash{{\SetFigFont{11}{13.2}{\rmdefault}{\mddefault}{\updefault}{\color[rgb]{0,0,0}$A_1$}%
}}}}
\put(376,-258){\makebox(0,0)[lb]{\smash{{\SetFigFont{11}{13.2}{\rmdefault}{\mddefault}{\updefault}{\color[rgb]{0,0,0}$A_2$}%
}}}}
\put(1468,250){\makebox(0,0)[lb]{\smash{{\SetFigFont{11}{13.2}{\rmdefault}{\mddefault}{\updefault}{\color[rgb]{0,0,0}$B_1$}%
}}}}
\put(885,250){\makebox(0,0)[lb]{\smash{{\SetFigFont{11}{13.2}{\rmdefault}{\mddefault}{\updefault}{\color[rgb]{0,0,0}$B_3$}%
}}}}
\put(1178,250){\makebox(0,0)[lb]{\smash{{\SetFigFont{11}{13.2}{\rmdefault}{\mddefault}{\updefault}{\color[rgb]{0,0,0}$B_2$}%
}}}}
\put(2040,-485){\makebox(0,0)[lb]{\smash{{\SetFigFont{9}{10.8}{\rmdefault}{\mddefault}{\updefault}{\color[rgb]{0,0,0}Unit size}%
}}}}
\put(376,-41){\makebox(0,0)[lb]{\smash{{\SetFigFont{11}{13.2}{\rmdefault}{\mddefault}{\updefault}{\color[rgb]{0,0,0}$A_3$}%
}}}}
\end{picture}%

\end{center}
\caption{The initial Mondrian tableau.}
\label{initial}
\end{figure}

\noindent {\bf The game.} To multiply two Schubert classes
$\sigma_{\lambda}$ and $\sigma_{\mu}$ in $G(k,n)$ we place the tableau
associated to $\lambda$ (respectively, $\mu$) at the southwest
(respectively, northeast) corner of an $n \times n$ square. The
squares in the $\lambda$ (respectively, $\mu$) tableau are all left
(respectively, right) aligned with respect to the $n \times n$ square.
We will denote the squares corresponding to $\lambda$ and $\mu$ by
$A_{i}$ and $B_{j}$, respectively.  Figure \ref{initial} shows the
initial tableau for the multiplication $\sigma_{2,1,1} \cdot
\sigma_{1,1,1}$ in $G(3,6)$.  

Initially the two Schubert cycles are defined with respect to two
transverse flags. If the intersection of the two Schubert cycles is
non-empty, then the Schubert cycles have to satisfy certain
conditions. A preliminary rule (MM rule) checks that these conditions
are satisfied. Then there are some simplifications that reduce the
problem to a smaller problem. The OS and S rules give these
simplifications.

 
$\bullet$ {\bf The MM (``must meet'') rule.}  We check that $A_i$
intersects $B_{k-i+1}$ in a square of side-length at least one for
every $i$ between $1$ and $k$.  If not, we stop:  the Schubert cycles
have empty intersection. The class of their intersection is zero.

 In a $k$-dimensional vector space $V^k$ every $i$-dimensional
subspace (such as $V^k \cap \mathrm{A}_i$) {\bf M}ust {\bf M}eet every
$(k-i+1)$-dimensional subspace (such as $V^k \cap \mathrm{B}_{k-i+1}$)
in at least a line. The intersection of two Schubert cycles is zero if
and only if the initial tableau formed by the two cycles does not
satisfy the MM rule.  

$\bullet$ {\bf The OS (``outer square'') rule.} We call the intersection
of $A_k$ and $B_k$ the {\bf O}uter {\bf S}quare of the tableau. We
replace every square with its intersection with the outer square. 

\begin{figure}[htbp]
\begin{center}
\begin{picture}(0,0)%
\includegraphics{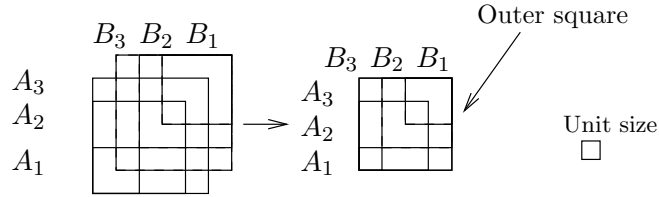}%
\end{picture}%
\setlength{\unitlength}{3947sp}%
\begingroup\makeatletter\ifx\SetFigFont\undefined%
\gdef\SetFigFont#1#2#3#4#5{%
  \reset@font\fontsize{#1}{#2pt}%
  \fontfamily{#3}\fontseries{#4}\fontshape{#5}%
  \selectfont}%
\fi\endgroup%
\begin{picture}(4048,1207)(376,-705)
\put(376,-547){\makebox(0,0)[lb]{\smash{{\SetFigFont{11}{13.2}{\rmdefault}{\mddefault}{\updefault}{\color[rgb]{0,0,0}$A_1$}%
}}}}
\put(376,-258){\makebox(0,0)[lb]{\smash{{\SetFigFont{11}{13.2}{\rmdefault}{\mddefault}{\updefault}{\color[rgb]{0,0,0}$A_2$}%
}}}}
\put(1468,250){\makebox(0,0)[lb]{\smash{{\SetFigFont{11}{13.2}{\rmdefault}{\mddefault}{\updefault}{\color[rgb]{0,0,0}$B_1$}%
}}}}
\put(885,250){\makebox(0,0)[lb]{\smash{{\SetFigFont{11}{13.2}{\rmdefault}{\mddefault}{\updefault}{\color[rgb]{0,0,0}$B_3$}%
}}}}
\put(1178,250){\makebox(0,0)[lb]{\smash{{\SetFigFont{11}{13.2}{\rmdefault}{\mddefault}{\updefault}{\color[rgb]{0,0,0}$B_2$}%
}}}}
\put(2191,-550){\makebox(0,0)[lb]{\smash{{\SetFigFont{11}{13.2}{\rmdefault}{\mddefault}{\updefault}{\color[rgb]{0,0,0}$A_1$}%
}}}}
\put(2191,-114){\makebox(0,0)[lb]{\smash{{\SetFigFont{11}{13.2}{\rmdefault}{\mddefault}{\updefault}{\color[rgb]{0,0,0}$A_3$}%
}}}}
\put(2917,104){\makebox(0,0)[lb]{\smash{{\SetFigFont{11}{13.2}{\rmdefault}{\mddefault}{\updefault}{\color[rgb]{0,0,0}$B_1$}%
}}}}
\put(2336,104){\makebox(0,0)[lb]{\smash{{\SetFigFont{11}{13.2}{\rmdefault}{\mddefault}{\updefault}{\color[rgb]{0,0,0}$B_3$}%
}}}}
\put(2627,104){\makebox(0,0)[lb]{\smash{{\SetFigFont{11}{13.2}{\rmdefault}{\mddefault}{\updefault}{\color[rgb]{0,0,0}$B_2$}%
}}}}
\put(2191,-333){\makebox(0,0)[lb]{\smash{{\SetFigFont{11}{13.2}{\rmdefault}{\mddefault}{\updefault}{\color[rgb]{0,0,0}$A_2$}%
}}}}
\put(3840,-305){\makebox(0,0)[lb]{\smash{{\SetFigFont{9}{10.8}{\rmdefault}{\mddefault}{\updefault}{\color[rgb]{0,0,0}Unit size}%
}}}}
\put(3301,378){\makebox(0,0)[lb]{\smash{{\SetFigFont{10}{12.0}{\rmdefault}{\mddefault}{\updefault}{\color[rgb]{0,0,0}Outer square}%
}}}}
\put(376,-41){\makebox(0,0)[lb]{\smash{{\SetFigFont{11}{13.2}{\rmdefault}{\mddefault}{\updefault}{\color[rgb]{0,0,0}$A_3$}%
}}}}
\end{picture}%

\end{center}
\caption{An application of the OS rule.}
\label{initial2}
\end{figure}

Since the $k$-planes are contained in both $\mathrm{A}_k$ and
$\mathrm{B}_k$, they must be contained in their intersection.  Figure
\ref{initial2} shows an application of the OS rule for the
intersection in Figure \ref{initial}.   

\begin{figure}[htbp]
\begin{center}
\begin{picture}(0,0)%
\includegraphics{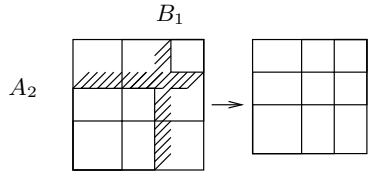}%
\end{picture}%
\setlength{\unitlength}{3947sp}%
\begingroup\makeatletter\ifx\SetFigFont\undefined%
\gdef\SetFigFont#1#2#3#4#5{%
  \reset@font\fontsize{#1}{#2pt}%
  \fontfamily{#3}\fontseries{#4}\fontshape{#5}%
  \selectfont}%
\fi\endgroup%
\begin{picture}(2262,1068)(676,-1124)
\put(676,-638){\makebox(0,0)[lb]{\smash{\SetFigFont{9}{10.8}{\rmdefault}{\mddefault}{\updefault}{\color[rgb]{0,0,0}$A_2$}%
}}}
\put(1599,-176){\makebox(0,0)[lb]{\smash{\SetFigFont{9}{10.8}{\rmdefault}{\mddefault}{\updefault}{\color[rgb]{0,0,0}$B_1$}%
}}}
\end{picture}

\end{center}
\caption{Adjusting the span of the linear constraints.}
\label{rule1}
\end{figure}

$\bullet$ {\bf The S (``span'') rule.} We check that $A_i$ and
$B_{k-i}$ touch or have a common square. If not, we remove the rows
and columns between these squares as shown in Figure \ref{rule1}.

This rule corresponds to the fact that a $k$-dimensional vector space
lies in the {\bf S}pan of any two of its subspaces of complementary
dimension whose only intersection is the origin. This rule removes any
basis element of the ambient vector space that is not needed in
expressing the $k$-planes parametrized by the intersection of the two
Schubert varieties.  

\noindent Once we have performed these preliminary steps, we will
inductively build a new flag (the $D$ flag) by degenerating the two
flags (the $A$ and $B$ flags).  At each stage of the game we will have
a partially built new flag (depicted by $D$ squares that arise as
intersections of $A$ and $B$ squares) and partially remaining $A$ and
$B$ flags (depicted by squares $A_i, \dots, A_k$ and $B_k$, $B_{k-i},
\dots, B_1$). After nesting the $D$ squares, we will increase the
dimension of the intersection of $\mathrm{A}_i$ with
$\mathrm{B}_{k-i}$ by one in order of increasing $i$. We will depict
this move in the Mondrian tableau by sliding $A_i$ anti-diagonally up
by one unit. As we specialize the flags, the intersection of the Schubert
varieties will break into irreducible components. Admissible Mondrian
tableaux depict the intermediate varieties that occur during the
process.  

A Mondrian tableau is {\it admissible} for $G(k,n)$ if the squares
 that constitute the tableau (except for the outer square) are uniquely
 labeled as an indexed $A$, $B$ or $D$ square such that

\begin{enumerate}
\item \label{1} The squares $A_k = B_k$ form the outer square.  They have
  side-length $m \leq n$ and contain the entire tableau.

\item \label{2} The $A$ squares are nested, distinct, left aligned and
  strictly contain all the $D$ squares. If the number of $D$ squares
  is $i-1 < k$, then the $A$ squares are $A_i, A_{i+1}, \dots, A_k$
  with the smaller index corresponding to the smaller squares.  (In
  particular, the total number of $A$ and $D$ squares is $k$.)

\item \label{3} The $B$ squares are nested, distinct and right
  aligned. They are labeled $B_k, B_{k-i}, B_{k-i-1}, \dots, B_1$,
  where a smaller square has the smaller index. The $A$ and $B$
  squares satisfy the MM and S rules. The $D$ squares may intersect
  $B_{k-i}$, but none are contained in $B_{k-i}$.

\item \label{4} The $D$ squares are distinct and labeled $D_1, \dots,
  D_{i-1}$.  They do not need to be nested; however, there can be at
  most one unnested $D$ square.  An {\it unnested} $D$ square is a $D$
  square that does not contain every $D$ square of smaller index. More
  precisely, if $D_j$ does not contain all the $D$ squares of smaller
  index, then it does not contain any of the $D$ squares of smaller
  index; it is contained in every $D$ square of larger index; and $D_h
  \subset D_s$ for every $h < s$ as long as $h$ and $s$ are different
  from $j$. All the $D$ squares of index lower than $j$ are to the
  southwest of $D_j$.  $D_{j-1}$ and $D_j$ share a common square or
  corner. Figure \ref{Dboxes} shows a typical configuration of $D$
  squares.

\begin{figure}[htbp]
\begin{center}
\begin{picture}(0,0)%
\includegraphics{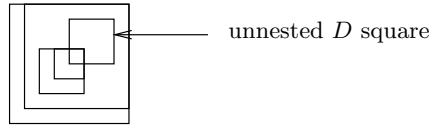}%
\end{picture}%
\setlength{\unitlength}{3947sp}%
\begingroup\makeatletter\ifx\SetFigFont\undefined%
\gdef\SetFigFont#1#2#3#4#5{%
  \reset@font\fontsize{#1}{#2pt}%
  \fontfamily{#3}\fontseries{#4}\fontshape{#5}%
  \selectfont}%
\fi\endgroup%
\begin{picture}(2361,774)(1714,-1498)
\put(3099,-951){\makebox(0,0)[lb]{\smash{{\SetFigFont{9}{10.8}{\rmdefault}{\mddefault}{\updefault}{\color[rgb]{0,0,0}unnested $D$ square}%
}}}}
\end{picture}%

\end{center}
\caption{A typical configuration of $D$ squares.}
\label{Dboxes}
\end{figure}    

\item \label{5} Let $S_1$ and $S_2$ be any two squares of the
    tableau. If the number of squares contained in their span but not
    contained in $S_1$ is $r$, then the side-length of $S_1$ is at
    least $r$ less than the side-length of their span.
\end{enumerate} 

\begin{figure}[htbp]
\begin{center}
\begin{picture}(0,0)%
\includegraphics{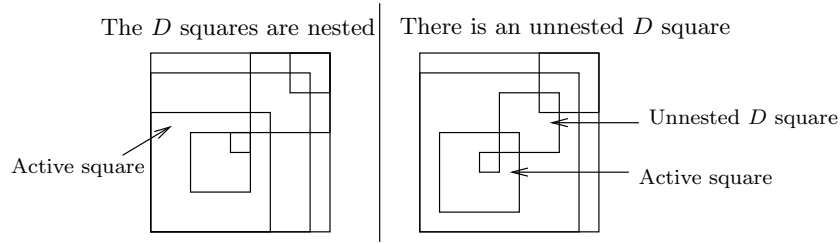}%
\end{picture}%
\setlength{\unitlength}{3947sp}%
\begingroup\makeatletter\ifx\SetFigFont\undefined%
\gdef\SetFigFont#1#2#3#4#5{%
  \reset@font\fontsize{#1}{#2pt}%
  \fontfamily{#3}\fontseries{#4}\fontshape{#5}%
  \selectfont}%
\fi\endgroup%
\begin{picture}(5108,1524)(-74,-748)
\put(2367,571){\makebox(0,0)[lb]{\smash{{\SetFigFont{9}{10.8}{\rmdefault}{\mddefault}{\updefault}{\color[rgb]{0,0,0}There is an unnested $D$ square}%
}}}}
\put(3868,-369){\makebox(0,0)[lb]{\smash{{\SetFigFont{8}{9.6}{\rmdefault}{\mddefault}{\updefault}{\color[rgb]{0,0,0}Active square}%
}}}}
\put(3931,  5){\makebox(0,0)[lb]{\smash{{\SetFigFont{8}{9.6}{\rmdefault}{\mddefault}{\updefault}{\color[rgb]{0,0,0}Unnested $D$ square}%
}}}}
\put(489,567){\makebox(0,0)[lb]{\smash{{\SetFigFont{9}{10.8}{\rmdefault}{\mddefault}{\updefault}{\color[rgb]{0,0,0}The $D$ squares are nested}%
}}}}
\put(-74,-306){\makebox(0,0)[lb]{\smash{{\SetFigFont{8}{9.6}{\rmdefault}{\mddefault}{\updefault}{\color[rgb]{0,0,0}Active square }%
}}}}
\end{picture}%

\end{center}
\caption{Two typical admissible Mondrian tableaux.}
\label{typmon}
\end{figure} 

\noindent Figure \ref{typmon} depicts two typical admissible Mondrian
tableaux.  Observe that the labels in an admissible Mondrian tableau
are determined from the picture.  The $B$ squares are the $k-i+1$
squares that are aligned with the north and east edges of the tableau.
The number of $A$ squares is equal to the number of $B$ squares, so
the $A$ squares are the largest $k-i+1$ squares that are aligned with
the south and west edges of the tableau. The rest of the squares are
$D$ squares.  Consequently, when we depict Mondrian tableaux we will
omit the labels.

To every admissible Mondrian tableau $M$, we can associate an
irreducible subvariety $\mathrm{M}$ of $G(k,n)$.  Define $AB_h$ to be
the intersection of $A_h$ and $B_{k-h+1}$ for $h=i+1, \dots, k$.
First, suppose $M$ does not contain any unnested $D$ squares.
$\mathrm{M}$ is defined as the closure of the locus of $k$-planes $V$
that satisfy
\begin{enumerate}
\item  $\dim(V \cap \mathrm{D}_s) = s$, for $1 \leq s
\leq i-1$, 

\item  $\dim(V \cap \mathrm{A}_i) = i$,
 
\item$\dim(V \cap \mathrm{AB}_h) = 1$ for $h=i+1, \dots,
k$,

\item $V$ is spanned by its intersection with $\mathrm{A}_i$ and
  $\mathrm{AB}_h$ for $h=i+1, \dots, k$.
\end{enumerate}
If $M$ has an unnested square $D_j$, we  modify the above by
requiring $\dim(V \cap \mathrm{D}_j) =1$. 

\begin{figure}[htbp]
\begin{center}
\begin{picture}(0,0)%
\includegraphics{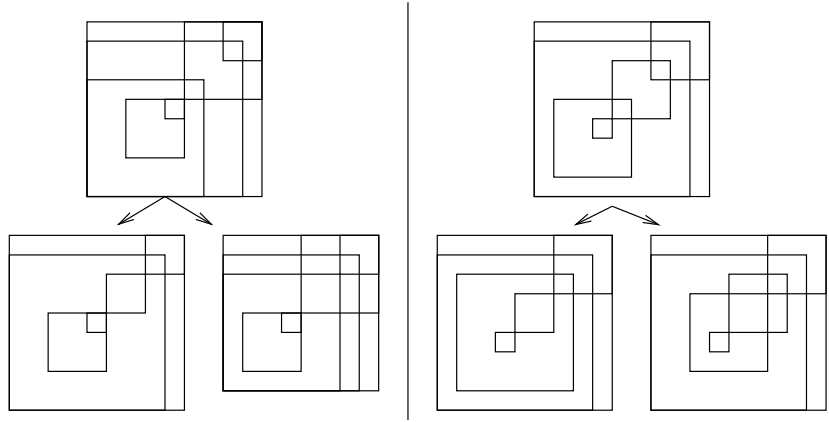}%
\end{picture}%
\setlength{\unitlength}{3947sp}%
\begingroup\makeatletter\ifx\SetFigFont\undefined%
\gdef\SetFigFont#1#2#3#4#5{%
  \reset@font\fontsize{#1}{#2pt}%
  \fontfamily{#3}\fontseries{#4}\fontshape{#5}%
  \selectfont}%
\fi\endgroup%
\begin{picture}(5154,2649)(409,-1948)
\end{picture}%

\end{center}
\caption{Simplifying the Mondrian tableaux in Figure \ref{typmon}.}
\label{move}
\end{figure}

\noindent {\bf The algorithm:} We now describe an algorithm that
simplifies admissible Mondrian tableaux.  Let $M$ be an admissible
Mondrian tableau with an outer square of side-length $m$. If $M$
contains no $A$ or $B$ squares (other than the outer square) and all
the $D$ squares are nested, then $M$ is a Mondrian tableau associated
to a Schubert variety. The algorithm terminates for $M$. Otherwise,
$M$ can be simplified as follows: If all the $D$ squares in $M$ are
nested, define the {\it active square} to be the smallest $A$ square
$A_i$. If $D_j$ is the unique unnested $D$ square in $M$, define the
{\it active square} to be $D_{j-1}$.  Move the active square
anti-diagonally up by one unit.  If there are any $D$ squares aligned
with the south and west edges of the active square, move them
anti-diagonally up by one unit with the active square. Keep all the
other squares fixed. Replace $M$ by the following two tableaux unless
the second tableau is not admissible (see Figure \ref{move}). In the
latter case replace $M$ with only Tableau 1.

$\bullet$ {\bf Tableau 1.}  If the active square is $A_i$, delete
$A_i$ and $B_{k-i}$. Draw their new intersection and label it $D_i$.
Keep their old span as the outer square.  If $D_i$ does not intersect
or touch $B_{k-i-1}$, slide all the $D$ squares anti-diagonally up
until $D_{i}$ touches $B_{k-i-1}$.  If the active square is $D_{j-1}$,
delete $D_{j-1}$ and $D_j$. Draw their new intersection and label it
$D_{j-1}$. Draw their old span and label it $D_j$. If $D_{j-2}$ does
not intersect or touch the new $D_{j-1}$, slide all the $D$ squares of
index $j-2$ anti-diagonally until $D_{j-2}$ touches $D_{j-1}$.  All
the remaining squares stay as in $M$.  

$\bullet$ {\bf Tableau 2.}  If the active square is $A_i$, we shrink
the outer square by one unit so that it passes along the new boundary
of $A_i$ and $B_{k-i}$ and we delete the column and row that lies
outside this square. The rest of the squares stay as in $M$.  If the
active square is $D_{j-1}$, we place the squares we move in their new
positions and keep the rest of the squares as in $M$.

Observe that if the active square is $A_i$, Tableau 2 is not
admissible if either $A_{i+1}$ has side-length one larger than $A_i$
or if $B_{k-i}$ has side-length $m-i$ (informally, if $A_i$ or
$B_{k-i}$ are as large as possible given $A_{i+1}$ and $B_{k}$).  If
the active square is $D_{j-1}$, then Tableau 2 is not admissible
either if the side-length of $D_j$ is not at least $j-1$ units smaller
than the side-length of the span of $D_j$ and $D_{j-1}$ or if
$D_{j-1}$ contains $D_j$ as a result of the move (informally, if
$D_{j-1}$ and $D_j$ are as large as possible). In these cases we
replace $M$ with only Tableaux 1.  

Geometrically, moving the active square corresponds to a
specialization of the flags defining the Schubert varieties. Let us
describe this explicitly in the case the active square is $A_i$ and
there are no $D$ squares abutting the south and west edges of $A_i$.
Let $s$ be the side-length of $A_i$ and suppose that initially $A_{i}$
and $B_{k-i}$ intersect in a square of side-length $r$. There is a
family of $s$-dimensional linear spaces $\mathrm{A}_i(t)$
parametrized by an open subset $0 \in U \subset \PP^1$ such that over
the points $t \in U$ with $ t \not= 0$, the dimension of intersection
$\mathrm{A}_i (t) \cap \mathrm{B}_{k-i}$ is equal to $r$ and when
$t=0$, the dimension of intersection $\mathrm{A}_i(0) \cap
\mathrm{B}_{k-i}$ is $r+1$.  Denoting the basis vectors represented by
the unit squares along the diagonal by $e_1, \dots, e_n$, we
can take this family to be
$$\mathrm{A}_i (t) = \ \mbox{the span of} \ \{(te_1 + (1-t) e_{s+1},
e_2, \dots, e_s \}.$$ When $t=1$, we have our original vector space
$\mathrm{A}_i$ represented by the old position of the square $A_i$.
When $t=0$, we have the new vector space $\mathrm{A}_i(0)$ represented
by the new position of the square $A_i$. When $t=0$, the intersection
of Schubert varieties defined with respect to the $A$ and $B$ flags
either remains irreducible or breaks into two irreducible components.
The two tableaux record these possibilities. The Littlewood-Richardson
rule may be phrased informally as: {\it If the $k$-planes in the limit
  do not intersect $\mathrm{A}_i (0) \cap \mathrm{B}_{k-i}$, then they
  must be contained in their new span.} 

More precisely, after we apply the MM, OS and S rules, the initial
Mondrian tableau is admissible.  It is clear that the move transforms
an admissible Mondrian tableau to one or two new admissible Mondrian
tableaux.  Therefore, we can continue applying the move to each of the
resulting tableau. After a cycle of moves (beginning with nested $D$
squares; forming a new $D$ square; and nesting the $D$ squares again),
we decrease the number of $A$ and $B$ squares each by one and we
increase the number of nested $D$ squares by one. If we continue
applying the move to every tableau that results from the initial
Mondrian tableau, after finitely many steps all the resulting tableaux
are tableaux associated to Schubert varieties. The
Littlewood-Richardson coefficient $c_{\lambda, \mu}^{\nu}$ of $G(k,n)$
is equal to the number of times the tableau corresponding to
$\sigma_{\nu}$ occurs at the end of this algorithm.

\begin{theorem}[\cite{coskun1:LR}, Thm.~ 3.1] \label{grmain} The
  Littlewood-Richardson coefficient $c_{\lambda, \mu}^{\nu}$ of
  $G(k,n)$ equals the number of times the Mondrian tableau associated
  to $\sigma_{\nu}$ results in a game of Mondrian tableaux starting
  with $\sigma_{\lambda}$ and $\sigma_{\mu}$ in an $n \times n$
  square.
\end{theorem}

\begin{figure}[htbp]
\begin{center}
  \epsfig{figure=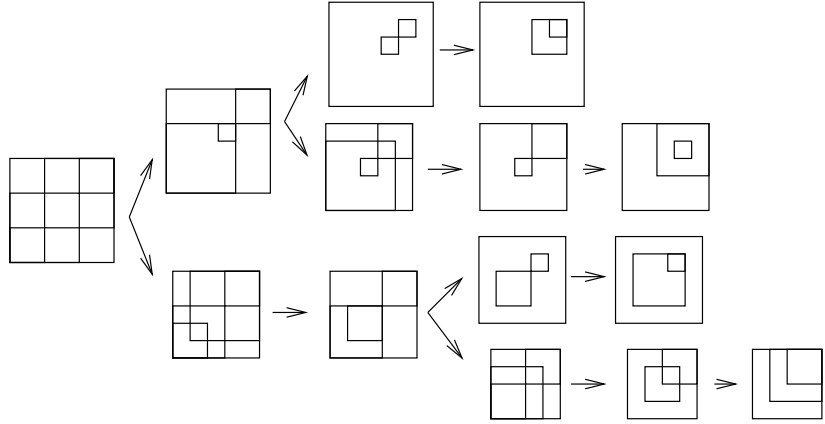}
\end{center}
\caption{The product $\sigma_{2,1}^2 = \sigma_{3,3} + 2
  \sigma_{3,2,1} + \sigma_{2,2,2}$ in $G(3,6)$.}
\label{grasex}
\end{figure}

In Figure \ref{grasex}, we compute $\sigma_{2,1} \cdot \sigma_{2,1}$
in $G(3,6)$ using the Mondrian tableaux rule.  When we move $A_1$,
there are two possibilities. We replace the initial tableau by the two
tableaux where we replace $A_1$ and $B_2$ by their new intersection
(and slide it up) and where we restrict the tableau to the new span of
$A_1$ and $B_2$. We continue resolving the first tableau by moving
$A_2$ and replacing it by two new tableaux. In the second tableau,
$B_2$ is as large as possible given the outer square. Hence, when we
move $A_1$ again, there is only one possibility. We then move $A_2$
and now there are two possibilities.  We replace the tableau with the
tableau where we take the intersection of $A_2$ and $B_1$ and with the
tableau where we restrict the tableau to the new span of $A_2$ and
$B_1$.  Continuing we conclude that
$$\sigma_{2,1}^2 = \sigma_{3,3} + 2 \sigma_{3,2,1} + \sigma_{2,2,2}.$$


More generally, we can define a {\it generalized Mondrian tableau} in
$G(k,n)$ to be any collection of $k$ squares centered along the
anti-diagonal of an $n \times n$ square satisfying the following two
properties:
\begin{enumerate}
\item None of the squares are equal to the span of the squares
  contained in them.

\item Let $S_1$ and $S_2$ be any two squares in the tableau. If the
  number of squares contained in their span but not contained in $S_1$
  is $r$, then the side-length of $S_1$ is at least $r$ less than the
  side-length of their span.
\end{enumerate}
We can associate an irreducible subvariety of the Grassmannian
$G(k,n)$ to such a tableau. Consider the locus of $k$-planes $V$
satisfying the following properties. For any square $S$ in the
tableau, $V$ intersects $\mathrm{S}$ in dimension equal to the number
of squares contained in $S$ (including itself).  If $S_1$ and $S_2$
are any two squares of the tableau, we further require $V \cap
\mathrm{S}_1$ and $V \cap \mathrm{S}_2$ to intersect only along the
subspaces represented by squares contained in both $S_1$ and $S_2$ and
otherwise to be independent.  The variety associated to the
generalized Mondrian tableau is the closure in $G(k,n)$ of the variety
parameterizing such $k$-planes. 

\cite{coskun1:LR} contains an algorithm for computing the classes of
varieties associated to generalized Mondrian tableaux. Observe that
the intersection of two Schubert varieties is a special case.
Replacing the $A$ and $B$ squares in the initial tableau with the
intersections $A_i \cap B_{k-i+1}$ for $i=1, \dots, k$ we obtain a
generalized Mondrian tableau. Here we will not discuss the rule that
expresses the classes of the varieties defined by generalized Mondrian
tableaux as a sum of Schubert varieties. However, we note that this
flexibility of Mondiran tableaux makes it possible to adapt them to
other contexts. 

\subsection{Painted Mondrian tableaux} \label{pmt}Mondrian tableaux are very
well-suited for computations in the cohomology of partial flag
varieties. For instance, a modification of the algorithm for the
Grassmannians yields a Littlewood-Richardson rule for two-step flag
varieties. In this subsection we will describe this rule and give
examples of how to express classes of subvarieties of $r$-step partial
flag varieties in terms of Schubert varieties.  
 
\begin{figure}[htbp]
\begin{center}
\begin{picture}(0,0)%
\includegraphics{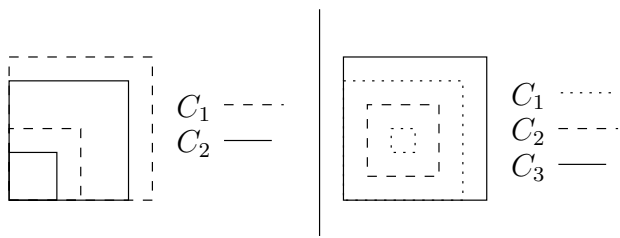}%
\end{picture}%
\setlength{\unitlength}{3947sp}%
\begingroup\makeatletter\ifx\SetFigFont\undefined%
\gdef\SetFigFont#1#2#3#4#5{%
  \reset@font\fontsize{#1}{#2pt}%
  \fontfamily{#3}\fontseries{#4}\fontshape{#5}%
  \selectfont}%
\fi\endgroup%
\begin{picture}(3849,1449)(889,-1648)
\put(4051,-811){\makebox(0,0)[lb]{\smash{{\SetFigFont{12}{14.4}{\rmdefault}{\mddefault}{\updefault}{\color[rgb]{0,0,0}$C_1$  }%
}}}}
\put(4051,-1036){\makebox(0,0)[lb]{\smash{{\SetFigFont{12}{14.4}{\rmdefault}{\mddefault}{\updefault}{\color[rgb]{0,0,0}$C_2$}%
}}}}
\put(4051,-1261){\makebox(0,0)[lb]{\smash{{\SetFigFont{12}{14.4}{\rmdefault}{\mddefault}{\updefault}{\color[rgb]{0,0,0}$C_3$}%
}}}}
\put(1951,-886){\makebox(0,0)[lb]{\smash{{\SetFigFont{12}{14.4}{\rmdefault}{\mddefault}{\updefault}{\color[rgb]{0,0,0}$C_1$  }%
}}}}
\put(1951,-1111){\makebox(0,0)[lb]{\smash{{\SetFigFont{12}{14.4}{\rmdefault}{\mddefault}{\updefault}{\color[rgb]{0,0,0}$C_2$}%
}}}}
\end{picture}%

\end{center}
\caption{Mondrian tableaux associated to $\sigma_{1,1,0,0}^{2,1,2,1}$
  in $Fl(2,4;6)$ and $\sigma_{2,1,0,0}^{1,2,1,3}$ in $Fl(2,3,4;6)$.}
\label{paintmon}
\end{figure}

In order to denote Schubert varieties in $\Flag$, we need to allow the
squares in a Mondrian tableau to have one of $r$ colors.  Let $C_1,
\dots, C_r$ be $r$ colors ordered by their indices.  A Mondrian
tableau associated to the Schubert cycle $\sigma_{\lambda_1, \dots,
  \lambda_{k_r}}^{\delta_1, \dots, \delta_{k_r}}$ is a collection of
$k_r$ nested squares in $r$ colors.  The $i$th square has side-length
$n-k_r + i - \lambda_i$ and color $C_{\delta_i}$. (In particular,
$k_1$ of the squares have color $C_1$; and $k_i - k_{i-1}$ of the
squares have color $C_i$ for $i > 1$.) Each square is labeled by an
$r$-vector of non-negative integers, where the $j$th index denotes the
number of squares of color less than or equal to $C_{j}$ in the square
(including possibly the square itself).  Figure \ref{paintmon} shows
two examples of painted Mondrian tableaux. When depicting Mondrian
tableaux for two-step flag varieties, we will always use dashed lines
for the color $C_1$ and solid lines for the color $C_2$.

As in the Grassmannian case, the labels of the squares are clear from
the picture. The smallest square has last ($r$th) index $1$. If a
square has last index $i$, the next larger square has last index
$i+1$. The square with last index $i$ has as its $j$th index the
number of squares of color at most $C_{j}$ contained in it.
Geometrically, the vector space $V_j$ in the tuple $(V_1, \dots, V_r)$
parametrized by the Schubert variety is required to intersect the
vector space represented by a square in dimension equal to the number
of squares of color at most $C_j$ in that square.  

We will now describe the Littlewood-Richardson rule for two-step flag
varieties. The strategy and many of the details are very similar to
the case of Grassmannians described above. Here we will focus mainly
on examples. The reader should refer to \cite{coskun1:LR} for more
details and to 
\verb+http://math.mit.edu/~coskun/gallery.html+
for more examples.  

\begin{figure}[htbp]
\begin{center}
\begin{picture}(0,0)%
\includegraphics{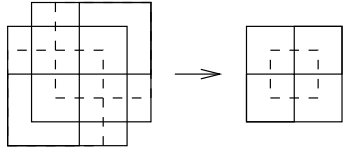}%
\end{picture}%
\setlength{\unitlength}{3947sp}%
\begingroup\makeatletter\ifx\SetFigFont\undefined%
\gdef\SetFigFont#1#2#3#4#5{%
  \reset@font\fontsize{#1}{#2pt}%
  \fontfamily{#3}\fontseries{#4}\fontshape{#5}%
  \selectfont}%
\fi\endgroup%
\begin{picture}(2124,924)(589,-823)
\end{picture}%

\end{center}
\caption{The initial tableau and the OS rule for the multiplication
  $\sigma_{1,1,1}^{2,1,2} \cdot \sigma_{1,1,1}^{2,1,2}$ in
  $Fl(1,3;6)$.}
\label{paintinit}
\end{figure}

As in the case of Grassmannians, in order to multiply two Schubert
cycles, we place the painted Mondrian tableaux associated to the two
Schubert cycles in opposite corners of an $n\times n$ square. We make
sure that the tableau at the southwest (respectively, northeast)
corner is left (respectively, right) justified. We refer to the
squares at the southwest corner as $A$ squares and the squares at the
northeast corner as $B$ squares. The initial position represents the
case when the flags defining the Schubert varieties are transverse.
See the left panel of Figure \ref{paintinit} for an example.

If the intersection of the two Schubert cycles is not zero, then the
defining flag elements have to satisfy certain intersection
conditions. The MM (must meet) rule ensures that these conditions are
satisfied.  First, the $i$th square in the southwest corner has to
intersect the $(k_2 - i +1)$th square in the northeast corner in a
square of side-length at least one for every $i$. Similarly, the $i$th
square of color $C_1$ in the southwest corner has to intersect the
$(k_1 - i +1)$th square of color $C_1$ in the northeast corner. The
intersection of two Schubert cycles is non-zero if and only if this
rule is satisfied. After checking this, we simplify the tableau
applying the OS (outer square) and S (span) rules. We apply the OS
rule for squares of both colors. We first restrict the tableau to the
intersection of the largest $A$ and $B$ squares. We then restrict all
the squares of color $C_1$ to the intersection of the largest $A$ and
$B$ squares of color $C_1$. The right panel in Figure \ref{paintinit}
shows an example. 

Once we have performed these preliminary operations, we move the $A$
squares in a specified order (with two exceptions the same order as
before). We thus eliminate the $\mathrm{A}$ and $\mathrm{B}$ flags and
build a new $\mathrm{D}$ flag represented by $D$ squares. At each
stage the painted Mondrian tableau corresponds to a subvariety of the
two-step flag variety.  Very generally, a painted Mondrian tableau is
a collection of $k_2$ squares (possibly not the span of consecutive
unit squares) of color $C_2$ that satisfy the two conditions for a
generalized Grassmannian tableau and $k_1$ squares of color $C_1$
consisting of the spans of the squares in color $C_2$ that satisfy the
conditions for a generalized Mondrian tableau for $G(k_1,k_2)$. The
tableaux that occur during our process can be viewed as a very special
subset of these tableaux called {\it admissible painted Mondrian
  tableaux}. We refer the reader to \cite{coskun1:LR} for their
precise description.

We can associate an irreducible subvariety of the two-step flag
variety $\Fl$ to a painted Mondrian tableau: Take the closure of the
locus of pairs $(V_1, V_2)$ such that 
\begin{enumerate}
\item $V_2$ intersects the vector space represented by a square of
  color $C_2$ in dimension equal to the number of squares of color
  $C_2$ contained in that square.

\item $V_1$ intersects the subspace of $V_2$ contained in a vector
  space represented by a square of color $C_1$ in dimension equal to
  the number of squares of color $C_1$ in that square.
\end{enumerate}
The dimension of such a variety has a simple expression in terms of
the sum of side-lengths and containment relations among the squares.
The dimension is the sum of two terms. The first term is the sum of
the side-lengths of the squares of color $C_2$ minus one for every
containment relation between them. The second term is the the number
of squares of color $C_2$ contained in every square of color $C_1$
minus the number of containment relations among the squares of color
$C_1$. These terms geometrically correspond to the dimensions of the
image and of the fiber, respectively, of the projection of the variety
from $Fl(k_1,k_2;n)$ to $G(k_2,n)$. 

\begin{figure}[htbp]
\begin{center}
\epsfig{figure=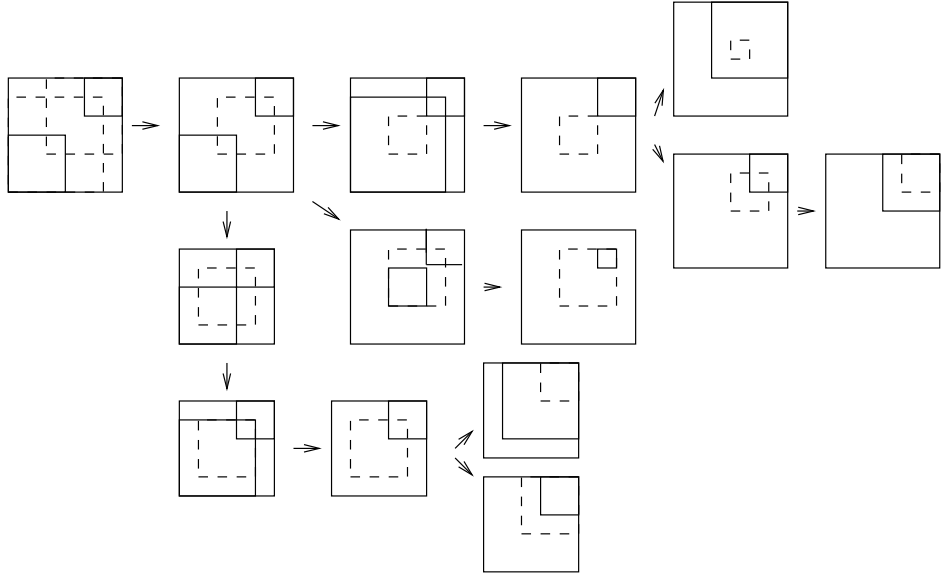}
\end{center}
\caption{Calculating the intersection  $\sigma_{1,0,0}^{2,1,2} \cdot
\sigma_{2,1,0}^{2,1,2}$ in $Fl(1,3;6)$.}
\label{twost1}
\end{figure}

We demonstrate the rule by calculating the intersection
$\sigma_{1,0,0}^{2,1,2} \cdot \sigma_{2,1,0}^{2,1,2}$ in $Fl(1,3;6)$
(see Figure \ref{twost1}).  We first move the smallest $A$ square (the
active square) anti-diagonally up by one unit. (Geometrically, this
move corresponds to specializing the third flag element of the
$\mathrm{A}$ flag $F_3^{\mathrm{A}}$.) After each move, we replace the
initial tableau by three new tableaux unless the dimension of the
variety associated to one or more of these tableaux is strictly
smaller than the initial dimension. If the dimension is smaller, we
discard that tableau. In Tableau I, we replace the active
square and the largest square of color $C_1$ by their new intersection
in color $C_1$. We complete the rest of the tableaux so that each
square continues to have the same number of squares and the same
number of squares of color $C_1$ as in the initial tableau. This is
depicted by the tableau to the right of the initial tableau.
(Geometrically, this possibility corresponds to the case where $V_1$
intersects the new intersection.)  In Tableau II, we draw the
new intersection of the active square and the largest $\mathrm{B}$
square with which the intersection increases in color $C_2$. We
complete the rest of the tableaux so that each square continues to
have the same number of squares and the same number of squares of
color $C_1$ as in the initial tableau.  This is depicted by the
tableau to the southeast of the initial tableau. (Geometrically, this
possibility corresponds to the case where $V_2$ intersects the new
intersection.)  In Tableau III, we restrict the outer square to
pass through the southwest corner of the active square. This is
depicted by the tableau directly below the initial tableau.
(Geometrically, neither $V_1$ nor $V_2$ intersect the new
intersection, so they are contained in the new span.)  

\begin{figure}[htbp]
\begin{center}
  \epsfig{figure=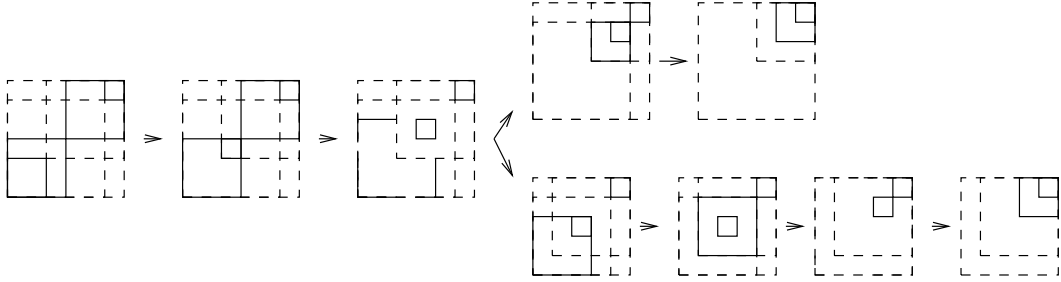}
\end{center}
\caption{Calculating the intersection $\sigma_{1,1,0,0}^{2,2,1,1} \cdot
  \sigma_{2,1,1,0}^{2,2,1,1}$ in $Fl(2,4;6)$. }
\label{twost2}
\end{figure}

The main new feature of the rule for two-step flag varieties is
squares called {\it fillers}. Fillers are squares of color $C_2$ that
occur as the intersection of the active square with a square of color
$C_1$ to its northeast. The newly formed intersection in the second
tableau above is a filler. Fillers affect the degeneration order: the
smallest filler takes precedence over the smallest $A$ square in the
order. Otherwise, the order follows the same pattern as in the
Grassmannian (see \cite{coskun1:LR}). Moving fillers may also cause
some of the squares containing it to become disconnected. We refer to
squares that are not the span of consecutive unit squares as {\it
  chopped squares}. See Figure \ref{twost2} below for an example.

We continue simplifying each of the three tableaux. Since the $D$
squares are nested, in the first and third tableau the active square
is the smallest $A$ square. In the second tableau the active square is
the filler $D$ square. Let us trace how the first tableau simplifies
(the top row of Figure \ref{twost1}). When we move $A_2$, only Tableau
II can occur. Now there is an unnested $D$ square ($D_2$). The active
square is $D_1$. When we move it, only Tableaux I and III occur. Note
that the possibility described by Tableau II would have dimension one
less than the initial dimension. Continuing the calculation one sees
that
$$\sigma_{1,0,0}^{2,1,2} \cdot
\sigma_{2,1,0}^{2,1,2} = \sigma_{3,1,0}^{1,2,2} +
\sigma_{2,2,0}^{1,2,2} + \sigma_{3,2,0}^{2,1,2} +
\sigma_{2,1,1}^{1,2,2} + \sigma_{2,2,1}^{2,1,2}.$$

To illustrate the new features of the two-step rule, we calculate
$\sigma_{1,1,0,0}^{2,2,1,1} \cdot \sigma_{2,1,1,0}^{2,2,1,1}$ in
$Fl(2,4;6)$ (see Figure \ref{twost2}).  Since the calculation is
similar to the previous case we only point out the new features. In
the second panel, the newly formed $D$ square is a filler. According
to the degeneration order, we first move the filler.  This move causes
the vector space corresponding to the smallest $A$ box to no longer be
the span of consecutive basis elements. This is depicted on the
Mondrian tableau as a chopped square. The third panel has an example
of such a square. Continuing the calculation one sees that  
$$\sigma_{1,1,0,0}^{2,2,1,1} \cdot
\sigma_{2,1,1,0}^{2,2,1,1}= \sigma_{1,1,2,2}^{2,2,1,1} +
\sigma_{2,2,2,0}^{2,2,1,1}. $$

In general, degenerating the flags in the order described decomposes
the intersection of two Schubert varieties in a two-step flag variety
to a union of Schubert varieties. This process may be encoded in an
explicit game of painted Mondrian tableaux. The main theorem then is
the following:

\begin{theorem}[\cite{coskun1:LR}, Thm.~4.1] \label{twostepbig}
Let $\sigma_{\lambda}$ and $\sigma_{\mu}$ denote two Schubert cycles
in the flag variety $Fl(k_1, k_2; n)$. Let their product be
$\sigma_{\lambda} \cdot \sigma_{\mu} = \sum_{\nu} c_{\lambda
\mu}^{\nu} \sigma_{\nu}$. The coefficient $c_{\lambda \mu}^{\nu}$ is
equal to the number of times the painted Mondrian tableau of
$\sigma_{\nu}$ occurs in a game of Mondrian tableau played by starting
with the Mondrian tableaux of $\sigma_{\lambda}$ and $\sigma_{\mu}$ in
an $n \times n$ square.
\end{theorem}

It is natural to wonder about the connection of this rule to the
geometric two-step puzzle conjecture of \S \ref{s:twostep}.  Because
in the theorem above the degenerations used are more general than in
Part~\ref{partone}, it is not clear if much direct connection is to be
expected.

\begin{figure}[htbp]
\begin{center}
\begin{picture}(0,0)%
\includegraphics{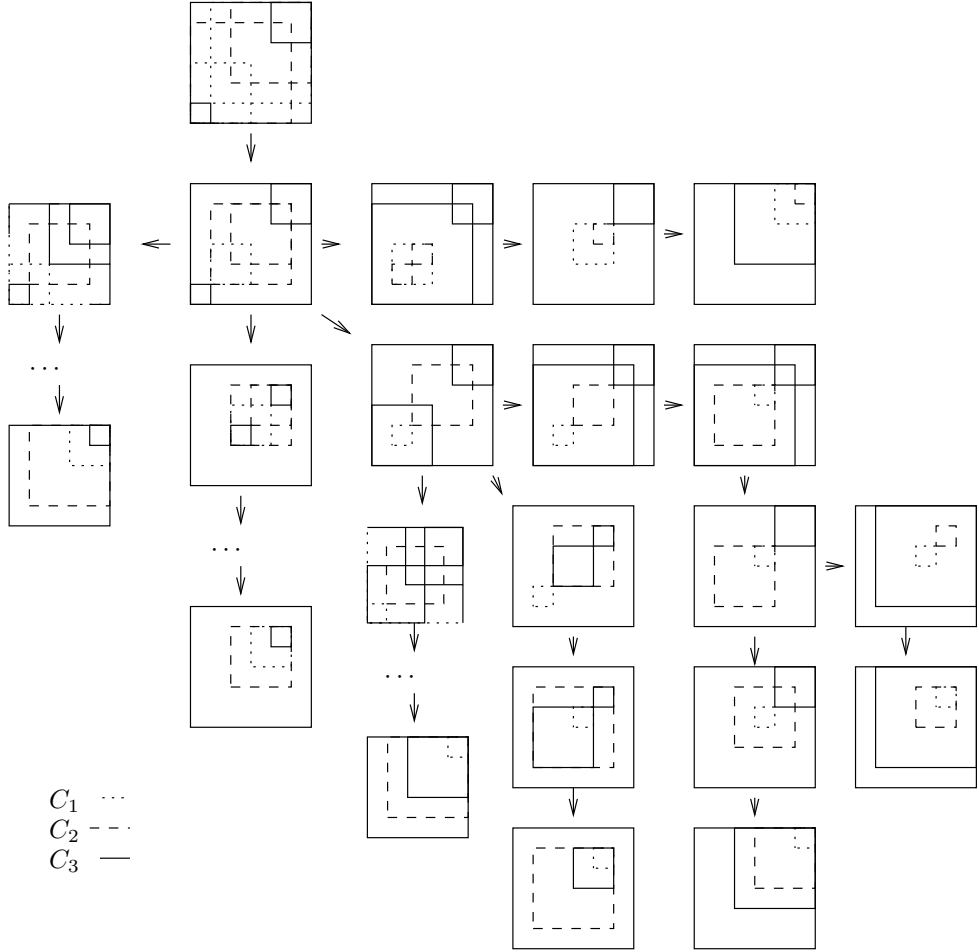}%
\end{picture}%
\setlength{\unitlength}{3947sp}%
\begingroup\makeatletter\ifx\SetFigFont\undefined%
\gdef\SetFigFont#1#2#3#4#5{%
  \reset@font\fontsize{#1}{#2pt}%
  \fontfamily{#3}\fontseries{#4}\fontshape{#5}%
  \selectfont}%
\fi\endgroup%
\begin{picture}(6099,5972)(1039,-5271)
\put(1304,-4373){\makebox(0,0)[lb]{\smash{{\SetFigFont{10}{12.0}{\rmdefault}{\mddefault}{\updefault}{\color[rgb]{0,0,0}$C_1$}%
}}}}
\put(1304,-4563){\makebox(0,0)[lb]{\smash{{\SetFigFont{10}{12.0}{\rmdefault}{\mddefault}{\updefault}{\color[rgb]{0,0,0}$C_2$}%
}}}}
\put(1304,-4753){\makebox(0,0)[lb]{\smash{{\SetFigFont{10}{12.0}{\rmdefault}{\mddefault}{\updefault}{\color[rgb]{0,0,0}$C_3$}%
}}}}
\put(1178,-1652){\makebox(0,0)[lb]{\smash{{\SetFigFont{10}{12.0}{\rmdefault}{\mddefault}{\updefault}{\color[rgb]{0,0,0}$\cdots$}%
}}}}
\put(2317,-2791){\makebox(0,0)[lb]{\smash{{\SetFigFont{10}{12.0}{\rmdefault}{\mddefault}{\updefault}{\color[rgb]{0,0,0}$\cdots$}%
}}}}
\put(3407,-3590){\makebox(0,0)[lb]{\smash{{\SetFigFont{10}{12.0}{\rmdefault}{\mddefault}{\updefault}{\color[rgb]{0,0,0}$\cdots$}%
}}}}
\end{picture}%

\end{center}
\caption{The product $\sigma_{2,1,0,0}^{3,1,2,3} \cdot
  \sigma_{1,0,0,0}^{3,2,1,3} $ in $Fl(1,2,4;6)$ equals
  $\sigma_{2,2,1,0}^{2,1,3,3} + \sigma_{2,2,0,0}^{1,2,3,3} +
  \sigma_{2,1,1,0}^{1,2,3,3} + \sigma_{2,2,1,0}^{1,3,2,3} +
  \sigma_{2,1,1,1}^{1,3,2,3} + \sigma_{2,2,2,0}^{3,1,2,3} +
  \sigma_{2,2,1,1}^{3,1,2,3}$.}
\label{3step}
\end{figure}

We conclude this section with an example of how Mondrian tableaux can
be used to calculate classes of subvarieties in other partial flag
varieties (see Figure \ref{3step} for an example in $Fl(1,2,4;6)$).
Informally, the procedure may
be described as follows. We move the active square anti-diagonally up
by one unit and replace the initial tableau with new tableaux. Each of
the new tableaux is obtained by either placing the active square in
its new position (and shrinking the outer square if necessary) or by
drawing the new intersection of the active square with the closest
square of color at most $C_j$ to its northeast in color $C_j$ and
filling the rest of the tableau so that each square has the same
number of squares of color at most $C_h$ as in the initial tableau for
every color $C_h$. We then discard the tableaux that represent
varieties with strictly smaller dimension.

\section{Quantum cohomology of Grassmannians and flag varieties}

\label{qcg}In this section we will describe a geometric method for computing the
small quantum cohomology of Grassmannians. The method is based on the
elegant observation of Buch, Kresch and Tamvakis that three-point
genus-zero Gromov-Witten invariants of certain homogeneous varieties
may be computed as the ordinary intersection products of convenient
auxiliary varieties.  In the case of Grassmannians these auxiliary
varieties are two-step flag varieties.  

Let $X$ be a smooth, complex projective variety. Let
$\beta \in H_2 (X, \ZZ)$ be the homology class of a curve. The
Kontsevich moduli space of $m$-pointed, genus-zero stable maps
$\Kgnb{0,m} (X, \beta)$ provides a useful compactification of the
space of rational curves on $X$ whose homology class is $\beta$.
Recall that $\Kgnb{0,m} (X, \beta)$ is the smooth, proper
Deligne-Mumford stack parameterizing the data of
\begin{enumerate}
\item[(i)] $C$, a proper, connected, at-worst-nodal curve of
  arithmetic genus $0$,
\item[(ii)]   $p_1,\dots,p_m$, an ordered sequence of distinct, smooth
  points of $C$,
\item[(iii)] and $f:C\rightarrow X$, a morphism with $f_*[C] = \beta$
  satisfying the following stability condition: every irreducible
  component of $C$ mapped to a point under $f$ contains at least $3$
  special points, i.e., marked points $p_i$ and nodes of $C$.
\end{enumerate}

The Kontsevich moduli space $\Kgnb{0,m}(X, \beta)$ admits $m$
evaluation morphisms to $X$, where the $i$th evaluation morphism
$ev_i$ maps $[C, p_1, \dots, p_m, f]$ to $p_i$.  The dimension of
$\Kgnb{0,m} (X, \beta)$ is $$c_1(T_X) \cdot \beta + \dim(X) + m - 3.$$
Given the classes of $m$ subvarieties $\gamma_1, \dots, \gamma_m$ of
$X$ whose codimensions add up to this dimension, the genus-zero
Gromov-Witten invariant of $X$ associated to the classes $\gamma_1,
\dots, \gamma_m$ and the curve class $\beta$ is defined as
$$I_{X, \beta} (\gamma_1, \dots, \gamma_m) : =
\int_{\Kgnb{0,m}(X,\beta)} ev_1^*(\gamma_1) \cup \cdots \cup
ev_m^*(\gamma_m).$$

When $X$ is a homogeneous variety, the Gromov-Witten invariant is
enumerative in the following sense:

\begin{proposition}
\cite[Lem.\ ~14]{fultonpan:kontsevich} \label{enumerative} Let
  $\Gamma_1, \dots, \Gamma_m$ be general  subvarieties of a
  homogeneous variety $X$ representing the Poincar\'e duals of the
  cohomology classes $\gamma_1, \dots, \gamma_m$, respectively.  The
  scheme theoretic intersection
  $$ev_1^{-1}(\Gamma_1) \cap  \cdots \cap ev_m^{-1}(\Gamma_m)$$ is a
  finite number of reduced points in $\mathcal{M}_{0,m}(X, \beta)$. Moreover,
  $$I_{X,\beta}(\gamma_1, \dots, \gamma_m)= \# \  ev_1^{-1}(\Gamma_1)
  \cap \cdots \cap ev_m^{-1}(\Gamma_m).$$
\end{proposition}

When $m=3$, the genus-zero Gromov-Witten invariants can be used as the
structure constants of a commutative and associative ring called the
small quantum cohomology ring.  The main theorem of \cite{bkt}
equates the three-point Gromov-Witten invariants of Grassmannians with
certain ordinary three-point intersections of two-step flag
varieties. Given a Schubert cycle $\sigma_{\lambda}$ in $G(k,n)$, there is a
special Schubert cycle $X_{\lambda}^{(d)}(F_{\bullet})$ in $Fl(k-d,
k+d; n)$ defined by
\begin{align*}
X_{\lambda}^{(d)}(F_{\bullet}) := \{ \ (V_1, V_2) \ | \  \dim (V_1
\cap F_{n-i-\lambda_{k-i}}) \geq k-d-i ,  \ \dim (V_2 \cap F_{n-k+j -
\lambda_j}) \geq j \ \}
\end{align*}
where $1 \leq i \leq k-d$ and $1 \leq j \leq k$. 

\begin{proposition}\label{two} \cite[Cor.\ ~1]{bkt} 
  Let $\lambda, \mu, \nu$ be partitions and $d \geq 0$ be an integer
  satisfying
\begin{equation}\label{eqn}
|\lambda| + |\mu| + |\nu| = k(n-k) + dn.
\end{equation}
Then the degree $d$ three-point Gromov-Witten invariants of $G(k,n)$
equal the ordinary three-point intersections of special Schubert
varieties  in the flag variety $Fl(k-d,k+d;n)$:
$$I_d (\sigma_{\lambda}, \sigma_{\mu}, \sigma_{\nu}) =
\int_{Fl(k-d,k+d;n)} [X_{\lambda}^{(d)}] \cup [X_{\mu}^{(d)}] \cup
[X_{\nu}^{(d)}].$$
\end{proposition}

\begin{figure}[htbp]
\begin{center}
\begin{picture}(0,0)%
\includegraphics{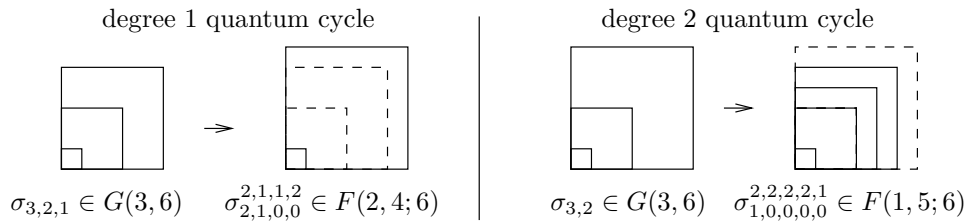}%
\end{picture}%
\setlength{\unitlength}{3947sp}%
\begingroup\makeatletter\ifx\SetFigFont\undefined%
\gdef\SetFigFont#1#2#3#4#5{%
  \reset@font\fontsize{#1}{#2pt}%
  \fontfamily{#3}\fontseries{#4}\fontshape{#5}%
  \selectfont}%
\fi\endgroup%
\begin{picture}(5708,1351)(526,-517)
\put(4238,711){\makebox(0,0)[lb]{\smash{{\SetFigFont{10}{12.0}{\rmdefault}{\mddefault}{\updefault}{\color[rgb]{0,0,0}degree 2 quantum cycle}%
}}}}
\put(526,-441){\makebox(0,0)[lb]{\smash{{\SetFigFont{10}{12.0}{\rmdefault}{\mddefault}{\updefault}{\color[rgb]{0,0,0}$\sigma_{3,2,1} \in G(3,6)$}%
}}}}
\put(1870,-441){\makebox(0,0)[lb]{\smash{{\SetFigFont{10}{12.0}{\rmdefault}{\mddefault}{\updefault}{\color[rgb]{0,0,0}$\sigma_{2,1,0,0}^{2,1,1,2} \in F(2,4;6)$}%
}}}}
\put(3919,-441){\makebox(0,0)[lb]{\smash{{\SetFigFont{10}{12.0}{\rmdefault}{\mddefault}{\updefault}{\color[rgb]{0,0,0}$\sigma_{3,2} \in G(3,6)$}%
}}}}
\put(5070,-441){\makebox(0,0)[lb]{\smash{{\SetFigFont{10}{12.0}{\rmdefault}{\mddefault}{\updefault}{\color[rgb]{0,0,0}$\sigma_{1,0,0,0,0}^{2,2,2,2,1}\in F(1,5;6)$}%
}}}}
\put(1102,711){\makebox(0,0)[lb]{\smash{{\SetFigFont{10}{12.0}{\rmdefault}{\mddefault}{\updefault}{\color[rgb]{0,0,0}degree 1 quantum cycle}%
}}}}
\end{picture}%

\end{center}
\caption{The quantum Mondrian tableaux associated to two Schubert varieties.}
\label{izquant}
\end{figure}

Combining the discussion in \S \ref{secflag} and Proposition \ref{two}
we obtain a Littlewood-Richardson rule for the small quantum
cohomology ring of the Grassmannian $G(k,n)$.  Given a Mondrian
tableau $\sigma_{\lambda}$ in $G(k,n)$ and an integer $d \leq k$, we
can associate a painted Mondrian tableau in $Fl(k-d, k+d; n)$ as
follows: The Mondrian tableau associated to the Schubert variety
$\sigma_{\lambda}$ consists of $k$ nested squares. We take the largest
$k-d$ squares (those of index $d+1, \dots, k$) and color them in
$C_1$. We color the remaining squares in $C_2$. Finally, we add $d$
squares of color $C_2$ at the largest available places in the flag
defining the Mondrian tableau of $\sigma_{\lambda}$ (see Figure
\ref{izquant} for two examples).  We call the resulting painted Mondrian
tableau the quantum Mondrian tableau of degree $d$ associated to
$\sigma_{\lambda}$. This tableau is none other than the painted
Mondrian tableau associated to the special Schubert variety
$X_{\lambda}^{(d)}$ of $Fl(k-d,k+d;n)$. 

Let $\sigma_{\lambda}, \sigma_{\mu}$ and $\sigma_{\nu}$ be three
Schubert cycles in $G(k,n)$ that satisfy the equality $$|\lambda| +
|\mu| + |\nu| = k (n-k) + dn.$$ Apply the algorithm described in the
previous section to the quantum Mondrian tableaux of degree $d$
associated to $\sigma_{\lambda}$ and $\sigma_{\mu}$ to express their
intersection as a sum of Schubert cycles in $Fl(k-d,k+d;n)$. Then
apply the algorithm to the quantum Mondrian tableau of degree $d$
associated to $\sigma_{\nu}$ and each of the summands of the previous
product. We have obtained the
following theorem.
\

\begin{theorem}[\cite{coskun1:LR}, Thm.\ 5.1] \label{gwmon}
The three-pointed Gromov-Witten invariant $I_d (\sigma_{\lambda},
\sigma_{\mu}, \sigma_{\nu})$ is equal to the number of times the point
class occurs as a result of applying the Littlewood-Richardson rule
for the two-step flag 
varieties to the quantum Mondrian tableau  of degree $d$ associated to
$\sigma_{\nu}$ and each outcome of the product of the quantum Mondrian
tableaux of degree $d$ associated to $\sigma_{\lambda}$ and
$\sigma_{\mu}$.
\end{theorem}
 
\begin{figure}[htbp]
\begin{center}
\begin{picture}(0,0)%
\includegraphics{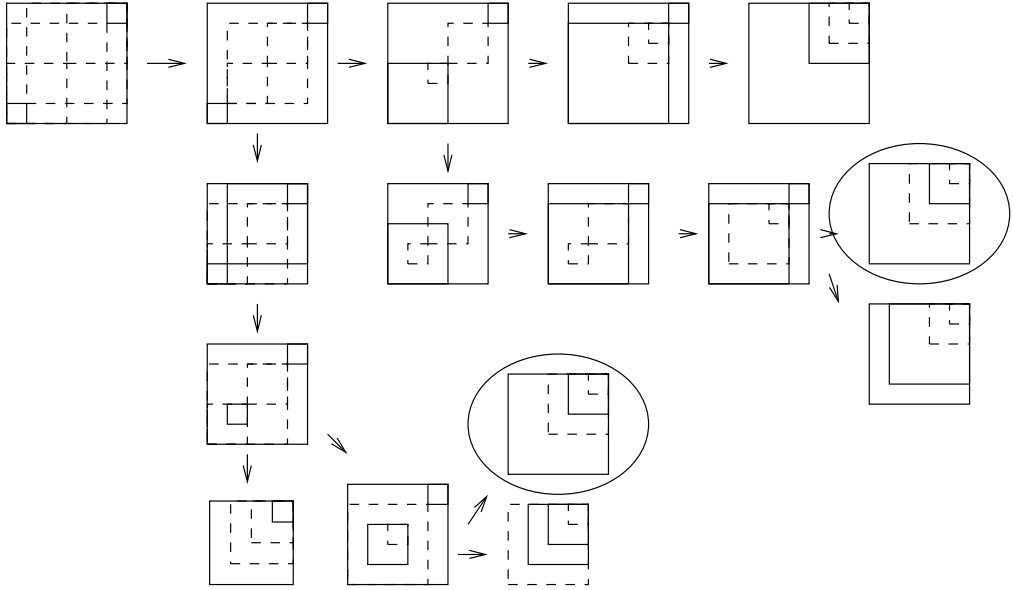}%
\end{picture}%
\setlength{\unitlength}{3947sp}%
\begingroup\makeatletter\ifx\SetFigFont\undefined%
\gdef\SetFigFont#1#2#3#4#5{%
  \reset@font\fontsize{#1}{#2pt}%
  \fontfamily{#3}\fontseries{#4}\fontshape{#5}%
  \selectfont}%
\fi\endgroup%
\begin{picture}(6320,3678)(589,-3523)
\end{picture}%

\end{center}
\caption{Computing the Gromov-Witten invariant $I_{G(3,6),d=1}(\sigma_{3,2,1},
\sigma_{3,2,1}, \sigma_{2,1})=2.$}
\label{quantexamp}
\end{figure}

We illustrate the use of Theorem \ref{gwmon} by computing the
Gromov-Witten invariant $$I_{G(3,6),d=1}(\sigma_{3,2,1},
\sigma_{3,2,1}, \sigma_{2,1})=2.$$ Figure \ref{quantexamp}
demonstrates the computation. The quantum cycle of $d=1$ associated to
$\sigma_{3,2,1}$ (respectively, $\sigma_{2,1}$) is
$\sigma_{2,1,0,0}^{2,1,1,2}$ (respectively,
$\sigma_{1,0,0,0}^{2,1,2,1}$).  In order to calculate the
Gromov-Witten invariant we have to find how many times
$\sigma_{2,2,2,1}^{1,2,1,2}$ (the dual of
$\sigma_{1,0,0,0}^{2,1,2,1}$) occurs in the square of the class
$\sigma_{2,1,0,0}^{2,1,1,2}$. An straightforward calculation with painted
Mondrian tableaux shows that the answer is 2.

%

Earlier, Gromov-Witten invariants could be computed algebraically
based on structure theorems for the quantum ring due to Bertram
\cite{bertram}.  The proofs of these results were drastically
simplified by Buch \cite{buch} using his powerful yet simple
``kernel-span'' technique.  Buch defined the {\em kernel\/} of a
rational curve in the Grassmannian $G(k,n)$ to be the intersection in
$\CC^n$ of all the $k$-planes corresponding to points on the curve,
and defined its {\em span\/} to be the linear span of these
$k$-planes.  When the curve has degree $d$, the dimension of the
kernel is at least $k-d$ while the span has dimension at most $k+d$.
By using this observation, Bertram's structure theorems can be proved
by applying classical Schubert calculus to the span of a curve.

The translation between three-pointed genus-zero Gromov-Witten
invariants and the ordinary products in the two-step flag variety
given in \cite{bkt} is obtained by proving that the map that sends a
rational curve to the pair $(V_1,V_2)$ of its kernel and span is
injective on the set of rational curves contributing to a
Gromov-Witten invariant $I_d(\sigma_\lambda, \sigma_\mu, \sigma_\nu)$,
and the image of this map is exactly the set of points in an
intersection $X^{(d)}_\lambda \cap X^{(d)}_\mu \cap X^{(d)}_\nu$ of
special Schubert varieties.  Furthermore, each rational curve can be
explicitly reconstructed from the pair $(V_1,V_2)$, and this
reconstruction exhibits the curve as a scroll in $\PP^{n-1}$ with
vertex $\PP(V_1)$ and span $\PP(V_2)$.  The problem of counting curves
thus transforms to finding the pairs of vertices and spans.

Unfortunately, only very special curves of degree $d \geq k$ in
 $G(k,n)$ have a non-empty ``kernel'' (the intersection of all the
 linear spaces parametrized by the curve). We can replace the kernel
 by a more natural invariant: the sequence of minimal degree
 subscrolls associated to the curve (see \cite{coskun:jumping}).  The
 advantage of minimal degree subscrolls is that every rational curve
 in the Grassmannian has an associated sequence of minimal degree
 subscrolls. When $d < k$ and the curve is general, we recover the
 vertex of the cone (equivalently, the kernel). 

Considering minimal degree subscrolls allows one to extend the
geometric point of view to the big quantum cohomology of
Grassmannians.  For instance, the geometry of scrolls leads to an
immediate proof of the vanishing of many $m$-pointed Gromov-Witten
invariants (see \cite{coskun:jumping}). Let $m$, $k$, $d$, $n$ be positive
integers satisfying $m \geq 3$, $2k \leq n$ and $d+k \leq n$.  Let
$\sigma_{\lambda^1}$, \dots, $\sigma_{\lambda^m}$ be $m$ Schubert
cycles, where the parts $\lambda_i^j$ of the partitions satisfy
$$\sum_{i,j} \lambda_i^j = dn + k(n-k) + m-3.$$ The Gromov-Witten
invariant $I_{G(k,n),d}(\sigma_{\lambda^1}, \dots, \sigma_{\lambda^m})
= 0$ unless $$\sum_{i,j}\max{(\lambda_i^j - d, 0)} \leq
(k+d)(n-k-d).$$ 

The geometric point of view also leads to a partial
Littlewood-Richardson rule for the big quantum cohomology of
Grassmannians. Using an algorithm similar to the one in
\cite{vakil:rationalelliptic}, one can give a positive rule for
computing some Gromov-Witten invariants of $G(k,n)$ (see
\cite{coskun1:degenerations}, Thm.\ 8.1). It would be interesting to
extend these results to obtain a Littlewood-Richardson rule for
arbitrary genus-zero Gromov-Witten invariants of $G(k,n)$. It would be
equally interesting to extend these results to other homogeneous
varieties, especially flag varieties.

\section{Linear spaces and a quadratic form}\label{sec-orth}
A modification of Mondrian tableaux may be used to calculate classes
of varieties in the isotropic setting as well. A precise description
is beyond the scope of this article. However, in the spirit of the
rest of Part \ref{parttwo}, we give a few illustrative examples. For
precise statements and details, we refer the reader to
\cite{coskun:quadrics}.

\begin{figure}[htbp]
\begin{center}
\begin{picture}(0,0)%
\includegraphics{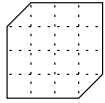}%
\end{picture}%
\setlength{\unitlength}{3947sp}%
\begingroup\makeatletter\ifx\SetFigFont\undefined%
\gdef\SetFigFont#1#2#3#4#5{%
  \reset@font\fontsize{#1}{#2pt}%
  \fontfamily{#3}\fontseries{#4}\fontshape{#5}%
  \selectfont}%
\fi\endgroup%
\begin{picture}(624,624)(514,-73)
\end{picture}%

\end{center}
\caption{A quadram of side-length $4$ depicting a smooth quadric
  surface in $\PP^3$.}
\label{quadram}
\end{figure}

Recall that the orthogonal Grassmannian $OG(k,n)$ can be interpreted
as the Fano variety of $(k-1)$-dimensional projective spaces on a smooth
quadric hypersurface $Q$ in $\PP^{n-1}$. In order to record
specializations of linear spaces on $Q$ we need to denote subquadrics
of $Q$. We will depict a quadric hypersurface in $\PP^r$ as a square
of side-length $r+1$ whose northwest and southeast corners have been
diagonally chopped off by one unit (see Figure \ref{quadram}). We will
refer to such shapes as {\it quadrams} (short for the diagram of a
quadric) and refer to $r+1$ as the side-length of the quadram. We will
also need to denote the singular loci of these quadrics. Consequently,
we will label our quadrams and write the label of the quadram in the
squares denoting the linear spaces along which the quadric is singular.

\begin{figure}[htbp]
\begin{center}
\begin{picture}(0,0)%
\includegraphics{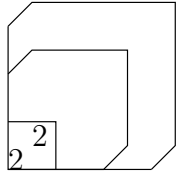}%
\end{picture}%
\setlength{\unitlength}{3947sp}%
\begingroup\makeatletter\ifx\SetFigFont\undefined%
\gdef\SetFigFont#1#2#3#4#5{%
  \reset@font\fontsize{#1}{#2pt}%
  \fontfamily{#3}\fontseries{#4}\fontshape{#5}%
  \selectfont}%
\fi\endgroup%
\begin{picture}(1074,1111)(589,-410)
\put(601,-361){\makebox(0,0)[lb]{\smash{{\SetFigFont{12}{14.4}{\rmdefault}{\mddefault}{\updefault}{\color[rgb]{0,0,0}$2$}%
}}}}
\put(751,-211){\makebox(0,0)[lb]{\smash{{\SetFigFont{12}{14.4}{\rmdefault}{\mddefault}{\updefault}{\color[rgb]{0,0,0}$2$}%
}}}}
\end{picture}%

\end{center}
\caption{The quadric diagram associated to $\sigma_{2}^{2,0}$ in $OG(3,7)$.}
\label{quadramexam}
\end{figure}

Suppose $n$ is odd. Set $m = (n-1)/2$.  Let $\lambda = \lambda_1 >
\cdots > \lambda_s$ and $\mu = \mu_{s+1} > \dots > \mu_k$ be two
strictly decreasing partitions as in \S \ref{othergps}. Let
$\sigma_{\lambda, \mu}$ be the corresponding Schubert cycle in
$OG(k,n)$. The {\it quadric diagram} associated to $\sigma_{\lambda,
  \mu}$ is a collection of $s$ nested squares $S_1 \subset \dots
\subset S_s$ of side-lengths $m+1-\lambda_1, \dots, m+1 -\lambda_s$,
$k-s$ nested quadrams $Q_{k-s} \subset \cdots \subset Q_1$ of
side-lengths $n-\mu_{s+1} , \dots, n-\mu_k$ containing all the squares
and a function $f$. The function $f$ associates to the unit squares in
the anti-diagonal of the quadric diagram the set of labels of the
quadrams which represent quadrics that are singular at that point.
When $n$ is even, we make the obvious modifications to this
definition. In particular, we have to distinguish between the
half-dimensional linear spaces that belong to the different
irreducible components. We will depict those that correspond to
$\lambda_s = 0$ in solid lines and those that correspond to $\mu_{s+1}
= m-1$ in dashed lines. Figure \ref{quadramexam} shows the quadric
diagram associated to $\sigma_{2}^{2,0}$ in $OG(3,7)$.

Geometrically, the $\PP^2$s parametrized by the Poincar\'e dual of
$\sigma_{2}^{2,0}$ intersect the line $l$ represented by the square
$S_1$. They intersect the subquadric everywhere tangent along $l$ in a
line. Note that this subquadric, denoted by $Q_2$, is singular along
$l$. Hence, we place the index ``$2$'' in the unit squares contained in
$S_1$. Finally, the $\PP^2$s are contained in the quadric represented
by the largest quadram $Q_1$. 

One can use quadric diagrams to calculate the classes of subvarieties
in orthogonal Grassmannians. For instance to multiply two Schubert
varieties in $OG(k,n)$ where the index $\mu$ does not contribute to
the discrepancy, we place the two Schubert cycles in opposite corners
of an $n\times n$ square. After initial manipulations, the quadric
diagram consists of squares that form a generalized Mondrian tableau
and nested quadrams containing the squares. Informally, we nest the
squares as in the Mondrian tableaux rule for ordinary
Grassmannians keeping track of the singularities of the quadrics
represented by the quadrams. Eventually all the squares and quadrams
are nested. However, this might still not correspond to a Schubert
variety. We further degenerate these varieties until the quadrics
represented by each quadram is as singular as possible. We give a few
examples demonstrating how the procedure works in practice, referring
the reader to \cite{coskun:quadrics} for details.

\begin{figure}[htbp]
\begin{center}
\begin{picture}(0,0)%
\includegraphics{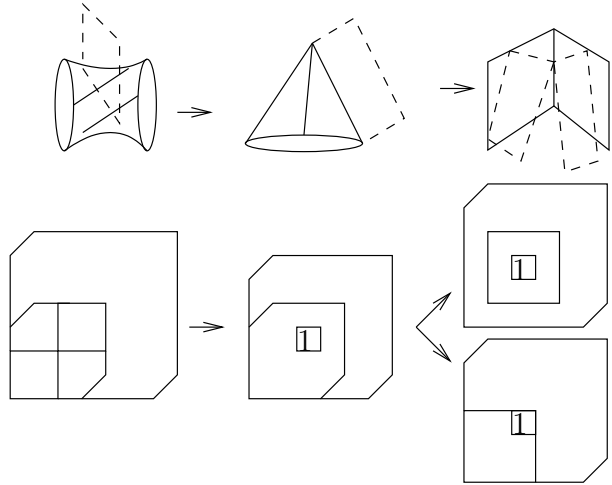}%
\end{picture}%
\setlength{\unitlength}{3947sp}%
\begingroup\makeatletter\ifx\SetFigFont\undefined%
\gdef\SetFigFont#1#2#3#4#5{%
  \reset@font\fontsize{#1}{#2pt}%
  \fontfamily{#3}\fontseries{#4}\fontshape{#5}%
  \selectfont}%
\fi\endgroup%
\begin{picture}(3785,3029)(589,-2548)
\put(2401,-1711){\makebox(0,0)[lb]{\smash{{\SetFigFont{12}{14.4}{\rmdefault}{\mddefault}{\updefault}{\color[rgb]{0,0,0}$1$}%
}}}}
\put(3751,-1261){\makebox(0,0)[lb]{\smash{{\SetFigFont{12}{14.4}{\rmdefault}{\mddefault}{\updefault}{\color[rgb]{0,0,0}$1$}%
}}}}
\put(3751,-2236){\makebox(0,0)[lb]{\smash{{\SetFigFont{12}{14.4}{\rmdefault}{\mddefault}{\updefault}{\color[rgb]{0,0,0}$1$}%
}}}}
\end{picture}%

\end{center}
\caption{Calculating the intersection $\sigma_2^{2,0} \cdot
  \sigma_2^{2,0} = 2 \sigma_{3,1}^1$ in $OG(3,7)$.}
\label{orth1}
\end{figure}

\begin{example}
  We calculate $\sigma_2^{2,0} \cdot \sigma_2^{2,0}$ in $OG(3,7)$.
  Figure \ref{orth1} shows the quadric diagrams and the corresponding
  projective geometry associated to this example.  Geometrically, this
  corresponds to finding the class of the variety of projective planes
  on a smooth quadric $Q$ in $\PP^6$ that intersect $2$ general lines
  $l_1$ and $l_2$. These lines span a $\PP^3$ that intersects $Q$ in a
  smooth quadric surface (depicted by $Q_2$ in the quadric diagram in
  the first panel). The planes must contain a line in the ruling
  opposite to $l_1$ and $l_2$ on this quadric surface. In order to
  calculate the class we specialize the conditions on the planes.  We
  move the line $l_1$ along the ruling of a quadric surface contained
  in $Q$ that intersects $l_2$.  Once we specialize $l_1$ to intersect
  $l_2$, the limit of the $\PP^3$s spanned by the two lines becomes
  tangent to the quadric at the point of intersection of $l_1$ and
  $l_2$. The intersection of the limit $\PP^3$ with the quadric $Q$ is
  a quadric cone (denoted by $Q_2$ in the second panel). The limit of
  the planes is the locus of planes that contain a line in this
  quadric cone (necessarily containing the cone point). Note that
  since the planes contain the cone point $p$, they have to be
  contained in the tangent space to $Q$ at $p$ (hence, we shrink $Q_1$
  by one unit in the second panel).  This is not a Schubert variety.
  To degenerate it into a Schubert variety, we take a pencil of
  $\PP^3$s that become more tangent to the quadric $Q$. In the limit
  the quadric cones break into a union of two planes. The limit of the
  lines is the locus of lines on each of the planes that pass through
  the limit of the vertices of the quadric cones. We conclude that
  $\sigma_2^{2,0} \cdot \sigma_2^{2,0} = 2 \sigma_{3,1}^1$.  It is not
  hard to see that both limits occur with multiplicity
  $1$.  

\end{example}

\begin{example} 
  The orthogonal Grassmannian $OG(2,6)$ and the two-step flag variety
  $Fl(1,3;4)$ are isomorphic. It is both instructive and fun to
  calculate the cohomology ring of $OG(2,6)$ using both quadric
  diagrams and the Mondrian tableaux rule for two-step flag varieties.
  Figure \ref{orthmon} shows two sample calculations. In the left
  panel, we calculate the class of lines in $OG(2,6)$ that intersect
  two general lines $l_1$, $l_2$. These are lines that are contained
  in the intersection of $Q$ with the $\PP^3$ spanned by the two
  lines. The degeneration is similar to the previous example. However,
  when the quadric cone breaks in this case, the two planes belong to
  two different irreducible components.  In the right panel, we
  calculate the intersection of Schubert varieties parameterizing
  lines that intersect a plane (one in each of the two rulings) and
  are contained in the tangent space at a point in the plane. The
  calculation is similar. In each panel we also show the corresponding
  Mondrian tableaux calculation below the quadric diagrams.

\begin{figure}[htbp]
\begin{center}
\begin{picture}(0,0)%
\includegraphics{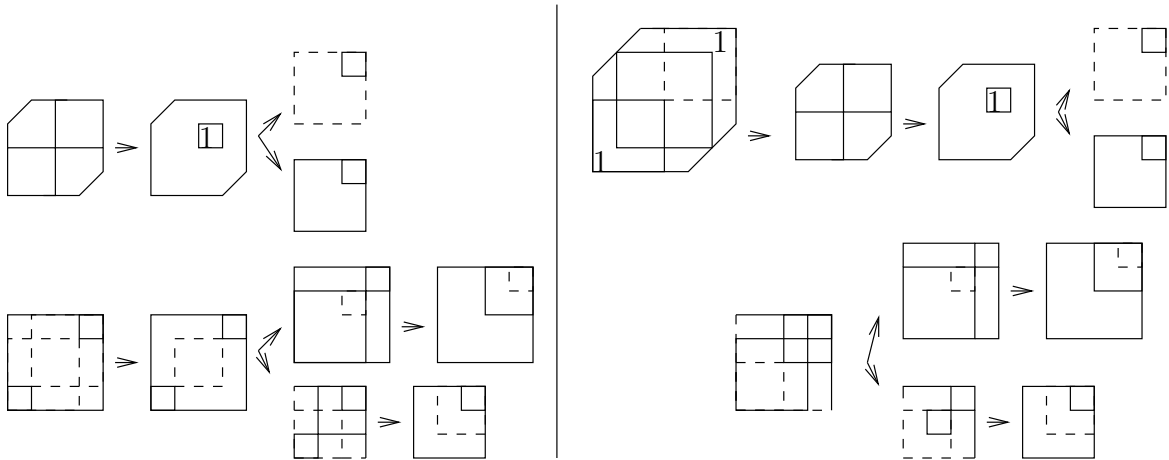}%
\end{picture}%
\setlength{\unitlength}{3947sp}%
\begingroup\makeatletter\ifx\SetFigFont\undefined%
\gdef\SetFigFont#1#2#3#4#5{%
  \reset@font\fontsize{#1}{#2pt}%
  \fontfamily{#3}\fontseries{#4}\fontshape{#5}%
  \selectfont}%
\fi\endgroup%
\begin{picture}(7299,2874)(439,-2398)
\put(1651,-436){\makebox(0,0)[lb]{\smash{{\SetFigFont{12}{14.4}{\rmdefault}{\mddefault}{\updefault}{\color[rgb]{0,0,0}$1$}%
}}}}
\put(6601,-211){\makebox(0,0)[lb]{\smash{{\SetFigFont{12}{14.4}{\rmdefault}{\mddefault}{\updefault}{\color[rgb]{0,0,0}$1$}%
}}}}
\put(4126,-586){\makebox(0,0)[lb]{\smash{{\SetFigFont{12}{14.4}{\rmdefault}{\mddefault}{\updefault}{\color[rgb]{0,0,0}$1$}%
}}}}
\put(4876,164){\makebox(0,0)[lb]{\smash{{\SetFigFont{12}{14.4}{\rmdefault}{\mddefault}{\updefault}{\color[rgb]{0,0,0}$1$}%
}}}}
\end{picture}%

\end{center}
\caption{Calculating the intersection $\sigma_{1}^0 \cdot \sigma_{1}^0
  = \sigma_{2,0} + \sigma_{2}^2$ in $OG(2,6)$, equivalently
  $\sigma_{1,0,0}^{2,1,2} \cdot
  \sigma_{1,0,0}^{2,1,2}=\sigma_{1,1,0}^{1,2,2} +
  \sigma_{1,1,1}^{2,1,2}$ in $Fl(1,3;4)$ (the left panel) and the
  intersection $\sigma_0^1 \cdot \sigma^{2,1}= \sigma_{2,0} +
  \sigma_{2}^2$ in $OG(2,6)$, equivalently $\sigma_{0,0,0}^{1,2,2}
  \cdot \sigma_{1,1,0}^{2,2,1} = \sigma_{1,1,0}^{1,2,2} +
  \sigma_{1,1,1}^{2,1,2} $ in $Fl(1,3;4)$ (the right panel).}
\label{orthmon}
\end{figure}
\end{example}

\end{document}